\documentclass[10pt]{amsart}
\usepackage{amssymb,amsmath,amsthm,bbm,amsfonts,amsbsy}
\usepackage[all]{xy}
\usepackage{tikz-cd}
\usepackage[shortlabels]{enumitem}
\numberwithin{equation}{section}
\newtheorem{thm}{Theorem}[section]

\newtheorem{defn}[thm]{Definition}

\newtheorem*{kerr}{Kerr Theorem}
\newtheorem*{TCM}{Tanaka-Chern-Moser Classification}

\theoremstyle{definition}
\newtheorem{definition}[thm]{Definition}
\newtheorem{rem}[thm]{Remark}

\theoremstyle{remark}

\allowdisplaybreaks[1]

\newcommand{\im}{\mathbf{i}}
\newcommand{\dd}{\text{d}}
\newcommand{\F}{\mathcal{F}}

\newcommand{\e}{\mathbf{e}}
\newcommand{\Q}{\mathcal{Q}}
\newcommand{\N}{\mathcal{N}}
\newcommand{\bb}{{\tt b}}
\newcommand{\h}{{\tt h}}
\newcommand{\vv}{{\tt v}}
\newcommand{\n}{{\tt n}}
\newcommand{\II}{\mathbf{II}}

\begin{document}

\title{3-folds CR-embedded in 5-dimensional real hyperquadrics}

\author{Curtis Porter}
\address{Katedra Matematiky, Univerzita Hradec Kr\'{a}lov\'{e}}
\email{cwp1729@gmail.com}


\begin{abstract}
E. Cartan's method of moving frames is applied to 3-dimensional manifolds $M$ which are CR-embedded in 5-dimensional real hyperquadrics $\Q$ in order to classify $M$ up to CR symmetries of $\Q$ given by the action of one of the Lie groups $SU(3,1)$ or $SU(2,2)$. In the latter case, the CR structure of $M$ derives from a shear-free null geodesic congruence on Minkowski spacetime, and the relationship to relativity is discussed. In both cases, we compute which homogeneous CR 3-folds appear in $\Q$.  
\end{abstract}

\maketitle
\tableofcontents

\section{Introduction}

For several physically significant solutions to Einstein's field equations in general relativity, spacetime is a 4-dimensional Lorentzian manifold that is foliated by a family of curves called a shear-free null geodesic congruence (SFNGC), which induces a CR structure on the 3-dimensional leaf space of the foliation. Conversely, a 3-dimensional CR structure can be ``lifted" to a spacetime admitting a SFNGC. The Robinson-Trautman metrics, for example, describe congruences which are hypersurface-orthogonal, and their corresponding CR structures are Levi-flat. These include models of electromagnetic and gravitational fields radiating along the foliating curves, generating wave-fronts orthogonal to their direction of propagation. Levi-nondegenerate CR structures, on the other hand, are associated with ``twisting" congruences, such as those appearing in Kerr's model of a rotating black hole. The geometry of SFNGCs is explained in \S\ref{SFNGCsec}, and a glimpse of their history in relativity is given in \S\ref{history}.

CR maps between manifolds establish a notion of \emph{intrinsic} CR equivalence, in terms of which Levi-flat CR structures are all locally the same in any fixed dimension. By contrast, Levi-nondegenerate CR manifolds $M$ are locally classified by Cartan's method of equivalence, which constructs a principal \emph{Cartan bundle} $\mathcal{B}\to M$ of (co)frames adapted to the CR structure of $M$ as well as a canonical \emph{Cartan connection} $\gamma$ on $\mathcal{B}$. When the curvature tensor $\dd\gamma+\tfrac{1}{2}[\gamma,\gamma]$ vanishes, $M$ is locally equivalent to a real hyperquadric $\Q$, and $\mathcal{B}$ is a special unitary Lie group of CR symmetries of $\Q$, with $\gamma$ playing the role of the Maurer-Cartan form on $\mathcal{B}$. See \S\ref{CRdefsec} for basic definitions and references. In dimension three, the real hyperquadric is the CR 3-sphere whose Lie group of CR symmetries is $SU(2,1)$. If $M$ is homogeneous under the action of its Lie group $\mathcal{G}$ of intrinsic CR symmetries, then either $\mathcal{G}=SU(2,1)$ and $M$ is necessarily flat, or else $\mathcal{G}$ is 3-dimensional and the curvature of $M$ may be zero, nonzero, or undefined (if $M$ is Levi-flat). In any case, $M$ is recognizable by the Lie algebra $\mathfrak{g}$ of $\mathcal{G}$ since 3-dimensional Lie algebras over $\mathbb{R}$ were classified by Bianchi. Homogeneous models of Levi-nondegenerate 3-folds are catalogued at the end of \S\ref{CR3sec}.

The properties of a SFNGC are maintained under conformal rescaling of the spacetime metric; indeed, a given SFNGC belongs to a substantially larger family of metrics than a conformal class, and the ambiguity in the choice of a representative is related to that of a choice of adapted coframing on the underlying CR manifold. It is therefore natural to ask which CR structures give rise to SFNGCs whose family of metrics contains one with special characteristics; e.g., a metric which is conformally flat. The answer to the latter question -- provided by a theorem attributed to Kerr -- is those 3-folds that may be embedded by a CR map into the 5-dimensional real hyperquadric $\Q\subset\mathbb{CP}^3$ whose CR symmetry group $SU(2,2)$ is infinitesimally isomorphic to the conformal symmetry group $SO(2,4)$ of compactified Minkowski spacetime. The Kerr Theorem is presented in \S\ref{kerr} in the framework of Penrose's twistor theory, emphasizing the role of these symmetry groups. Unsurprisingly, homogeneous 3-folds are noteworthy in this context, both mathematically and physically. Their symmetry simplifies their structure equations, making them ideal subjects to test any criteria for embeddability in $\Q$. Those that are embedded with symmetries in $SU(2,2)$ correspond to SFNGC with conformal symmetries, while those that admit no embedding at all correspond to SFNGC in spacetimes with curvature.  The proof of the Kerr Theorem is sketched in \S\ref{kerrproof}. 

There is little added expense in generalizing beyond the context of the Kerr Theorem. To wit, there are two 5-dimensional real hyperquadrics, differing in the signature of their (rank-2) Levi forms, which is either definite or split. They are both called $\Q$; when necessary we recognize them individually by their CR symmetry groups: $SU(3,1)$ for the CR 5-sphere and $SU(2,2)$ for the projectivized null twistors, respectively, or $SU_\star$ to refer to both. They are both fixed by a parabolic subgroup $\mathcal{P}\subset SU_\star$ stabilizing a null line in $\mathbb{C}^4$. We consider CR embeddings of 3-folds $M$ into both homogeneous spaces $\Q=SU_\star/\mathcal{P}$. Of especial interest are embeddings that preserve intrinsic CR symmetries of $M$: 
\begin{defn}\label{equivembed}
Let $\Q$ be the 5-dimensional real hyperquadric with CR symmetry group $SU_\star$, where $SU_\star$ is one of $SU(2,2)$ or $SU(3,1)$, and let $M$ be a homogeneous CR 3-fold whose Lie group of CR symmetries is $\mathcal{G}$. A CR embedding $f:M\to\Q$ is \emph{equivariant} if for every $x\in M$ and $\varsigma\in \mathcal{G}$, there exists $\varsigma_\star\in SU_\star$ such that $f(\varsigma(x))=\varsigma_\star(f(x))$. If such an embedding exists, $M$ is \emph{equivariantly embeddable} in $\Q$.
\end{defn}

Curry and Gover (\cite{CGbook}) developed a general framework for analyzing CR embeddings using tractor calculus, which replaces the Cartan bundle $\mathcal{B}$ and Cartan connection $\gamma$ of $M$ with an associated vector bundle and differential operator, though the distinction is purely aesthetic since $\mathcal{B}$ and $\gamma$ are recoverable from the tractor bundle and connection. As such, when the ambient CR manifold under consideration is a hyperquadric $\Q$, employing the Curry-Gover formalism to study embedded CR submanifolds $M\subset \Q$ should be equivalent to applying Cartan's method of moving frames (\cite{GriffithsMF}, \cite{CFB}), which is the strategy of the present paper. As an application of exterior differential systems, the method of moving frames is similar in technique to the method of equivalence, except that it involves adapting (co)frames of an ambient homogeneous space to the geometry of embedded submanifolds. Accordingly, moving frames classify submanifolds up to \emph{extrinsic} equivalence determined by the symmetry group of the homogeneous space. This extrinsic action may discriminate CR submanifolds which are intrinsically equivalent. In particular, Levi-flat manifolds are not necessarily locally identical, and the local classification of Levi-nondegenerate CR 3-folds is more refined.

In \S\ref{hermframesec}, the Cartan bundle of $\Q$ is labeled $\mathcal{H}$, in part to distinguish it from that of an abstract, Levi-nondegenerate 3-fold $M$, but also because it is constructed as \emph{Hermitian frames} or bases of $\mathbb{C}^4$. This approach makes explicit the isomorphism $\mathcal{H}\cong SU_\star$. Its Cartan connection -- the Maurer-Cartan form of $SU_\star$ -- is denoted $\mu$ so there is no confusion with $\gamma$, the Cartan connection of $M$. In \S\ref{CR3sec}, the method of equivalence is implemented to build $\mathcal{B}\to M$ and $\gamma$ by the standard procedure that incorporates all intrinsic symmetries of $M$. We begin adapting frames in $\mathcal{H}$ over $M\subset\Q$ in \S\ref{firstadaptsec}. Informally speaking, this process divides $\mu$ into two pieces: one that contains $\gamma$ as though $\mathcal{H}$ contains $\mathcal{B}$, and another that acts on the ``normal bundle," encoded in a bilinear form of rank at most two over $M$ that we dub the \emph{second fundamental form} and denote $\II$ in \S\ref{IIsec} by analogy to the study of Riemannian hypersurfaces in Euclidean space. 

The classification of Levi-nondegenerate $M\subset\Q$ in \S\ref{LNclasssec} parallels their intrinsic classification in \S\ref{CR3sec}, including the distinction between flat 3-folds (\S\ref{LNflatsec}) and those $M$ with nonzero curvature $\dd\gamma+\tfrac{1}{2}[\gamma,\gamma]$ (\S\ref{LNcurvedsec}). An essential result of the paper is Theorem \ref{curvedtheorem}, which details how the second fundamental form of an embedded, Levi-nondegenerate 3-fold relates to its intrinsic CR invariants in the case that the curvature tensor is nonvanishing. This characterizes up to symmetries in $SU_\star$ the abstract, curved 3-folds which admit embeddings in $\Q$ in the spirit of the Bonnet Theorem. The praxis of Theorem \ref{curvedtheorem} is complicated in general, but it simplifies when $M$ is homogeneous, and even more so for equivariant embeddings as in Definition \ref{equivembed}:
 
\begin{thm}\label{curvedhomogthm}
Let $\Q$ be the 5-dimensional real hyperquadric with CR symmetry group $SU_\star$, where $SU_\star$ is one of $SU(2,2)$ or $SU(3,1)$, and let $M$ be a homogeneous, Levi-nondegenerate CR 3-fold with nonzero curvature whose Lie algebra of infinitesimal CR symmetries is $\mathfrak{g}$. $M$ is locally equivariantly embeddable in $\Q$ if and only if one of the following is true:
\begin{enumerate}
\item $SU_\star=SU(3,1)$ and $\mathfrak{g}=\mathfrak{su}(2)$; every member of the parameter family of these models is realized in $\Q$ with rank$(\II)=2$.

\item $SU_\star=SU(2,2)$ and $\mathfrak{g}=\mathfrak{su}(1,1)$; every member of the parameter family of these models is realized in $\Q$ with rank$(\II)=2$.

\item $SU_\star=SU(2,2)$ and $\mathfrak{g}=\mathfrak{su}(2)$; a single member of the parameter family of these models is realized in $\Q$ with rank$(\II)=1$.
\end{enumerate}
\end{thm}

Since Levi-nondegenerate 3-folds with vanishing curvature are all locally intrinsically equivalent to the 3-sphere, there is no question that they are locally embeddable in either hyperquadric. Still, the method of moving frames serves to distinguish embeddings that are inequivalent under the action of $SU_\star$ on $\Q$, and as always the rank of $\II$ is a valuable invariant of this action. Those $M\subset\Q$ with $\II$ of submaximal rank are all equivariantly embeddable. The maximally symmetric case $\II=0$ realizes the full symmetry group of the 3-sphere, whereas rank$(\II)=1$ describes models with submaximal symmetry. By contrast, embeddings with rank$(\II)=2$ are classified by structure equations \eqref{flatgenSE} which only evince an equivariant embedding in the hyperquadric with symmetry group $SU(2,2)$.
  
\begin{thm}\label{flathomogthm}
Let $\Q$ be the 5-dimensional real hyperquadric with CR symmetry group $SU_\star$, where $SU_\star$ is one of $SU(2,2)$ or $SU(3,1)$, and let $M$ be a homogeneous, Levi-nondegenerate CR 3-fold with zero curvature whose Lie group of CR symmetries is $\mathcal{G}$ with Lie algebra $\mathfrak{g}$. $M$ is locally equivariantly embeddable in $\Q$ if and only if one of the following is true:
\begin{enumerate}
\item $\mathcal{G}=SU(2,1)$ and $M$ is the CR 3-sphere. $M$ is embedded as an orbit of $U(2,1)\subset SU_\star$ with $\II=0$.

\item $SU_\star=SU(3,1)$ and $M$ is the unique flat model with $\mathfrak{g}=\mathfrak{su}(2)$. $M$ is locally embedded as an orbit of $\hat{\mathcal{G}}\subset SU(3,1)$ with rank$(\II)=1$, where $\hat{\mathcal{G}}$ is an extension of $\mathcal{G}$ by a central action of $U(1)\subset SU(3,1)$.

\item $SU_\star=SU(2,2)$ and $M$ is the unique flat model with $\mathfrak{g}=\mathfrak{su}(1,1)$. $M$ is locally embedded as an orbit of $\hat{\mathcal{G}}\subset SU(2,2)$ with rank$(\II)=1$, where $\hat{\mathcal{G}}$ is an extension of $\mathcal{G}$ by a central action of $U(1)\subset SU(2,2)$.

\item $SU_\star=SU(2,2)$ and $M$ is the unique flat model with $\mathfrak{g}$ as the extension of $\mathbb{R}^2$ by the derivation $\left[\begin{smallmatrix}3&1\\1&3\end{smallmatrix}\right]$. $M$ is embedded with rank$(\II)=2$.
\end{enumerate}
\end{thm}

For $\Q$ whose Levi form has split signature, it is also possible that $M\subset\Q$ is Levi-flat, in which case there is no Cartan bundle or connection to speak of. All such $M$ are locally the same, intrinsically. Even so, the reduction of $\mathcal{H}\cong SU(2,2)\to M$ in \S\ref{firstadaptsec} and the definition of the second fundamental form $\II$ in \S\ref{LFIIsec} are meaningful, and in \S\ref{LFclasssec} we classify Levi-flat $M\subset\Q$ up to the action of $SU(2,2)$ on $\Q$. The maximally symmetric case $\II=0$ has the structure equations of the 10-dimensional parabolic subgroup $\mathcal{R}\subset SU(2,2)$ that stabilizes a partial flag given by a null line in a 3-plane spanned by null lines. When rank$(\II)=1$, $\mu$ is fully reduced to a coframing of $M$ with classifying structure equations \eqref{LFrankone}, including one homogeneous model. Finally, the maximal rank case branches based on first-order behavior of $\II$, where the generic subcase is fully reduced and classified by \eqref{LFr2u1}, and the alternative exhibits the Maurer-Cartan equations of $\mathfrak{sl}_2\mathbb{R}\oplus\mathfrak{su}(p,q)$ for either $(p,q)=(2,0)$ or $(p,q)=(1,1)$.

\begin{thm}\label{LFhomogthm}
Let $\Q$ be the 5-dimensional real hyperquadric with CR symmetry group $SU(2,2)$, and let $M$ be a Levi-flat CR 3-fold embedded in $\Q$ with second fundamental form $\II$. $M$ is homogeneous for its Lie group of CR symmetries induced by the action of $SU(2,2)$ on $\Q$ if and only if one of the following is true: 
\begin{enumerate}
\item $\II=0$ and $M$ contains a complex line which is null for the Levi form of $\Q$. In this case the group if CR symmetries of $M$ is the 10-dimensional parabolic subgroup $\mathcal{R}\subset SU(2,2)$ that stabilizes a partial flag given by a null line in a 3-plane spanned by null lines.

\item Rank$(\II)=1$ and the algebra of infinitesimal CR symmetries of $M$ is isomorphic to the extension of $\mathbb{R}^2$ by the derivation $\sqrt[3]{3}\left[\begin{smallmatrix}1&0\\0&-1\end{smallmatrix}\right]$.

\item Rank$(\II)=2$ and the algebra of infinitesimal CR symmetries of $M$ is $\mathfrak{sl}_2\mathbb{R}\oplus\mathfrak{su}(p,q)$ for either $(p,q)=(2,0)$ or $(p,q)=(1,1)$.

\item Rank$(\II)=2$ and the algebra of infinitesimal CR symmetries of $M$ is isomorphic to the extension of $\mathbb{R}^2$ by the derivation $\left[\begin{smallmatrix}1&0\\0&2\end{smallmatrix}\right]$.
\end{enumerate}
\end{thm}

\vspace{\baselineskip}

\noindent{\bf Acknowledgements:} This project was initiated while the author participated in the Fall 2017 Simons Semester \textit{Symmetry and Geometric Structures} hosted by the Mathematics Institute of the Polish Academy of Sciences. The author also acknowledges the Czech Science Foundation (GA\v{C}R) for support via the program GA\v{C}R 19-14466Y.

\vspace{\baselineskip}

\section{CR Structures}\label{CRsec}

\subsection{Basic Definitions, Adapted Coframings}\label{CRdefsec}

For any fiber bundle $\pi:E\to M$, $E_x=\pi^{-1}(x)$ denotes the fiber of $E$ over $x\in M$ and $\Gamma(E)$ denotes the sheaf of smooth (local) sections of $E$. If $E$ is a vector bundle, $\mathbb{C}E$ is its complexification whose fibers are $\mathbb{C}E_x=E_x\otimes_\mathbb{R}\mathbb{C}$. We use bold text for the constants $\im=\sqrt{-1}$ and $\mathbf{e}$, the natural exponential.

Here, \emph{CR structure} refers specifically to a \emph{hypersurface-type} CR structure $(M,D,J)$, which is a $(2n+1)$-dimensional smooth manifold $M$ equipped with a corank-$1$ distribution $D\subset TM$ carrying an \emph{almost-complex structure}
\begin{align*}
&J:D\to D,
&J^2=-\mathbbm{1},
\end{align*}
where $J_x:D_x\to D_x$ is linear for every $x\in M$, and $\mathbbm{1}$ is the identity map on the fibers of $D$. The induced action of $J$ on the complexified bundle splits 
\begin{align*}
\mathbb{C}D=H\oplus\overline{H},
\end{align*}
where the \emph{CR bundle} $H$ is the $\im$-eigenspace and the \emph{anti-CR bundle} $\overline{H}$ the $(-\im)$-eigenspace of $J$. The \emph{CR dimension} of $M$ is $\text{rank}_\mathbb{C}H=n$. 

Given two CR structures $(M_1,D_1,J_1)$ and $(M_2,D_2,J_2)$, a \emph{CR map} between them is a smooth map $f:M_1\to M_2$ whose pushforward $f_*:TM_1\to TM_2$ satisfies $f_*D_1\subset D_2$ and $f_*\circ J_1=J_2\circ f_*$; in other words, $f_*H_1\subset H_2$. $M_1$ and $M_2$ are \emph{CR equivalent} if there exists a CR map between them that is a diffeomorphism. Often we are merely concerned with \emph{local} equivalence maps, defined on some neighborhood of any given point:

\begin{defn}\label{CRembed} Let $M_1,M_2$ be CR manifolds of CR dimension $n_1,n_2$, respectively, where $n_1\leq n_2$. $M_1$ is said to be (locally) \emph{CR-embeddable} in $M_2$ if for each $x\in M_1$, there is a neighborhood $N_x\subset M_1$ and a CR map $f:N_x\to M_2$ which is a CR equivalence onto its image; such $f$ is a \emph{CR embedding}. If $n_2=n_1$, a CR embedding is a (local) CR equivalence, and if $M_2=M_1$ a local CR equivalence is a \emph{(local) CR symmetry} of $M_1$.
\end{defn}

All CR structures in this paper are \emph{CR-integrable}; i.e., sections of the (anti-)CR bundle are closed under the Lie bracket of vector fields,
\begin{align}\label{CRint}
[\Gamma(H),\Gamma(H)]\subset\Gamma(H)\Longleftrightarrow
[\Gamma(\overline{H}),\Gamma(\overline{H})]\subset\Gamma(\overline{H}).
\end{align}
The failure of integrability of the underlying real distribution $D$ is measured by the Levi form, 
\begin{align}\label{levi}
&\left.\begin{array}{rl}\ell:H_x\times H_x&\to\mathbb{C}T_xM/\mathbb{C}D_x\\
(y_1,y_2)&\mapsto\im[Y_1,\overline{Y}_2](x)\mod \mathbb{C}D_x\end{array}\right\}
&&Y_i\in\Gamma(H),
&&Y_i(x)=y_i\quad(i=1,2).
\end{align}
$M$ is \emph{Levi-flat} if $\ell$ vanishes identically. The Newlander-Nirenberg Theorem implies that Levi-flat CR manifolds are locally CR equivalent to $\mathbb{R}\times\mathbb{C}^n$.
 
CR structures can locally be encoded into an adapted coframing. Writing $D^\bot\subset T^*M,\overline{H}^\bot\subset\mathbb{C}T^*M$ for the annihilators of $D$ and $\overline{H}$, a \emph{0-adapted coframing} is given by a collection of 1-forms
\begin{align}\label{zeroadapted}
&\varphi^0\in\Gamma(D^\bot),
&\varphi^j\in\Gamma(\overline{H}^\bot)\quad 1\leq j\leq n,
&&\text{such that}
&&\varphi^0\wedge\Big(\bigwedge_{j=1}^n\varphi^j\Big)\wedge\Big(\bigwedge_{j=1}^n\overline{\varphi}^j\Big)\neq 0.
\end{align}
Equivalently, a 0-adapted coframing is a local section of the bundle $\pi:\F^0\to M$ whose fiber over $x\in M$ consists of \emph{0-adapted coframes}, which are linear isomorphisms
\begin{align}\label{zeroadaptedbundle}
\F^0_x=\{\varphi_x:T_xM\stackrel{\simeq}{\longrightarrow} \mathbb{R}\oplus\mathbb{C}^n\ |\ \varphi_x(D_x)=\mathbb{C}^n,\ \varphi_x\circ J_x=\im\varphi_x\}.
\end{align}
We call $\varphi^0$ a \emph{characteristic form}, while $\varphi^j$ and $\overline{\varphi}^j$ are CR and anti-CR forms, respectively. The CR integrability condition \eqref{CRint} is expressed
\begin{align}\label{CRintforms}
&\dd\varphi^i\equiv0\mod\{\varphi^0,\dots,\varphi^n\},
&0\leq i\leq n.
\end{align}
In particular, a characteristic form is real-valued, so using the summation convention we can write
\begin{align}\label{levimatrix}
&\dd\varphi^0\equiv\im\ell_{jk}\varphi^j\wedge\overline{\varphi}^k\mod\{\varphi^0\},
&\ell_{kj}=\overline{\ell}_{jk}\in C^\infty(M,\mathbb{C});
&&1\leq j,k\leq n,
\end{align}
where $\ell_{jk}(x)$ is a local representation of the Levi form \eqref{levi} as an $n\times n$ Hermitian-symmetric matrix. The signature $(p,q)$ of this matrix is an invariant of $M$ under CR equivalence (modulo $(p,q)\sim(q,p)$).

Of course, the 0-adapted coframing $\{\varphi^0,\varphi^j\}$ is not uniquely determined by \eqref{zeroadapted}, but only up to a transformation of the form 
\begin{align*}
&\left[\begin{array}{cc}u&0\\b&a\end{array}\right]\left[\begin{array}{c}\varphi^0\\\varphi^j\end{array}\right],
&0\neq u\in C^\infty(M),\quad a\in C^\infty(M,GL_n\mathbb{C}),\quad b\in C^\infty(M,\mathbb{C}^n).
\end{align*}
Equivalently, the bundle $\pi:\F^0\to M$ carries a natural $G_0$-principal action on its fibers \eqref{zeroadaptedbundle},
\begin{align}\label{Gzero}
G_0=\left.\left\{\left[\begin{array}{cc}u&0\\b&a\end{array}\right]\in GL(\mathbb{R}\oplus\mathbb{C}^n)\ \right|\ 
0\neq u\in\mathbb{R},\ a\in GL_n\mathbb{C},\ b\in\mathbb{C}^n \right\}.
\end{align} 
Thus we see that a CR structure is an example of a G-structure (\cite[Def.2.1]{BGG}); there is a tautologically defined 1-form $\Phi\in\Omega^1(\F^0,\mathbb{R}\oplus\mathbb{C}^n)$,
\begin{align}\label{tautform}
\Phi|_{\varphi_x}=\varphi_x\circ\pi_*,
\end{align}
and any local equivalence $f:M_1\to M_2$ between CR manifolds lifts canonically to a diffeomorphism $\hat{f}:\F_1^0\to\F_2^0$ between their 0-adapted coframe bundles in a manner that identifies their tautological forms, $\hat{f}^*\Phi_2=\Phi_1$. To find all local invariants of a CR structure via Cartan's method of equivalence (\cite{Gardner}), one attempts to complete the tautological form to a full coframing of a principal bundle over $M$ by choosing a complementary pseudoconnection form (\cite[Def.2.2]{BGG}). Such a choice depends on reducing the structure group of the coframe bundle as much as possible by successively adapting frames to higher order. 

For example, if $M$ is not Levi-flat, we could define a 1-adapted coframing to be a 0-adapted coframing which has the additional property that the matrix entries \eqref{levimatrix} of the Levi form take constant, specified values (such as $\ell$ being diagonalized with $p$ positive ones and $q$ negative ones on the diagonal). This reduces $\F^0$ to the subbundle of 1-adapted coframes whose structure group $G_1$ is matrices \eqref{Gzero} with the additional constraint,
\begin{align*}
&G_1\subset G_0: &\overline{a}^t\ell a=u\ell.
\end{align*}  
This reduction is not meaningful in the Levi-flat case; indeed, we have already noted that there are no local invariants for Levi-flat CR manifolds. In general, the degree of (non)degeneracy of the Levi form has substantial bearing on the application of the method of equivalence.

The opposite extreme of Levi-flatness is \emph{Levi-nondegeneracy}, when $\ell$ has signature $(p,q)$, $p+q=n$. In CR dimension $n=1$, Levi-nondegeneracy is the same as pseudo-convexity, and the corresponding equivalence problem was solved by Cartan (\cite{CartanCR}, \cite{Jacobowitz}). The general case was treated by Tanaka (\cite{Tan62}) using his modified version of Cartan's method that would later provide a valuable framework for understanding all parabolic geometries with canonical Cartan connections (\cite{Tan79}, \cite{CapSlovak}). Chern-Moser (\cite{ChernMoser}) offered an alternative solution using the standard method and emphasizing the link between the intrinsic geometry of CR manifolds and the extrinsic analysis of normal forms. The solution to the Levi-nondegenerate equivalence problem may be stated as follows.

\begin{TCM}
Let $M$ be a hypersurface-type CR manifold of dimension $2n+1$ whose Levi form has signature $(p,q)$, $p+q=n$. 
\begin{itemize}
\item There exists a canonically defined principal bundle $\mathcal{B}\to M$ whose structure group is isomorphic to the parabolic subgroup $\mathcal{P}\subset SU(p+1,q+1)$ given by the stabilizer of a complex line in $\mathbb{C}^{n+2}$ which is null for a Hermitian form $\h$ of signature $(p+1,q+1)$.

\item There exists a canonical Cartan connection $\gamma\in\Omega^1(\mathcal{B},\mathfrak{su}(p+1,q+1))$ whose curvature tensor $\dd\gamma+\tfrac{1}{2}[\gamma,\gamma]\in\Omega^2(\mathcal{B},\mathfrak{su}(p+1,q+1))$ and its covariant derivatives determine a complete set of local invariants of $M$.

\item The algebra of infinitesimal symmetries of $M$ has dimension $\leq n^2+4n+3$, and the upper bound is only achieved where curvature locally vanishes. In this case, $M$ is locally CR equivalent to the hyperquadric $\Q\subset\mathbb{CP}^{n+1}$ given by the complex projectivization of the $\h$-null cone in $\mathbb{C}^{n+2}$; i.e., $\Q=SU(p+1,q+1)/\mathcal{P}$.
\end{itemize}
\end{TCM}  

The real hyperquadric $\Q$ is the ``flat model" of Levi-nondegenerate CR geometry in the sense that it is locally characterized by a vanishing curvature tensor. When $\dim M=5$, a nondegenerate Levi form either has definite signature $(2,0)$ or split signature $(1,1)$, and the Cartan connection $\gamma$ takes values in $\mathfrak{su}(3,1)$ or $\mathfrak{su}(2,2)$, respectively. Thus, for the flat models $M=\Q$, the principal bundle $\mathcal{B}$ is isomorphic to one of the Lie groups $SU(3,1)$ or $SU(2,2)$, and $\gamma$ is exactly the Maurer-Cartan form of $\mathcal{B}$.  

\begin{rem}\label{MQbundlesrem}
The proof of the Tanaka-Chern-Moser Classification is constructive, building the \emph{Cartan Bundle} $\mathcal{B}$ out of normalized jets of local CR coframings so that the Cartan connection $\gamma$ prolongs the tautological form \eqref{tautform} to a full coframing of $\mathcal{B}$. By design, a CR map $f:M_1\to M_2$ lifts uniquely to a smooth map $\hat{f}:\mathcal{B}_1\to\mathcal{B}_2$ such that $\hat{f}^*\gamma_2=\gamma_1$, and conversely, a smooth map $\hat{f}:\mathcal{B}_1\to\mathcal{B}_2$ such that $\hat{f}^*\gamma_2=\gamma_1$ descends uniquely to a CR map $f:M_1\to M_2$. We are concerned with CR embeddings as in Definition \ref{CRembed} where $M_1=M$ is 3-dimensional and $M_2=\Q$ is a 5-dimensional hyperquadric. To avoid confusion, the Cartan bundle and connection of $M$ -- constructed in \S\ref{CR3sec} -- are called $\F^3$ and $\gamma$, while those of $\Q$ -- exhibited in \S\ref{hermframesec} -- are called $\mathcal{H}$ and $\mu$. 
\end{rem}

\vspace{\baselineskip}

\subsection{3-dimensional, Levi-nondegenerate CR Manifolds}\label{CR3sec}

This section closely follows \cite{BryantCR3} with only minor changes to notation, and omitting several details. Fix $n=1$ so that $\dim M=3$ and a local 0-adapted coframing \eqref{zeroadapted} consists of an $\mathbb{R}$-valued characteristic form $\varphi^0$ and a $\mathbb{C}$-valued CR form $\varphi^1$ such that $\varphi^0\wedge\varphi^1\wedge\overline{\varphi}^1\neq 0$. This coframing is a local section of the bundle $\pi:\F^0\to M$ of 0-adapted coframes, and as such it determines a local trivialization of $\F^0$ over which the tautological form \eqref{tautform} is 
\begin{align}\label{Philocal}
&\Phi=\left[\begin{array}{c}\kappa\\\eta\end{array}\right]=\left[\begin{array}{cc}u&0\\b&a\end{array}\right]\pi^*\left[\begin{array}{c}\varphi^0\\\varphi^1\end{array}\right],
&0\neq u\in C^\infty(\F^0),\quad 0\neq a\in C^\infty(\F^0,\mathbb{C}),\quad b\in C^\infty(\F^0,\mathbb{C}),
\end{align}
with the functions $u,a,b$ acting as $G_0$-valued fiber coordinates for $\F^0$. 

CR integrability \eqref{CRintforms} is automatic in dimension three, and in particular \eqref{levimatrix} reads
\begin{align*}
&\dd\varphi^0\equiv\im\ell\varphi^1\wedge\overline{\varphi}^1\mod\{\varphi^0\},
&\ell\in C^\infty(M).
\end{align*}
Levi-nondegeneracy says $\ell$ is non-vanishing, so in this case we reduce to the bundle $\F^1\subset\F^0$ of 1-adapted coframes with $\ell=1$, which reduces the structure group $G_0$ to
\begin{align}\label{Gone}
G_1=\left.\left\{\left[\begin{array}{cc}|a|^2&0\\b&a\end{array}\right]\in GL(\mathbb{R}\oplus\mathbb{C})\ \right|\ 
a\in \mathbb{C}\setminus\{0\},\ b\in\mathbb{C} \right\}.
\end{align}
After pulling back the tautological form \eqref{Philocal} along the inclusion $\F^1\hookrightarrow\F^0$, its exterior derivative can be expressed in terms of a pseudoconnection form taking values in the Lie algebra $\mathfrak{g}_1$ of \eqref{Gone},
\begin{align}\label{FoneSE}
&\dd \left[\begin{array}{c}\kappa\\\eta\end{array}\right]=
-\left[\begin{array}{cc}\alpha_0+\overline{\alpha}_0&0\\\beta_0&\alpha_0\end{array}\right]\wedge
\left[\begin{array}{c}\kappa\\\eta\end{array}\right]+
\left[\begin{array}{c}\im\eta\wedge\overline{\eta}\\0\end{array}\right];
&\alpha_0,\beta_0\in\Omega^1(\F^1,\mathbb{C}).
\end{align}
However, the structure equations \eqref{FoneSE} do not uniquely determine $\alpha_0$ and $\beta_0$ as they remain the same after a replacement 
\begin{align}\label{abprime}
&\left[\begin{array}{c}\alpha'_0\\\beta'_0\end{array}\right]=
\left[\begin{array}{c}\alpha_0\\\beta_0\end{array}\right]+
\left[\begin{array}{cc}s^1&0\\s^2&s^1\end{array}\right]\left[\begin{array}{c}\kappa\\\eta\end{array}\right],
&s^1,s^2\in C^\infty(\F^1,\mathbb{C}).
\end{align}

For any $\alpha_0',\beta_0'$ of the form \eqref{abprime}, the 1-forms $\kappa,\eta,\overline{\eta},\alpha'_0,\overline{\alpha}'_0,\beta'_0,\overline{\beta}'_0$ are called a 1-adapted coframing of $\F^1$. The bundle $\hat{\pi}:\hat{\F}^1\to\F^1$ of 1-adapted coframes of $\F^1$ features a tautological $\mathbb{R}\oplus\mathbb{C}\oplus\mathfrak{g}_1$-valued form whose $\mathbb{R}\oplus\mathbb{C}$-valued components are simply the $\hat{\pi}$ pullback of $\Phi$ (we recycle the names of the individual 1-forms), 
\begin{align*}
\hat{\pi}^*\Phi=\left[\begin{array}{c}\kappa\\\eta\end{array}\right]\in\Omega^1(\hat{\F}^1,\mathbb{R}\oplus\mathbb{C}),
\end{align*}
and whose $\mathfrak{g}_1$-valued components $\left[\begin{smallmatrix}\alpha+\overline{\alpha}&0\\\beta&\alpha\end{smallmatrix}\right]\in\Omega^1(\hat{\F}^1,\mathfrak{g}_1)$ satisfy the ``lifted" structure equations,
\begin{align}\label{hatFSE}
&\dd \left[\begin{array}{c}\kappa\\\eta\end{array}\right]=
-\left[\begin{array}{cc}\alpha+\overline{\alpha}&0\\\beta&\alpha\end{array}\right]\wedge
\left[\begin{array}{c}\kappa\\\eta\end{array}\right]+
\left[\begin{array}{c}\im\eta\wedge\overline{\eta}\\0\end{array}\right];
&\alpha,\beta\in\Omega^1(\hat{\F}^1,\mathbb{C}).
\end{align}
In particular, from \eqref{abprime} we see that 
\begin{align*}
&\alpha=\hat{\pi}^*\alpha_0-s^1\kappa,
&\beta=\hat{\pi}^*\beta_0-s^2\kappa-s^1\eta,
&&s^1,s^2\in C^\infty(\hat{\F}^1,\mathbb{C}),
\end{align*}
where $s^1,s^2$ now serve as fiber coordinates for $\hat{\pi}:\hat{\F}^1\to\F^1$. Differentiating the structure equations \eqref{hatFSE} yields  
\begin{align*}
\dd\left[\begin{array}{c}\alpha\\\beta\end{array}\right]=
-\left[\begin{array}{cc}\sigma^1_0&0\\\sigma^2_0&\sigma^1_0\end{array}\right]\wedge\left[\begin{array}{c}\kappa\\\eta\end{array}\right]+
\left[\begin{array}{c}-\im\beta\wedge\overline{\eta}-2\im\overline{\beta}\wedge\eta+R\eta\wedge\overline{\eta} \\
-\beta\wedge\overline{\alpha}\end{array}\right],
\end{align*}
for some $R\in C^\infty(\hat{\F}^1)$, with $\sigma_0^1,\sigma_0^2,\kappa,\eta,\alpha,\beta$ and their conjugates furnishing a coframing of $\hat{\F}^1$. The identity $\dd^2\alpha=0$ then reveals that we can restrict to a subbundle $\F^2\subset\hat{\F}^1$ whose sections are 2-adapted coframings defined by $R=0$, which reduces the real dimension of the fibers over $\F^1$ by one and forces $s^1$ and $\sigma=\sigma^1_0$ to be strictly $\mathbb{R}$-valued. Next, the same identity shows that we can reduce further to 3-adapted coframings corresponding to a subbundle $\F^3\subset\F^2$ where $s^2=0$, hence the real fiber dimension of $\F^3\to\F^1$ is one. 

The coframing of $\F^3$ given by the complex forms $\eta,\alpha,\beta$ and their conjugates, along with the real forms $\kappa,\sigma$, is globally defined on $\F^3$ and uniquely determined by the structure equations 
\begin{equation}\label{threeCRSE}
\begin{aligned}
\dd\kappa&=\im\eta\wedge\overline{\eta}-(\alpha+\overline{\alpha})\wedge\kappa,\\
\dd\eta&=-\beta\wedge\kappa-\alpha\wedge\eta,\\
\dd\alpha&=-\sigma\wedge\kappa-\im\beta\wedge\overline{\eta}-2\im\overline{\beta}\wedge\eta,\\
\dd\beta&=-\sigma\wedge\eta+\overline{\alpha}\wedge\beta+S\kappa\wedge\overline{\eta},\\
\dd\sigma&=(\alpha+\overline{\alpha})\wedge\sigma+\im\beta\wedge\overline{\beta}+\kappa\wedge(P\overline{\eta}+\overline{P}\eta),
\end{aligned}
\end{equation}
where $S,P\in C^\infty(\F^3,\mathbb{C})$ have differential identities
\begin{equation}\label{Bianchi}
\begin{aligned}
\dd S&=S(3\overline{\alpha}+\alpha)+U\kappa+P\eta+Q\overline{\eta},\\
\dd P&=P(3\overline{\alpha}+2\alpha)-\im S\overline{\beta}+W\kappa+R\eta+V\overline{\eta},
\end{aligned}
\end{equation}
for some $U,Q,W,V\in C^\infty(\F^3,\mathbb{C})$ and $R\in C^\infty(\F^3)$. $\F^3$ realizes the Cartan bundle $\mathcal{B}$ of $M$ whose existence is guaranteed by the Tanaka-Chern-Moser Classification; the Cartan connection is 
\begin{align}\label{Mcarcon}
\gamma=\left[\begin{array}{ccc}
-\tfrac{1}{3}(2\alpha+\overline{\alpha})&-\im\overline{\beta}&-\im\sigma\\
\eta&\tfrac{1}{3}(\alpha-\overline{\alpha})&\im\beta\\
-\im\kappa&\overline{\eta}&\tfrac{1}{3}(\alpha+2\overline{\alpha})\end{array}\right],
\end{align}
so that 
\begin{align*}
&\overline{\gamma}^t\h+\h\gamma=0,
&\h=\left[\begin{array}{crc}0&0&-1\\0&1&0\\-1&0&0\end{array}\right],
\end{align*}
and $\gamma$ is indeed $\mathfrak{su}(2,1)$-valued. The equations \eqref{Bianchi} engender the Bianchi identities of the curvature tensor $\dd\gamma+\gamma\wedge\gamma$. Furthermore, when $S=0\Rightarrow P=0$ so that curvature locally vanishes, \eqref{threeCRSE} are exactly the Maurer-Cartan equations of $\mathfrak{su}(2,1)$, as previously discussed.

If the curvature tensor of $M$ never vanishes, we can adapt to even higher order (see \cite[\S3.4.4]{BryantCR3}). First, differentiate the identities \eqref{Bianchi} to obtain 
\begin{equation}\label{dRQUVW}
\begin{aligned}
\dd R&= R(3\overline{\alpha}+3\alpha)+R_0'\kappa+R_1'\overline{\eta}+\overline{R}_1'\eta,\\
\dd Q&= Q(4\overline{\alpha}+\alpha)-5\im S\beta+U_1'\kappa+(V-\im U)\eta+Q'\overline{\eta},\\
\dd U&= U(4\overline{\alpha}+2\alpha)+4S\sigma+P\beta+Q\overline{\beta}+W\eta+U_1'\overline{\eta}+U_2'\kappa,\\
\dd V&\equiv V(4\overline{\alpha}+2\alpha)-\im S\sigma-4\im P\beta-\im Q\overline{\beta}+(R_1'-\im W)\eta +V'\kappa &\mod\{\overline{\eta}\},\\
\dd W&=W(4\overline{\alpha}+3\alpha)+5P\sigma+R\beta+(V-\im U)\overline{\beta}+(R_0'-\im|S|^2)\eta+V'\overline{\eta}+W'\kappa,
\end{aligned}
\end{equation}
for some $R_0'\in C^\infty(\F^3)$ and $R_1',Q',U_1',U_2',V',W'\in C^\infty(\F^3,\mathbb{C})$. For later use, we differentiate the first line of \eqref{dRQUVW} in order to record
\begin{equation}\label{dR'}
\begin{aligned}
\dd R_1'&\equiv R_1'(4\overline{\alpha}+3\alpha)-3\im R\beta+(R_0''-\tfrac{\im}{2}R_0')\eta+R_1''\kappa
&&\mod\{\overline{\eta}\},
\end{aligned}
\end{equation}
with additional $R_0''\in C^\infty(\F^3)$ and $R_1''\in C^\infty(\F^3,\mathbb{C})$.

Now observe that the identity \eqref{Bianchi} for $\dd S$ implies that if $S$ is nowhere zero, we can reduce to $\F^4\subset\F^3$ where $S=1$, and over $\F^4$ we have 
\begin{align}\label{alphareduction}
\alpha=\frac{1}{8}\Big( (U-3\overline{U})\kappa+ (P-3\overline{Q})\eta +(Q-3\overline{P})\overline{\eta}\Big).
\end{align}
At this point, the equation \eqref{Bianchi} for $\dd P$ shows that there is a subbundle $\F^5\subset\F^4$ on which $P=0$ and 
\begin{align}\label{betareduction}
\beta=\im(\overline{W}\kappa+\overline{V}\eta+R\overline{\eta}).
\end{align}
Finally, by the third line of \eqref{dRQUVW} we can reduce to 6-adapted coframes $\F^6\subset\F^5$ defined by $\Re U=\tfrac{1}{2}(U+\overline{U})=0$, with 
\begin{align}\label{sigmareduction}
&\sigma=-\frac{1}{4}\Re(U(4\overline{\alpha}+2\alpha)+Q\overline{\beta}+W\eta+U_1'\overline{\eta}+U_2'\kappa),
&\text{subject to}
&&\eqref{alphareduction},\eqref{betareduction}.
\end{align}

After these reductions to $\F^6\subset\F^3$, what remains of \eqref{threeCRSE} is  
\begin{equation}\label{FsixSE}
\begin{aligned}
\dd\kappa&=\im\eta\wedge\overline{\eta}-2\kappa\wedge(\overline{A}\eta+A\overline{\eta}),\\
\dd\eta&=A\eta\wedge\overline{\eta}+\im \kappa\wedge(B\eta+ C\overline{\eta}),
\end{aligned}
\end{equation}
where
\begin{align}\label{ABCreduction}
&A=\frac{Q}{8},
&B=\overline{V}+\frac{\im}{2}U,
&&C=R.
\end{align}
The structure equations (\ref{FsixSE}, \ref{ABCreduction}) for a 6-adapted coframing satisfy $\dd^2\kappa=\dd^2\eta=0$ by virtue of the identities \eqref{dRQUVW}. When $M$ is homogeneous under the action of its CR symmetry group, $A,B,C$ are constant and the differential conditions $\dd^2\kappa=\dd^2\eta=0$ simplify to algebraic relations 
\begin{align}\label{Fsixdsquared}
&B=\overline{B},
&AB=\overline{A}C.
\end{align}

Conversely, suppose $M$ is a 3-dimensional CR manifold admitting a 1-adapted coframing $\kappa,\eta$ satisfying \eqref{FsixSE}, where $A,B,C$ are constant with $C\in\mathbb{R}$. Such $M$ is locally homogeneous, and the identities $\dd^2\kappa=\dd^2\eta=0$ once again imply \eqref{Fsixdsquared}. Moreover, the 1-adapted coframing $\kappa,\eta$ determines a section $M\to\F^3$ along which the Cartan connection forms pull back to
\begin{equation}\label{ABCforms}
\begin{aligned}
\alpha&=-\im(|A|^2+\tfrac{3}{4}B)\kappa-3\overline{A}\eta+A\overline{\eta},\\
\beta&=\tfrac{1}{3}(AB+4\overline{A}C-4A^2\overline{A})\kappa+\im(\tfrac{1}{4}B-|A|^2)\eta+\im C\overline{\eta},\\
\sigma&=(\tfrac{5}{3}C(A^2+\overline{A}^2)-\tfrac{1}{3}|A|^4-\tfrac{13}{6}B|A|^2+\tfrac{1}{16}B^2-C^2)\kappa+\tfrac{1}{3}\Re\Big(\im(10CA-5B\overline{A}-4A\overline{A}^2)\eta\Big),
\end{aligned}
\end{equation}
according to the structure equations \eqref{threeCRSE}, the penultimate of which necessitates the relation
\begin{align}\label{hiddenrelation}
40A^3\overline{A}-10A^2B-28|A|^2C+9BC=6S.
\end{align}

\begin{rem}\label{flatvsnon} When $M$ is locally flat ($S=0\Rightarrow P=0$), higher-order functions in \eqref{dRQUVW} and \eqref{dR'} vanish in turn. On the other hand, if \eqref{FsixSE} constitutes a 6-adapted coframing ($S=1, P=0, \overline{U}=-U$) of a non-flat, homogeneous CR 3-fold, the section $M\to\F^3$ determined by $\kappa,\eta$ pulls back the higher-order coefficients to
\begin{align*}
&\begin{array}{l}
R=C, \\ 
Q=8A, \\ 
U=-\im(2|A|^2+\tfrac{3}{2}B), \\ 
V=\tfrac{1}{4}B-|A|^2, \\ 
W=\tfrac{\im}{3}(B\overline{A}+4CA-4A\overline{A}^2),\\
R_1'=6CA,\\
R_1''=\im\overline{A}C(4C-10A^2)-\tfrac{7\im}{2}ABC,
\end{array}
\begin{array}{l}
R_0'=0,\\
R_0''=33C|A|^2-\tfrac{3}{4}BC,\\
Q'=88A^2-5C,\\
U_1'=\tfrac{\im}{3}(20C\overline{A}-49BA-92A^2\overline{A}),\\
U_2'=8|A|^4-\tfrac{5}{2}B^2+4C^2-\tfrac{1}{3}C(52A^2+20\overline{A}^2),\\
V'=-\im(9|A|^4-\tfrac{3}{2}|A|^2B-\tfrac{1}{3}C(37A^2+5\overline{A}^2)+\tfrac{5}{16}B^2+C^2),\\
W'=\tfrac{1}{3}(4|A|^2(5AC-B\overline{A}-4\overline{A}|A|^2)+7ABC+2\overline{A}(B^2-2C^2)).
\end{array}
\end{align*}
\end{rem}

To summarize, we see that \eqref{FsixSE} -- with constants $A\in\mathbb{C}$, $B,C\in\mathbb{R}$ constrained by \eqref{Fsixdsquared} and \eqref{hiddenrelation} -- locally characterizes all 3-dimensional, Levi-nondegenerate CR manifolds which are homogeneous under the action of their CR symmetry groups. In the non-flat case $S\neq0$, these are the structure equations for the maximal CR symmetry algebras tangent to these symmetry groups. Homogeneous, Levi-nondegenerate CR 3-folds are classified (locally and globally) in \cite{CartanCR}, which exhibits local hypersurface realizations of each model labeled by capital Latin letters. Nurowski and Tafel (\cite{NurTaf}) offer alternative coordinate realizations and make explicit the reliance of Cartan's arguments on Bianchi's classification of 3-dimensional real Lie algebras, which are labeled with Roman numerals as in, e.g., \cite[\S8.2]{exactsolutionsbook}.

Bianchi's type I algebra is abelian; it is simply the vector space $\mathbb{R}^3$, which cannot serve as the symmetry algebra of a Levi-nondegenerate CR 3-fold. Indeed, for coordinates $(x,z)\in\mathbb{R}\times\mathbb{C}$, the 0-adapted coframing $\varphi^0=\dd x$, $\varphi^1=\dd z$ has trivial structure equations. Type II is the Heisenberg Lie algebra -- the negatively-graded part of $\mathfrak{su}(2,1)$ with respect to the grading induced by the parabolic subgroup $\mathcal{P}\subset SU(2,1)$ appearing in the Tanaka-Chern-Moser Classification -- whose structure equations are
\begin{align}\tag{II, A}
\eqref{FsixSE} \text{ with } A=B=C=0, \text{ subject to \eqref{hiddenrelation} with } S=0.
\end{align}
As the twin labels suggest, this is Cartan's A model, the CR 3-sphere, whose full group of CR symmetries is $SU(2,1)$. Skipping to the end of the list, Bianchi's VIII and IX are $\mathfrak{su}(1,1)$ and $\mathfrak{su}(2)$, respectively, each of which is the symmetry algebra of a parameter-family of homogeneous models:
\begin{align}\tag{VIII, C}
&(\ref{FsixSE}, \ref{hiddenrelation}) \text{ with } A=0,\ &B=-1,\ &&&C=0,\ &&S=0;\\
\tag{VIII, K}
&(\ref{FsixSE}, \ref{hiddenrelation}) \text{ with } A=0,\ &B<0,\ &&&C=\frac{2}{3B},\ &&S=1;\\
\tag{IX, D} 
&(\ref{FsixSE}, \ref{hiddenrelation}) \text{ with } A=0,\ &B=1,\ &&&C=0,\ &&S=0;\\
\tag{IX, L} 
&(\ref{FsixSE}, \ref{hiddenrelation}) \text{ with } A=0,\ &B>0,\ &&&C=\frac{2}{3B},\ &&S=1.
\end{align}
The rest of the algebras on Bianchi's list may be represented as extensions of the abelian Lie algebra $\mathbb{R}^2$ by a single, nontrivial derivation; i.e., a nonzero $2\times 2$ matrix. Type V extends $\mathbb{R}^2$ by a diagonal matrix, and the resulting structure equations are Levi-flat (see \S\ref{LFclasssec}). Type IV is represented as the extension by $\left[\begin{smallmatrix}1&1\\0&1\end{smallmatrix}\right]$, and with a suitable choice of basis its structure equations are 
\begin{align}\tag{IV, F}
&(\ref{FsixSE}, \ref{hiddenrelation}) \text{ with } A=2\sqrt[4]{6},\ &B=C=9\sqrt{6},\ &&S=1.
\end{align}
The remaining Bianchi types are parameter-families of algebras, indexed by $0< t\in\mathbb{R}$, each member of which has a unique homogeneous model. Type $\text{VII}_t$ is the extension of $\mathbb{R}^2$ by $\left[\begin{smallmatrix}t&1\\-1&t\end{smallmatrix}\right]$, which is rank-2 for every $t$, and its corresponding homogeneous model is
\begin{align}\tag{$\text{VII}_t$, H}
&(\ref{FsixSE}, \ref{hiddenrelation}) \text{ with } A=\frac{2t\sqrt{6}}{\sqrt[4]{6t^4+60t^2+54}},\ &B=C=\frac{54t^2+6}{\sqrt{6t^4+60t^2+54}},\ &&S=1.
\end{align}
Type $\text{VI}_t$ is the extension of $\mathbb{R}^2$ by $\left[\begin{smallmatrix}t&m\\m&t\end{smallmatrix}\right]$, where $m\in\mathbb{R}$ is determined by the curvature $S$. Homogeneous models are described in general by
\begin{align}\label{sixEgeneral}
&(\ref{FsixSE}, \ref{hiddenrelation}) \text{ with } A=\iota t\quad (\iota=1\text{ or }\im),\ &B=\iota^2C=\frac{9t^2-m^2}{4},\ &&S=\frac{\iota^2}{96}(t^2-m^2)(t^2-9m^2).
\end{align}
Note some exceptional values of $t$ ($\iota=m=1$) with zero curvature:
\begin{align}\tag{III=$\text{VI}_1$, B}
&(\ref{FsixSE}, \ref{hiddenrelation}) \text{ with } A=1,\ &B=C=2,\ &&S=0;\\
\tag{$\text{VI}_3$, E}
&(\ref{FsixSE}, \ref{hiddenrelation}) \text{ with } A=3,\ &B=C=20,\ &&S=0.
\end{align}
For the remaining values of $t$, we have 
\begin{align}\tag{$\text{VI}_t$, E}
&0<t\neq 1,3;
&\eqref{sixEgeneral} \text{ with } m \text{ such that } S=1.
\end{align}

\vspace{\baselineskip}

\section{Physical Motivation}\label{physics}

\subsection{Some History}\label{history}

The subject matter of this article is closely related to relativistic theories of radiation, both electromagnetic and gravitational. In general relativity, spacetime is an oriented 4-dimensional manifold $\mathcal{S}$. The distribution of mass-energy in $\mathcal{S}$ is encoded in a symmetric 2-tensor field, and Einstein's field equations say that this field prescribes the Ricci curvature of a Lorentzian metric; we suppose that $\mathcal{S}$ admits a solution $g\in\bigodot^2T^*\mathcal{S}$. An electromagnetic (EM) field is a 2-form $F\in\Omega^2(\mathcal{S})$ that may be interpreted as ($-\im$ times) the curvature of a principal connection on a $U(1)$-bundle over $\mathcal{S}$. $F$ is \emph{null} if it is $g$-orthogonal to itself and its Hodge dual. Null EM fields are associated with electromagnetic radiation. A gravitational field is the curvature of $g$, called \emph{null} if the Weyl tensor has Petrov type N (see \cite[Ch.5]{Oneill} for a pleasant introduction to the Petrov classification) -- the most degenerate type among non-conformally-flat metrics. More generally, a metric is \emph{algebraically special} if its Weyl tensor is at all degenerate; i.e., if it is of any Petrov type besides I. It is helpful to think of the motivation for the present work in the context of the search for a theoretical framework for gravitational radiation analogous to that of electromagnetic radiation.

The study of gravitational waves was initiated by Einstein in the beginning of the twentieth century; a brief history with references is given in \cite{TrautHistory}. For our purposes, it suffices to join the story \textit{in medias res}, when Trautman showed in \cite{Traut58} that gravitational fields satisfying a Sommerfeld radiation condition are asymptotically null. The next year, Robinson reported to the Royaumont Conference that null EM and gravitational fields determine a foliation of spacetime by a family of curves known as a shear-free null geodesic congruence, or SFNGC. These are discussed in detail in \S\ref{SFNGCsec}. Robinson also proved the converse for electromagnetic fields (\cite{robnullEM}); i.e., a SFNGC gives rise to a null EM field. Spacetimes admitting SFNGCs seemed to be natural candidates for a model of gravitational radiation, though the work \cite{Sachs} of Sachs established that these were not restricted to null gravitational fields. Indeed, joint work \cite{GoldbergSachs} with Goldberg would show that, away from sources of mass-energy -- that is, in a vacuum spacetime with a Ricci-flat metric -- every non-flat, algebraically special metric admits a SFNGC tangent to principal null directions of algebraic multiplicity $>1$. 

Robinson and Trautman (\cite{RTGRwaves}) produced a class of metrics corresponding to hypersurface-orthogonal SFNGC, including a model for radiation with spherical wavefronts. Then Kerr sought metrics corresponding to SFNGC that were not hypersurface-orthogonal (\cite{Kerr}), and in the process generalized the Schwarzschild solution to incorporate angular momentum, generating a model for spinning black holes (\cite{Oneill}). The Kerr metrics have Petrov type D.

Kerr's name is also attached to a theorem relating SFNGC of flat (Minkowski) spacetime to the objects of study of this article. In \S\ref{kerr}, we offer a geometric overview of the correspondence between subsets of Minkowski spacetime and those of a 5-dimensional real hyperquadric in the spirit of Penrose's Twistor program (\cite{Pen67,WardWells}), emphasizing the various symmetry groups involved. Then \S\ref{SFNGCsec} delves into SFNGC for general spacetimes and explains their connection to CR geometry, following \cite{RTflows} and \cite{NurTraut}. Finally, a sketch of the proof of the Kerr Theorem appears in \S\ref{kerrproof}, using explicit coordinate calculations as in \cite[\S5]{Tafel}.

\vspace{\baselineskip}

\subsection{The Kerr Theorem}\label{kerr}

$\mathbb{R}^n$ equipped with a symmetric, bilinear form $\bb$ of signature $(p,q)$, $p+q=n$, will be denoted $\mathbb{R}^{p,q}$. The complex-linear extension of $\bb$ to $\mathbb{C}^n$ is also called $\bb$, but in the complexification its signature is no longer significant. Hence, we reserve the notation $\mathbb{C}^{p,q}$ for when $\mathbb{C}^n$ carries a Hermitian form $\h$ of signature $(p,q)$, $p+q=n$, unrelated to any underlying real form. 

The setting of special relativity is Minkowski spacetime $\mathbb{M}$, an affine space with modeling vector space $\mathbb{R}^{1,3}$. Therefore, the Lorentz group $O(1,3)$ -- and its affine extension, the Poincar\'{e} group -- plays a central role in relativistic theories. However, when a relativistic theory (such as the electrodynamics expressed in Maxwell's equations) exhibits conformal invariance, the corresponding group of symmetries is larger. 

Conformal compactification of $\mathbb{R}^{p,q}$ is achieved by affixing a point ``at infinity" for each one in the $\bb$-null cone in order that inversion may be globally defined. The resulting quadric is the real-projectivization of the $\hat{\bb}$-null cone in $\mathbb{R}^{p+1,q+1}$, whose bilinear form wears a hat to distinguish it from that of $\mathbb{R}^{p,q}$. The group of (oriented) conformal symmetries of compactified Minkowski spacetime $\mathbb{M}^c$ is the symmetry group of the $\hat{\bb}$-null cone in $\mathbb{R}^{2,4}$; i.e., $SO(2,4)$.  

The Pl\"{u}cker embedding sends the Grassmannian $Gr(2,\mathbb{C}^4)$ of complex 2-planes in $\mathbb{C}^4$ into the complex-projective space $\mathbb{P}(\Lambda^2\mathbb{C}^4)=\mathbb{CP}^5$, and its image is the quadric given by the projectivization of the $\hat{\bb}$-null cone in $\mathbb{C}^6$. This may be considered a geometric analog of the Lie algebra isomorphism $\mathfrak{sl}_4\mathbb{C}\cong\mathfrak{so}_6\mathbb{C}$. Moreover, the Grassmannian $Gr^0(2,\mathbb{C}^{2,2})\subset Gr(2,\mathbb{C}^4)$ of totally $\h$-null 2-planes embeds onto $\mathbb{M}^c$ by analogy to the isomorphism $\mathfrak{su}(2,2)\cong\mathfrak{so}(2,4)$.

Both $\mathbb{CP}^3$ and $Gr(2,\mathbb{C}^4)$ are partial flag manifolds associated to $\mathbb{C}^4$; to these we add $F_{1,2}\mathbb{C}^4$ consisting of pairs $(l,\Pi)$ of a complex line and plane (respectively) satisfying $l\subset\Pi\subset\mathbb{C}^4$. With the projection maps $\lambda(l,\Pi)=l$ and $\pi(l,\Pi)=\Pi$ we obtain the double fibration 
\begin{align}\label{Cdoublefiber}
\xymatrix@=1em{
&\ar_{\lambda}@{->}[dl] F_{1,2}\mathbb{C}^4\ar^{\pi}@{->}[dr] &\\
\mathbb{CP}^3& &Gr(2,\mathbb{C}^4)
}
\end{align}
lying at the heart of Penrose's Twistor theory, which concerns the correspondence between subsets of $\mathbb{CP}^3$ and $Gr(2,\mathbb{C}^4)$ via the images of $\lambda\circ\pi^{-1}$ and $\pi\circ\lambda^{-1}$. To clarify some of the physical motivation for this framework, we restrict to $\h$-isotropic flags $F_{1,2}^0\mathbb{C}^{2,2}\subset F_{1,2}\mathbb{C}^4$, so that \eqref{Cdoublefiber} becomes  
\begin{align}\label{Rdoublefiber}
\xymatrix@=1em{
&\ar_{\lambda}@{->}[dl] F^0_{1,2}\mathbb{C}^{2,2}\ar^{\pi}@{->}[dr] &\\
\Q& &\mathbb{M}^c=Gr^0(2,\mathbb{C}^{2,2}),
}
\end{align}
where $\Q\subset\mathbb{CP}^3$ is the 5-dimensional real hyperquadric given by the complex projectivization of the $\h$-null cone $\N\subset\mathbb{C}^{2,2}$. The trajectory of a massless particle in $\mathbb{M}$ is tangent to a $\bb$-null (affine) line, and each such line corresponds to a point in $\Q$. Physicists refer to a foliation of $\mathbb{M}$ by null lines as a \emph{null congruence}, the relevance of which to the present work is stated in the 

\begin{kerr}
A null congruence of $\mathbb{M}$ corresponds to a CR submanifold of $\Q$ if and only if it is shear-free.
\end{kerr}

The Kerr Theorem first appeared in print in \cite[\S VIII]{Pen67}, where it is stated that a shear-free null congruence is representable in $\mathbb{CP}^3$ as the intersection of $\Q$ with a complex-analytic surface (or a limiting case of such intersections); see also \cite[Ch.6]{PenBook2}. The proof sketch in \S\ref{kerrproof} makes this construction explicit. The version we've stated is closer to \cite[Thm.7]{NurTraut}.

\vspace{\baselineskip}

\subsection{Shear-Free Null Geodesic Congruences}\label{SFNGCsec}

In this section we follow \cite{RTflows} and \cite{NurTraut}.
Let $\mathcal{S}$ be a smooth, 4-dimensional manifold with a line bundle $K\subset T\mathcal{S}$ whose fibers are spanned by a nowhere-vanishing vector field $k\in\Gamma(K)$, which determines a smooth flow 
\begin{align*}
\phi:I\times \mathcal{S}\to \mathcal{S},
\end{align*}
where $I\subseteq\mathbb{R}$ is some open interval containing zero. For fixed $x\in \mathcal{S}$ and variable $t\in I$, $\phi(t,x)$ is the integral curve of $k$ passing through $x$ when $t=0$, and $\mathcal{S}$ is foliated by these flow curves. For fixed $t\in I$,
\begin{equation}\label{flowmap}
\begin{aligned}
\phi_t:\mathcal{S}&\to \mathcal{S}\\
x&\mapsto \phi(t,x)
\end{aligned}
\end{equation} 
is a diffeomorphism whose pushforward $\phi_{t*}:T\mathcal{S}\to T\mathcal{S}$ satisfies 
\begin{align}\label{Kflow}
\phi_{t*}K_x=K_{\phi(t,x)},
\end{align}
and therefore descends to a well-defined map on the quotient bundle $T\mathcal{S}/K\to \mathcal{S}$. Thus, the family $\{\phi_t:t\in I\}$ of diffeomorphisms provides linear isomorphisms between the spaces $T_{\phi(t,x)}\mathcal{S}/K_{\phi(t,x)}$ for any fixed $x\in \mathcal{S}$, and we see that the quotient bundle $T\mathcal{S}/K$ has the same fibers as the tangent bundle of the leaf space $M $; i.e., the 3-dimensional quotient manifold of equivalence classes $[x]$ of points $x\in \mathcal{S}$, where two points are equivalent if they lie in the same leaf of the foliation (the same flow curve),
\begin{align}\label{Tleaf}
&T_{[x]}M \cong T_{\phi(t,x)}\mathcal{S}/K_{\phi(t,x)}
&\forall t\in I.
\end{align} 

\begin{rem}\label{orientrem}
In general, a 4-manifold need not admit a globally defined, non-vanishing tangent vector field, nor should the entire leaf space of a foliation necessarily inherit a global manifold structure. However, our considerations are local in nature and we will continue to implicitly assume that $\mathcal{S}$ is such that our constructions are well-defined. In particular, we may also take $\mathcal{S}$ to be orientable. If $\omega\in\Omega^4(\mathcal{S})$ is a volume form, then the contraction $k\lrcorner\omega\in\Omega^3(\mathcal{S})$ vanishes on $K$ and descends to a 3-form on $T\mathcal{S}/K$. Note that $k\lrcorner\omega$ does not determine a well-defined volume form on $M $ unless $\mathcal{L}_k\omega=0$, where $\mathcal{L}_k$ denotes the Lie derivative along $k$. However, the sign of $k\lrcorner\omega$ on any ordered basis of \eqref{Tleaf} is sufficient to determine whether a volume form on $M $ is positively or negatively oriented relative to $k\lrcorner\omega$, and so determines a choice of orientation on $M $.
\end{rem} 

Suppose that $\mathcal{S}$ is equipped with a non-degenerate metric $g\in\bigodot^2T^*\mathcal{S}$. For the moment, we make no assumptions about the signature of $g$. The one-form $\kappa\in\Omega^1(\mathcal{S})$ dual to $k$ has as its kernel a rank-3 distribution 
\begin{align*}
&\kappa=k\lrcorner g
&\leadsto
&&\ker\kappa=K^\bot\subset T\mathcal{S}.
\end{align*} 
\begin{definition}
The flow of $k$ is \emph{conformally geodesic} if it preserves the distribution $K^\bot$,
\begin{align*}
&\phi_{t*}K_x^\bot=K_{\phi(t,x)}^\bot
&\forall t\in I, x\in \mathcal{S}.
\end{align*}
Equivalently, the flow of $k$ is \emph{conformally geodesic} when
\begin{align}\label{Lkkappa}
&\kappa\wedge\phi_t^*\kappa=0
&\Rightarrow
&&\kappa\wedge\mathcal{L}_k\kappa=0.
\end{align}
\end{definition}
Hence, a conformally geodesic flow not only identifies the fibers of $K$ along a flow curve as in \eqref{Kflow}, but also the fibers of $K^\bot$. The implications of this for the leaf space $M $ depend on the metric properties of $k$. If $g(k_x,k_x)\neq 0$ for every $x\in \mathcal{S}$, then $T\mathcal{S}=K\oplus K^\bot$ and $T_{[x]}M \cong K^\bot_{\phi(t,x)}$ for every $t\in I$. On the other hand, if $g$ has mixed signature and $g(k,k)=0$, then $K\subset K^\bot$ and $M $ inherits a well-defined, rank-2 distribution 
\begin{align}\label{Dleaf}
D\subset TM \quad\text{with fibers}\quad D_{[x]}\cong K^\bot_{\phi(t,x)}/K_{\phi(t,x)}\quad\forall t\in I.
\end{align}
We also have when $k$ is null that $\kappa$ descends to the quotient bundle $T\mathcal{S}/K$, and the additional condition \eqref{Lkkappa} that $k$ is conformally geodesic further implies that $\kappa$ determines a well-defined, non-vanishing one-form (of the same name) on $M $, which annihilates \eqref{Dleaf}.

\begin{rem}\label{georeparam}
Condition \eqref{Lkkappa} is always invariant under conformal scaling of $\kappa$, and when $k$ is null it is even invariant under scaling of $k$ by a non-vanishing function, which effects a reparameterization of the flow curves of $k$. 
\end{rem}

Henceforth, we restrict to the case that $g$ has Lorentzian signature $(1,3)$ and $k$ is $g$-null with a conformally geodesic flow. The foliation of $\mathcal{S}$ by flow curves is now called a \emph{null geodesic congruence}, the fibers of the quotient bundle $K^\bot/K$ are called \emph{screen spaces}, and the geometry of the null congruence may be understood intuitively in terms of the following illustration regarding optical scalars (\cite[\S5.7]{Oneill}). In relativity, light propagates in null directions; imagine a beam of light casting the shadow of an opaque disk onto a 2-dimensional screen placed orthogonal to its (null) direction. As the screen is moved along the flow curve, this circular image might be rotated, enlarged, or distorted into an ellipse of greater eccentricity. If the latter, non-conformal distortion does not occur, the null congruence is \emph{shear-free}. The precise geometric definition applies to arbitrary conformally geodesic flows.  

\begin{definition}
A conformally geodesic flow is \emph{shear-free} if it preserves the conformal class of $g$ restricted to $K^\bot$; i.e., for any $t\in I$ and $x\in \mathcal{S}$, there is some $s\in\mathbb{R}$, $s>0$ such that
\begin{align}\label{shearfree}
\phi_{t}^*(g|_{K^\bot_{\phi(t,x)}})=sg|_{K^\bot_x},
\end{align}
so that in particular,
\begin{align}\label{Lkg}
\mathcal{L}_kg=ag+\kappa\odot\alpha
\end{align}
for some $a\in C^\infty(\mathcal{S})$ and $\alpha\in\Omega^1(\mathcal{S})$.
\end{definition}

\begin{rem}\label{SFNGCclassrem}
Using general properties of the Lie derivative, it is straightforward to confirm that \eqref{Lkg} is maintained under rescaling of $k$ by a non-vanishing function, albeit for different $a,\alpha$. Along with Remark \ref{georeparam}, this shows that a shear-free null geodesic congruence (SFNGC) is independent of the choice of $k$ spanning $K$. Note that if $k$ is $g$-null, it is also $\tilde{g}$-null, where 
\begin{align}\label{SFNGCclass}
&\tilde{g}=fg+\kappa\odot\xi,
&0<f\in C^\infty(\mathcal{S}),\ \xi\in\Omega^1(\mathcal{S}),
\end{align}
and $\tilde{\kappa}=k\lrcorner\tilde{g}$ is a rescaling of $\kappa$. Here again, the properties of the Lie derivative show that $\tilde{g}$ satisfies \eqref{Lkg} whenever $g$ does, so the class \eqref{SFNGCclass} of metrics associated to given SFNGC is manifestly larger than a conformal class of metrics.
\end{rem}

For null $k$, $g|_K=0$ and we see from \eqref{shearfree} that a SFNGC determines a well-defined conformal structure on the subbundle \eqref{Dleaf} of the leaf space. As such, we can define an almost-complex structure on $M $,
\begin{align*}
&J:D\to D,
&J^2=-\mathbbm{1},
\end{align*}
by taking $J_{[x]}$ to be a rotation by $\frac{\pi}{2}$ in $D_{[x]}$. (There are two choices for the direction of the rotation -- clockwise or counterclockwise -- in each $D_{[x]}$. Take the one which is positively oriented for the orientation induced by the semi-Riemannian volume form on $\mathcal{S}$; see Remark \ref{orientrem}.)

Thus we see that a SFNGC induces a CR structure on the 3-dimensional leaf space $M $, with $\kappa$ serving as a characteristic form. To this we may add a CR form $\eta\in\Omega^1(M ,\mathbb{C})$ so that $\kappa,\eta,\overline{\eta}$ is a 0-adapted CR coframing. CR integrability is automatic in dimension three, but pseudo-convexity is not. In the Levi-flat case, $\kappa\wedge\dd\kappa=0$ and $M $ is foliated by complex curves; the original curves of our SFNGC are \emph{hypersurface-orthogonal}, as one would expect from a spacetime featuring radiating wave fronts. More interesting from a CR perspective is the Levi-nondegenerate case corresponding to ``twisting" congruences $\kappa\wedge\dd\kappa\neq 0$, including the Kerr spacetime which describes a rotating black hole.

\vspace{\baselineskip}

Conversely, suppose that $M $ is a 3-dimensional CR manifold with a 0-adapted coframing $\kappa,\eta,\overline{\eta}$, and set $\mathcal{S}=\mathbb{R}\times M $. We use the same names $\kappa,\eta,\overline{\eta}$ to denote their pullbacks along the projection $\mathcal{S}\to M $. Take $k\in\Gamma(T\mathcal{S})$ to be $k=\frac{\partial}{\partial r}$ where $r$ is the Cartesian coordinate of $\mathbb{R}$, and choose any $\rho\in\Omega^1(\mathcal{S})$ with $\rho(k)=1$; i.e., $\rho\equiv\dd r\mod\{\kappa,\eta,\overline{\eta}\}$. The metric
\begin{align}\label{crmetric}
g=\kappa\odot\rho-\eta\odot\overline{\eta}
\end{align}
has signature $(1,3)$ and satisfies $g(k,k)=0$ as well as $\kappa=k\lrcorner g$. The flow curves of $k$ are the $r$-coordinate curves of $\mathcal{S}$, and Lie derivatives along $k$ can be computed via H. Cartan's formula, yielding 
\begin{align}\label{Lkcr}
&\mathcal{L}_k\kappa=\mathcal{L}_k\eta=\mathcal{L}_k\overline{\eta}=0,
&\mathcal{L}_k\rho\equiv0\mod\{\kappa,\eta,\overline{\eta}\},
\end{align}
whence both conditions \eqref{Lkkappa} and \eqref{Lkg} are verified. This establishes a correspondence
\begin{align}\label{SFNGCCR}
\{\text{SFNGC on 4-manifolds}\}\stackrel{\text{(local)}}{\longleftrightarrow}\{\text{CR structure on 3-manifolds}\}
\end{align}

Now suppose we submit our coframing on $M $ to a 0-adapted transformation,
\begin{align}\label{crmetrictrans}
&\left[\begin{array}{c}\kappa'\\\eta'\end{array}\right]=
\left[\begin{array}{cc}u&0\\b&a\end{array}\right]\left[\begin{array}{c}\kappa\\\eta\end{array}\right];
&u\in C^\infty(M ),\ a,b\in C^\infty(M ,\mathbb{C}),\ u,a\neq 0,
\end{align}
and write the metric $\tilde{g}$ as in \eqref{crmetric}. In terms of our original coframing, we obtain
\begin{align*}
\tilde{g}&=\kappa'\odot\rho-\eta'\odot\overline{\eta}'\\
%&=u\kappa\odot\rho-|a|^2\eta\odot\overline{\eta}-\kappa\odot(|b|^2\kappa+a\overline{b}\eta+\overline{a}b\overline{\eta})\\
&=|b|^2g+\kappa\odot((u-|b|^2)\rho-|b|^2\kappa-a\overline{b}\eta-\overline{a}b\overline{\eta})\\
&=fg+\kappa\odot\xi
\end{align*}
as in \eqref{SFNGCclass}. Following our initial selection of $\rho$, the ambiguity of the metric $g$ due to our choice \eqref{crmetrictrans} of coframing on $M $ is measured by 5 real functions of 3 variables, rather than the 5 functions of 4 variables apparent in the full class of metrics discussed in Remark \ref{SFNGCclassrem}. However, if we allow the fiber coordinates $u,a,b$ of our G-structure \eqref{crmetrictrans} to vary with $r$, then the structure group of our bundle of 0-adapted CR frames exactly parameterizes the class of metrics associated to this SFNGC (note that the Lie derivatives along $k$ of our CR coframing will no longer vanish identically as in \eqref{Lkcr} if $u,a,b$ depend on $r$; \eqref{Lkkappa} and \eqref{Lkg} will hold nonetheless).

\vspace{\baselineskip}

The correspondence \eqref{SFNGCCR} raises several questions, the first of which is presented as Problem 1 in \cite{NurTraut}, and the second of which was communicated to the author by Pawe\l{} Nurowski:
\begin{itemize}
\item Which CR structures lift to a SFNGC whose class \eqref{SFNGCclass} of metrics contains a solution to Einstein's vacuum field equations?

\item The Goldberg-Sachs theorem says that there are two SFNGC associated to a metric of Petrov type D; are the two corresponding CR structures always equivalent?

\item Which CR structures lift to a SFNGC whose class \eqref{SFNGCclass} of metrics contains one that is (conformally) flat?
\end{itemize}
The Kerr Theorem offers the answer to the final question: \emph{those that are embedded within the real hyperquadric $\Q$}. The present article attempts to answer the inevitable follow-up question: \emph{which are those?}

\vspace{\baselineskip}

\subsection{Kerr Theorem Proof Sketch}\label{kerrproof}

We argue as in \cite[\S5]{Tafel}. Remark \ref{SFNGCclassrem} says that we can scale the vector field $k\in\Gamma(K)$ tangent to our SFNGC at will, so we are less occupied with null \emph{vectors} tangent to Minkowski spacetime $\mathbb{M}$ than we are with null \emph{directions}. The projectivized $\bb$-null cone in $\mathbb{R}^{1,3}$ is the (Riemann) 2-sphere $\mathbb{CP}^1$, hence a single stereographic coordinate $\zeta\in\mathbb{C}$ suffices to parameterize all null tangent directions in each $T_x\mathbb{M}$, with the exception of one direction ``at infinity." In standard coordinates $(x_0,x_1,x_2,x_3)\in\mathbb{M}$, the metric is diagonal,
\begin{align*}
g=\dd x_0\odot\dd x_0-\dd x_1\odot\dd x_1-\dd x_2\odot\dd x_2-\dd x_3\odot\dd x_3.
\end{align*}
Introducing null and complexified coordinates
\begin{align*}
&u=x_0-x_3,
&v=x_0+x_3,
&&w=x_1+\im x_2,
&&\overline{w}=x_1-\im x_2,
\end{align*}
brings $g$ into the form
\begin{align*}
g=\dd u\odot\dd v-\dd w\odot\dd\overline{w}.
\end{align*}
A general null vector field and its dual form are, up to real scale, 
\begin{align*}
&k=\frac{\partial}{\partial v}-\zeta\frac{\partial}{\partial w}-\overline{\zeta}\frac{\partial}{\partial \overline{w}}+|\zeta|^2\frac{\partial}{\partial u},
&\kappa=\dd u+\zeta\dd\overline{w}+\overline{\zeta}(\dd w+\zeta\dd v),
&&\zeta\in C^\infty(\mathbb{M},\mathbb{C}),
\end{align*}
while the SFNGC ``at infinity" is given by the $u$-coordinate lines. In the latter case, the remaining coordinates $v,w,\overline{w}$ descend to the leaf-space of $u$-coordinate lines, which is the Levi-flat $\mathbb{R}\times\mathbb{C}$, and this corresponds to a CR structure in $\Q$ that is tangent to a complex curve.

In the general case we can write
\begin{align*}
&g=\kappa\odot\dd v-\eta\odot\overline{\eta},
&\eta=\dd w+\zeta\dd v,
\end{align*}
and after computing Lie derivatives,
\begin{align*}
&\mathcal{L}_k\kappa=\dd\zeta(k)\dd\overline{w}+\dd\overline{\zeta}(k)\dd w+(\overline{\zeta}\dd\zeta(k)+\zeta\dd\overline{\zeta}(k))\dd v,
&\mathcal{L}_k\eta=\overline{\mathcal{L}_k\overline{\eta}}=-\dd\zeta+\dd\zeta(k)\dd v,
\end{align*}
we see that conditions \eqref{Lkkappa} and \eqref{Lkg} hold when
\begin{align}
0&=\dd\zeta(k)=\dd\overline{\zeta}(k),\tag*{conformally geodesic:}\\
0&=\kappa\wedge\eta\wedge\mathcal{L}_k\eta,\tag*{shear-free:}
\end{align}
where the second becomes equivalent to 
\begin{align}\label{minksf}
\dd(u+\zeta\overline{w})\wedge\dd(w+\zeta v)\wedge\dd\zeta=0.
\end{align}
Name the three $\mathbb{C}$-valued functions 
\begin{align*}
&z_1=u+\zeta\overline{w},
&z_2=w+\zeta v,
&&z_3=\zeta,
\end{align*}
and observe that 
\begin{align}\label{Qlocalhypersurfaceeqn}
\im(z_1-\overline{z}_1+z_2\overline{z}_3-z_3\overline{z}_2)=0.
\end{align}
If $Z_0,Z_1,Z_2,Z_3$ are coordinates for $\mathbb{C}^{2,2}$ with the Hermitian form
\begin{align*}
\h(Z,W)=\im(Z_1\overline{W}_0-Z_0\overline{W}_1+Z_2\overline{W}_3-Z_3\overline{W}_2),
\end{align*}
then \eqref{Qlocalhypersurfaceeqn} describes the projectivization in $\mathbb{CP}^3$ of the $\h$-null cone $\h(Z,Z)=0$ in the affine coordinate neighborhood $Z_0\neq0$, via projective coordinates $[Z_0:Z_1:Z_3:Z_4]=[1:z_1:z_2:z_3]$. Moreover, if $\zeta=z_3$ is implicitly defined by $H(z_1,z_2,z_3)=0$, where $H$ is holomorphic (and not constant) in $z_1,z_2,z_3$, then the 3-form $\dd z_1\wedge\dd z_2\wedge\dd z_3$ vanishes on the subbundle $\dd H=0$ of the complexified tangent bundle of $\mathbb{CP}^3$, and over the hyperquadric $\Q$ locally defined by \eqref{Qlocalhypersurfaceeqn}, this is exactly the shear-free condition \eqref{minksf}. The level set $H=0$ is a complex-analytic surface in $\mathbb{CP}^3$ whose intersection with the real hyperquadric $\Q$ defines a 3-dimensional CR submanifold of $\Q$. 

\vspace{\baselineskip}

For the remaining details, consult \cite[\S VIII]{Pen67}, \cite[\S5]{Tafel}, or \cite[Ch.6]{PenBook2}.

\vspace{\baselineskip}

\section{Moving Frames Over Embedded 3-folds}\label{movingframesec}

%%%%%%%%%%%%%%%%%%%%%%%%%%%%%%%%%%%%%%%%%%%%%%%%%%%%%%%%%%%%%%%%%%%%%
\subsection{Hermitian Frames of $\mathbb{C}^4$}\label{hermframesec}
%%%%%%%%%%%%%%%%%%%%%%%%%%%%%%%%%%%%%%%%%%%%%%%%%%%%%%%%%%%%%%%%%%%%%

Let $\underline{e}=(e_0,e_1,e_2,e_3)$ denote the standard basis of column vectors for $\mathbb{C}^4$ and recall that $\e$ is the natural exponential and $\im=\sqrt{-1}$. Fix index ranges and constants 
\begin{align*}
&0\leq i,j\leq 3, 
&\epsilon=\pm1,
&&\delta_\epsilon=\left\{\begin{smallmatrix}0,&\epsilon=1\\1,&\epsilon=-1\end{smallmatrix}\right.
&&\Rightarrow \epsilon=(-1)^{\delta_\epsilon}.
\end{align*}
The Hermitian form $\h$ of signature $(3-\delta_\epsilon,1+\delta_\epsilon)$ acts on vectors ${\tt z}=z^ie_i$ and ${\tt w}=w^je_j$ via
\begin{align*}
\h({\tt w},{\tt z})%&=\left[\begin{array}{cccc}\overline{w}^0&\overline{w}^1&\overline{w}^2&\overline{w}^3\end{array}\right]\left[\begin{array}{r}-\im z^3\\ z^1\\\epsilon z_2\\\im z^0\end{array}\right]\\
%&=\overline{{\tt w}}^t\left[\begin{array}{cccc}0&0&0&\im\\0&1&0&0\\0&0&\epsilon&0\\-\im&0&0&0\end{array}\right]{\tt z}\\
&=\im(\overline{w}^0z^3-\overline{w}^3z^0)+\overline{w}^1z^1+\epsilon\overline{w}^2z^2.
\end{align*}
A \emph{Hermitian frame} is an ordered, complex basis $\underline{\vv}=(\vv_0,\vv_1,\vv_2,\vv_3)$ of $\mathbb{C}^4$ such that 
\begin{align}\label{hvivj}
\h(\vv_i,\vv_j)=(\pm1)^{\delta_\epsilon}\left[\begin{array}{cccc}0&0&0&\im\\0&1&0&0\\0&0&\epsilon&0\\-\im&0&0&0\end{array}\right].
\end{align}
Denote by $\hat{\mathcal{H}}$ the collection of all Hermitian frames, and note that $\hat{\mathcal{H}}$ is identified with the unitary group $U(3-\delta_\epsilon,1+\delta_\epsilon)$ by fixing $\underline{e}$ as the identity and taking $\underline{\vv}$ to be the matrix whose column vectors are the basis vectors of $\underline{\vv}$. 

\begin{rem}\label{deltaepsilonrem1}
The notation $(\pm1)^{\delta_\epsilon}$ in \eqref{hvivj} means the sign is allowed to change when $\epsilon=-1$ but not when $\epsilon=1$. This is to avoid privileging either of $\vv_1,\vv_2$ as necessarily positive-definite when $\epsilon=-1$. Both vectors are called ``orthonormal" regardless of the value of $\epsilon$.  
\end{rem}

The symbol $\vv_i$ will also refer to the $\mathbb{C}^4$-valued function mapping a Hermitian frame to its $i^{\text{th}}$ basis vector: $\vv_i\in C^\infty(\hat{\mathcal{H}},\mathbb{C}^4)$. These functions are differentiated via the Maurer-Cartan (MC) forms of $U(3-\delta_\epsilon,1+\delta_\epsilon)$,
\begin{align}\label{dv}
&\dd \vv_i=\mu^j_i\vv_j,
&\mu\in\Omega^1(\hat{\mathcal{H}},\mathfrak{u}(3-\delta_\epsilon,1+\delta_\epsilon)).
\end{align}  
In our representation of this Lie algebra, we can write
\begin{align}\label{UMCform}
\mu=\left[\begin{array}{cccc}
\lambda&-\im\overline{\xi}&-\im\overline{\phi}_2&\psi\\
\eta&\im\rho&-\overline{\phi}_1&\xi\\
\zeta&\epsilon\phi_1&\im\tau&\epsilon\phi_2\\
\kappa&\im\overline{\eta}&\epsilon\im\overline{\zeta}&-\overline{\lambda}\end{array}\right],
\end{align}
where $\kappa,\psi\in\Omega^1(\hat{\mathcal{H}})$ and the rest are $\mathbb{C}$-valued, so that $\overline{\mu}^t\h+\h\mu=0$. The MC equations $\dd\mu=-\mu\wedge\mu$ are
\begin{equation}\label{UMCeq}
\begin{aligned}
\dd\kappa&=\im\eta\wedge\overline{\eta}+(\lambda+\overline{\lambda})\wedge\kappa+\epsilon\im\zeta\wedge\overline{\zeta},\\
\dd\eta&=(\lambda-\im\rho)\wedge\eta-\xi\wedge\kappa+\overline{\phi}_1\wedge\zeta,\\
\dd\psi&=\psi\wedge(\lambda+\overline{\lambda})-\im\xi\wedge\overline{\xi}-\epsilon\im\phi_2\wedge\overline{\phi}_2,\\
\dd\xi&=\psi\wedge\eta+\xi\wedge(\overline{\lambda}+\im\rho)+\epsilon\overline{\phi}_1\wedge\phi_2,\\
\dd\lambda&=\im\overline{\xi}\wedge\eta+\im\overline{\phi}_2\wedge\zeta-\psi\wedge\kappa,\\
\dd\rho&=\epsilon\im\phi_1\wedge\overline{\phi}_1-\overline{\xi}\wedge\eta-\xi\wedge\overline{\eta},\\
\dd\phi_1&=\epsilon\im\zeta\wedge\overline{\xi}+\im(\rho-\tau)\wedge\phi_1-\im\phi_2\wedge\overline{\eta},\\
\dd\phi_2&=\epsilon\psi\wedge\zeta+\phi_2\wedge(\overline{\lambda}+\im\tau)+\xi\wedge\phi_1,\\
\dd\zeta&=(\lambda-\im\tau)\wedge\zeta-\epsilon\phi_1\wedge\eta-\epsilon\phi_2\wedge\kappa,\\
\dd\tau&=\zeta\wedge\overline{\phi}_2+\overline{\zeta}\wedge\phi_2-\epsilon\im\phi_1\wedge\overline{\phi}_1.
\end{aligned}
\end{equation}

Define $\det:\hat{\mathcal{H}}\to U(1)\subset\mathbb{C}$ as usual by
\begin{align*}
\vv_0\wedge \vv_1\wedge \vv_2\wedge \vv_3=\det(\underline{\vv})e_0\wedge e_1\wedge e_2\wedge e_3,
\end{align*}
and let $\mathcal{H}\subset\hat{\mathcal{H}}$ denote the collection of \emph{oriented} Hermitian frames satisfying $\det(\underline{\vv})=1$. From any Hermitian frame $\underline{\vv}$ one obtains an oriented Hermitian frame in a variety of ways; e.g.,   
\begin{align}\label{orient}
\underline{\vv}\mapsto(\vv_0,\vv_1,\overline{\det(\underline{\vv})}\vv_2,\vv_3),
\end{align}
in this case preserving the vectors $\vv_0,\vv_1,\vv_3$. $\mathcal{H}$ is identified with the special unitary group $SU(3-\delta_\epsilon,1+\delta_\epsilon)$. Keeping the same names after pulling back the MC forms \eqref{UMCform} along the inclusion $\mathcal{H}\hookrightarrow\hat{\mathcal{H}}$ exhibits the trace-free condition of the special unitary Lie algebra $\mathfrak{su}(3-\delta_\epsilon,1+\delta_\epsilon)$,
\begin{align}\label{trace-free}
\im\rho+\im\tau+\lambda-\overline{\lambda}=0.
\end{align} 

If two Hermitian frames $\underline{\vv}$ and $\tilde{\underline{\vv}}$ share the same $\vv_0=\tilde{\vv}_0$, then they differ by a transformation
\begin{equation}\label{P0}
\begin{aligned}
&\left[\begin{array}{cccc}\vv_0&\tilde{\vv}_1&\tilde{\vv}_2&\tilde{\vv}_3\end{array}\right]=\left[\begin{array}{cccc}\vv_0&\vv_1&\vv_2&\vv_3\end{array}\right]p_0,\\
&
p_0=\left[\begin{array}{cccr}
1&-\im(a_1\overline{c}_1+\epsilon a_2\overline{c}_2)
&\epsilon\im\e^{\im r}(\overline{a}_2\overline{c}_1-\overline{a}_1\overline{c}_2)
&c_0(\pm1)^{\delta_\epsilon} \\
0&a_1&-\epsilon\e^{\im r}\overline{a}_2&c_1(\pm1)^{\delta_\epsilon}\\
0&a_2&\e^{\im r}\overline{a}_1&c_2(\pm1)^{\delta_\epsilon}\\
0&0&0&(\pm1)^{\delta_\epsilon}\end{array}\right]
\begin{array}{c}
a_1,a_2,c_1,c_2\in\mathbb{C},\\
|a_1|^2+\epsilon|a_2|^2=(\pm1)^{\delta_\epsilon},\\
c_0=t-\tfrac{\im}{2}(|c_1|^2+\epsilon|c_2|^2),\\
r,t\in\mathbb{R}.\end{array}
\end{aligned}
\end{equation}
(See Remark \ref{deltaepsilonrem1} regarding the sign ambiguity $(\pm1)^{\delta_\epsilon}$ when $\epsilon=-1$.) All such $p_0\in U(3-\delta_\epsilon,1+\delta_\epsilon)$ form a subgroup that we call $\mathcal{P}_0$. Apparently, $\det(p_0)=\e^{\im r}$.

More generally, if the first basis vectors of two Hermitian frames $\underline{\vv}$ and $\tilde{\underline{\vv}}$ are $\mathbb{C}$-collinear -- i.e., $\tilde{\vv}_0=l\vv_0$ for some $l\in\mathbb{C}\setminus\{0\}$ -- then
\begin{align}\label{Pl}
&\tilde{\underline{\vv}}=\underline{\vv} p_lp_0;
&p_l=\left[\begin{array}{cccc}l&0&0&0\\0&1&0&0\\0&0&1&0\\0&0&0&\overline{l}^{-1}\end{array}\right],
\quad l\in\mathbb{C}\setminus\{0\},
&&p_0 \text{ as in \eqref{P0}.}
\end{align}
All such $p_l\in U(3-\delta_\epsilon,1+\delta_\epsilon)$ form an abelian subgroup that we call $\mathcal{P}_l$. We also name 
\begin{align}\label{PPlP0}
&\hat{\mathcal{P}}=\mathcal{P}_l\mathcal{P}_0,
&\mathcal{P}=\hat{\mathcal{P}}\cap SU(3-\delta_\epsilon,1+\delta_\epsilon);
\end{align}
the latter consists of multiples $p_lp_0$ with the $r$-coordinate of \eqref{P0} constrained by $\e^{\im r}=\overline{l}l^{-1}$, and it is exactly the parabolic subgroup $\mathcal{P}\subset SU(3-\delta_\epsilon,1+\delta_\epsilon)$ mentioned in the Tanaka-Chern-Moser Classification for CR-dimension $n=2$.

We conclude this section by recording the dual version of the $\hat{\mathcal{P}}$-transformation of Hermitian frames; namely, if we use superscripts to denote the dual coframe of $\underline{\vv}$ (i.e., $\vv^i(\vv_j)=\delta^i_j$) then the dual coframes of two Hermitian frames $\tilde{\underline{\vv}}$ and $\underline{\vv}$ with $\tilde{\vv}_0=l\vv_0$ are related by
\begin{align}\label{coframeP0Pl}
&\left[\begin{array}{c}\tilde{\vv}^0\\\tilde{\vv}^1\\\tilde{\vv}^2\\\tilde{\vv}^3\end{array}\right]=
{p_0}^{-1}{p_l}^{-1}
\left[\begin{array}{c}\vv^0\\\vv^1\\\vv^2\\\vv^3\end{array}\right];
&p_lp_0\in\hat{\mathcal{P}}.
\end{align}

\vspace{\baselineskip}

\subsection{First Adaptations}\label{firstadaptsec}

Let $\N\subset\mathbb{C}^4$ be the null-cone of $\h$,
\begin{align*}
&\N=\{\vv\in\mathbb{C}^4\setminus\{\left[\begin{smallmatrix}0&0&0&0\end{smallmatrix}\right]^t\}:\h(\vv,\vv)=0\}
&\Longrightarrow
&&T_{\vv}\N=\{{\tt w}\in\mathbb{C}^4:\Re(\h(\vv,{\tt w}))=0\}.
\end{align*}
We recognize the following distinguished subbundles of $T\N$ by their fibers,
\begin{align*}
&L_{\vv}=\{l\vv:l\in \mathbb{C}\},
&\subset
&&L_{\vv}^\bot=\{{\tt w}\in\mathbb{C}^4:\h(\vv,{\tt w})=0\}.
\end{align*}
Define the projection
\begin{align*}
h_\N:\hat{\mathcal{H}}&\to\N\\
\underline{\vv}&\mapsto\vv_0,
\end{align*}
which identifies $\N$ with the homogeneous space $U(3-\delta_\epsilon,1+\delta_\epsilon)/\mathcal{P}_0$. For each $\underline{\vv}=(\vv_0,\vv_1,\vv_2,\vv_3)\in h_{\N}^{-1}(\vv_0)$, we have spanning sets
\begin{align}\label{TNframed}
&\langle \vv_0\rangle_\mathbb{C}=L_{\vv_0},
&\langle \vv_0,\vv_1,\vv_2\rangle_\mathbb{C}=L_{\vv_0}^\bot,
&&\langle \vv_0,\vv_1,\vv_2,\im \vv_0,\im \vv_1,\im \vv_2,\vv_3\rangle_\mathbb{R}=T_{\vv_0}\N.
\end{align}
Conversely, one can assign an adapted basis of $T_{\vv_0}\N$ to each $\vv_0\in\N$ in order to define a section $s:\N\to\hat{\mathcal{H}}$ as follows: take $\vv_0$ itself to span $L_{\vv_0}$, choose two orthonormal vectors $\vv_1,\vv_2\in L_{\vv_0}^\bot$ (varying smoothly in a neighborhood of $\vv_0$), and then $\vv_3$ is uniquely determined to complete the Hermitian frame $s(\vv_0)=(\vv_0,\vv_1,\vv_2,\vv_3)$. 

\begin{rem}\label{deltaepsilonrem2}
Remark \ref{deltaepsilonrem1} explains why the $(\pm1)^{\delta_\epsilon}$ ambiguity in \eqref{hvivj} and \eqref{P0} allows us to choose any non-null $\vv_1\in L_{\vv_0}^\bot$ after $\vv_0$ is fixed.
The order of the basis elements in a frame reflects the ascending filtration $L\subset L^\bot\subset T\N$, but of course it is possible to first choose $\vv_3$ with $\h(\vv_0,\vv_3)=\im(\pm1)^{\delta_\epsilon}$  and then take orthonormal $\vv_1,\vv_2$ in the orthogonal complement of $L_{\vv_0}\oplus\langle\vv_3\rangle_\mathbb{R}$.
\end{rem}

With a section $s:\N\to\hat{\mathcal{H}}$, we pull back $\dd\vv_0\in\Omega^1(\hat{\mathcal{H}},\mathbb{C}^4)$ from \eqref{dv} to get
\begin{align}\label{dv0_pullback}
s^*\dd \vv_0=s^*\lambda \vv_0+s^*\eta \vv_1+s^*\zeta \vv_2+s^*\kappa \vv_3,
\end{align}
but $s^*\vv_0$ is just the identity map on $\N$, so its differential is the identity on $T\N$. Thus,  
\begin{align}\label{taut_pullback}
&s^*\lambda=\vv^0,
&s^*\eta=\vv^1,
&&s^*\zeta=\vv^2,
&&s^*\kappa=\vv^3.
\end{align}
Citing \eqref{coframeP0Pl} with $l=1$, we can say in other words that $\kappa\in\Omega^1(\hat{\mathcal{H}})$ and $\eta,\zeta,\lambda\in\Omega^1(\hat{\mathcal{H}},\mathbb{C})$ are semi-basic, tautological 1-forms for the (co)frame bundle fibration 
\begin{align}\label{hatHNfib}
\mathcal{P}_0\hookrightarrow U(3-\delta_\epsilon,1+\delta_\epsilon)\stackrel{h_\N}{\longrightarrow}\N.
\end{align}

Write $\mathbb{P}:\mathbb{C}^4\setminus\{\left[\begin{smallmatrix}0&0&0&0\end{smallmatrix}\right]^t\}\to\mathbb{CP}^3$ for the canonical complex projection, so the real hyperquadric is the image $\Q=\mathbb{P}(\N)$, and $T_{\mathbb{P}(\vv_0)}\Q=\mathbb{P}_*(T_{\vv_0}\N)\cong L_{\vv_0}^*\otimes T_{\vv_0}\N/L_{\vv_0}$ for $\vv_0\in\N$. The latter can be made explicit for a frame $\underline{\vv}\in h_\N^{-1}(\vv_0)$, namely
\begin{align}\label{TQ}
T_{\mathbb{P}(\vv_0)}\Q=\langle\mathbb{P}_*\vv_1,\mathbb{P}_*\vv_2\rangle_\mathbb{C}\oplus\langle\mathbb{P}_*\vv_3\rangle_\mathbb{R}\cong\langle\vv^0\otimes\vv_1,\vv^0\otimes\vv_2\rangle_\mathbb{C}\oplus\langle\vv^0\otimes\vv_3\rangle_\mathbb{R}.
\end{align}

\begin{rem}\label{deltaepsilonrem3}
By \eqref{Pl} and \eqref{coframeP0Pl} we see that for $\tilde{\underline{\vv}}=\underline{\vv}p_l\in h_\N^{-1}(\mathbb{P}^{-1}(\vv_0))$, $\tilde{\vv}^0\otimes\tilde{\vv}_3=|l|^{-2}\vv^0\otimes\vv_3$. When $\epsilon=-1$, a $p_0$-transformation \eqref{P0} with $(\pm1)^{\delta_\epsilon}=-1$ changes the sign of $\vv^0\otimes\vv_3$ (see Remark \ref{deltaepsilonrem2}).
\end{rem}

It's clear from \eqref{TNframed} that the fibers of $T\N\to \N$ vary along those of $\N\to\Q$. Nonetheless, Remark \ref{deltaepsilonrem3} reassures us that \eqref{TQ} is well-defined. Moreover, it always holds that $L_{\vv_0}=L_{\tilde{\vv}_0}$ and $L_{\vv_0}^\bot=L_{\tilde{\vv}_0}^\bot$ when $\mathbb{P}(\vv_0)=\mathbb{P}(\tilde{\vv}_0)$. In particular, for any $\vv_0\in \N$, $L_{\vv_0}=\ker\mathbb{P}_*|_{T_{\vv_0}\N}$, and there is a well-defined, $\mathbb{R}$-corank-1 subbundle
\begin{align}\label{QD}
D_{\mathbb{P}(\vv_0)}=\mathbb{P}_*L_{\vv_0}^\bot\subset T_{\mathbb{P}(\vv_0)}\Q.
\end{align}
Scalar multiplication by $\im$ in $\mathbb{C}^4$ defines an endomorphism $J:L_{\vv_0}^\bot\to L_{\vv_0}^\bot$ satisfying $J^2=-\mathbbm{1}$, of which $L_{\vv_0}$ is an invariant subspace, hence $J$ descends to a well-defined almost-complex structure $J:D\to D$.

The Tanaka-Chern-Moser Classification for $n=2$ ensures that $\mathcal{H}$ realizes the Cartan bundle of the CR structure $(\Q,D,J)$, as encoded in the fibration
\begin{align}\label{CBQ}
&\mathcal{P}\hookrightarrow SU(3-\delta_\epsilon,1+\delta_\epsilon)\stackrel{h_\Q}{\longrightarrow}\Q;
&h_\Q=\mathbb{P}\circ h_\N.
\end{align}
It will be convenient to work instead with $\hat{\mathcal{P}}\to\hat{\mathcal{H}}\stackrel{h_\Q}{\longrightarrow}\Q$ in order to factor through \eqref{hatHNfib}:
\begin{displaymath}
\begin{tikzcd}        
    \mathcal{P}_0 \arrow[hook]{r}{\underline{\vv}\mapsto\underline{\vv}p_0}
    	& U(3-\delta_\epsilon,1+\delta_\epsilon)\arrow{d}{h_\N}\\
     \mathcal{P}_l\cong\mathbb{C}\setminus\{0\}\ar{r}{\vv_0\mapsto l\vv_0}
     	&\N\arrow{d}{\mathbb{P}}\\
     	&\Q
\end{tikzcd}
\end{displaymath}
To this end, let $\varsigma:\Q\to\hat{\mathcal{H}}$ be a section that factors through $s:\N\to\hat{\mathcal{H}}$ as in \eqref{taut_pullback}. Then $\varsigma^*\kappa\in\Omega^1(\Q)$ is a characteristic form annihilating \eqref{QD}, which along with $\varsigma^*\eta,\varsigma^*\zeta\in\Omega^1(\Q,\mathbb{C})$ and their conjugates furnishes a 1-adapted coframing whose Levi form \eqref{levimatrix} is represented 
\begin{align}\label{Qlevi}
\ell=\left[\begin{array}{cc}1&0\\0&\epsilon\end{array}\right]
\end{align}
in the basis \eqref{TQ} given by $\underline{\vv}=\varsigma(\mathbb{P}(\vv_0))$. Remark \ref{deltaepsilonrem3} explains why the representation \eqref{Qlevi} of $\ell$ does not depend on the choice of $(\pm1)^{\delta_\epsilon}$ in \eqref{hvivj}, \eqref{P0} when $\epsilon=-1$. 
%In other words $\ell(\mathbb{P}_*\vv_1,\mathbb{P}_*\vv_1)=1$ regardless of $\h(\vv_1,\vv_1)=(\pm1)^{\delta_\epsilon}$.

If $M$ is a 3-dimensional manifold CR-embedded in $\Q$, then there is a rank-2, $J$-invariant subbundle $D_M\subset D$ tangent to $M$. Let $\hat{M}=\mathbb{P}^{-1}(M)\subset\N$ be the cone over $M$ and $L_M\subset T\hat{M}$ be $\mathbb{P}_*^{-1}(D_M)$. Name the restriction and projections
\begin{align*}
&\mathcal{H}^0=h_{\Q}^{-1}(M)\subset\hat{\mathcal{H}};
&h_{\hat{M}}=h_\N|_{\mathcal{H}^0},
&&h_M=h_\Q|_{\mathcal{H}^0}.
\end{align*}
The Hermitian frames $\mathcal{H}^0$ are ``zero-adapted to $M$" and the fibers of $h_M$ have (real) dimension $\dim\hat{\mathcal{P}}=11$ -- see \eqref{PPlP0}. We will exploit the degrees of freedom \eqref{P0} and \eqref{Pl} in the structure group $\hat{\mathcal{P}}$ of $\mathcal{H}^0$ to reduce to subbundles of frames which are increasingly adapted to the CR geometry of $M$.

Since $\Q$ is Levi-nondegenerate, $D\subset T\Q$ is a rank-4 contact distribution, so there cannot be a 3-dimensional submanifold of $\Q$ which is tangent to $D$. That is to say $\varsigma^*\kappa$ cannot vanish identically on $TM$ or, equivalently, $s^*\kappa|_{T\hat{M}}\neq0$. Therefore, we restrict to those Hermitian frames 
\begin{equation}\label{TMv3}
\underline{\vv}\in h_{\hat{M}}^{-1}(\vv_0) \quad\text{ with } \vv_3\in T_{\vv_0}\hat{M}.
\end{equation}
Further adaptation branches based on the (non)vanishing of \eqref{Qlevi} on $D_M$. If $M$ is Levi-nondegenerate, $\ell|_{D_M}\neq 0$ and we consider $\underline{\vv}\in\mathcal{H}^0$ with $L_M=\langle \vv_0,\vv_1\rangle_\mathbb{C}$, which along with \eqref{TMv3} reduces $\mathcal{H}^0$ and its structure group over $\hat{M}$:
\begin{align}\label{LNfirstadapt}
&\mathcal{H}^1=\left\{\underline{\vv}\in\mathcal{H}^0: \begin{array}{c}L_M=\langle \vv_0,\vv_1\rangle_\mathbb{C}\\ T\hat{M}=L_M\oplus\langle\vv_3\rangle_\mathbb{R}\end{array}\right\},
&\mathcal{P}_1=\left\{p_0\in\mathcal{P}_0 \text{ as in } \eqref{P0} : \begin{array}{r}a_2=c_2=0\\ (\pm1)^{\delta_\epsilon}=1\end{array}\right\}.
\end{align}
When $\epsilon=-1$ it is also possible that $\ell_{D_M}=0$. Such $M$ is Levi-flat and contains a complex curve tangent to the $\mathbb{P}_*$ image of an $\h$-null distribution of complex rank two in $T\hat{M}$. In this case we arrange for $L_M=\langle \vv_0,\vv_1+\vv_2\rangle_\mathbb{C}$ combined with \eqref{TMv3} to obtain 
\begin{align}\label{LFfirstadapt}
&\mathcal{H}^1=\left\{\underline{\vv}\in\mathcal{H}^0: \begin{array}{c}L_M=\langle \vv_0,\vv_1+\vv_2\rangle_\mathbb{C}\\ T\hat{M}=L_M\oplus\langle\vv_3\rangle_\mathbb{R}\end{array}\right\}
&\mathcal{P}_1=\left\{p_0\in\mathcal{P}_0 \text{ as in } \eqref{P0} : 
\begin{array}{c}\epsilon=-1, (\pm1)^{\delta_{-1}}=1\\a_1+a_2=\e^{\im r}(\overline{a}_2+\overline{a}_1)\\ c_2=c_1\end{array}\right\}.
\end{align}
In every case, we can say $\forall \underline{\vv},\tilde{\underline{\vv}}\in h_{\hat{M}}^{-1}(\vv_0)\cap\mathcal{H}^1\ \exists p_1\in\mathcal{P}_1 \text{ such that } \tilde{\underline{\vv}}=\underline{\vv}p_1$, so that the fibers of $h_M|_{\mathcal{H}^1}$ have real dimension $\dim\mathcal{P}_l\mathcal{P}_1=7$. Furthermore, our adaptations in both cases are undisturbed by the projection \eqref{orient}, so we may descend as in \eqref{CBQ} to oriented frames over $M$
\begin{align}\label{H2orient}
&\mathcal{P}_2\to\mathcal{H}^2\stackrel{h_M}{\longrightarrow}M,
&\mathcal{P}_2=\mathcal{P}_l\mathcal{P}_1\cap\mathcal{P},
&&\mathcal{H}^2=\mathcal{H}^1\cap\mathcal{H}.
\end{align}

\vspace{\baselineskip}

\subsection{Second Fundamental Form $(\II)$} \label{IIsec}
All higher-order adaptations of moving frames over $M\subset\Q$ will be controlled by the MC equations \eqref{UMCeq}. Specifically, we pull back the MC form $\mu$ as in \eqref{dv} along the inclusion $\mathcal{H}^1\hookrightarrow\mathcal{H}^0$ and explore the differential consequences of the algebraic relations imposed on the individual 1-forms \eqref{UMCform} over $\mathcal{H}^1$. Obviously we must separately consider 3-folds $M$ which are Levi-nondegenerate \eqref{LNfirstadapt} or Levi-flat \eqref{LFfirstadapt}, but in both cases the reduction \eqref{H2orient} will force the relation \eqref{trace-free} (our notation suppresses the pullback along the inclusion $\mathcal{H}^2\hookrightarrow\hat{\mathcal{H}}$ and keeps the same names for the MC forms with every such reduction). Moreover, in both cases the method of moving frames exploits the remaining degrees of freedom in $\mathcal{P}_2$ to normalize a symmetric, $\mathbb{C}$-bilinear tensor
\begin{align}\label{II}
\II:{\bigodot}^2(T\mathcal{H}^2/T\mathcal{P}_2)\to \mathbb{C},
\end{align}
which we dub the \emph{Second Fundamental Form} by analogy to the study of hypersurfaces in Euclidean space. (One might consider the the Levi form $\ell$ of $\Q$ restricted to $D_M$ to be a first fundamental form of sorts, but for $\dim M=3$ this is simply a real function whose only import is its non-vanishing.)

\vspace{\baselineskip}

\subsubsection{Levi-Nondegenerate 3-folds}\label{LNIIsec}  Our adaptation \eqref{LNfirstadapt} together with \eqref{dv0_pullback} implies that 
\begin{align}\label{zetazero}
\zeta=0  \text{ on }\mathcal{H}^2.
\end{align}
Applying Cartan's Lemma to the equation \eqref{UMCeq} for $\dd\zeta$ yields
\begin{align}\label{IIcoeff}
&\left[\begin{array}{c}\phi_1\\\phi_2\end{array}\right]=
\left[\begin{array}{cc}a&b\\b&c\end{array}\right]\left[\begin{array}{c}\eta\\\kappa\end{array}\right]
&\text{ for some } a,b,c\in C^\infty(\mathcal{H}^2,\mathbb{C}).
\end{align}

\begin{rem}\label{rrecycle}
It behooves us to remember that $\mathcal{P}_2$ \eqref{H2orient} is in part defined by $a_2=0$ (which forces $|a_1|=1$) in \eqref{P0}, as well as $\e^{\im r}=\overline{l}l^{-1}$ in \eqref{P0}, \eqref{Pl}. Consequently, there is no crime in recycling notation by setting $a_1=\e^{\im r}$ in \eqref{P0} so that $\im\rho$ in \eqref{UMCform} measures the infinitesimal generator of this $U(1)$-action on $\mathcal{H}^2$. 
\end{rem}

Differentiating \eqref{IIcoeff} using $\dd\phi_1,\dd\phi_2$ from \eqref{UMCeq} reveals
\begin{equation}\label{dII}
\begin{aligned}
\dd\left[\begin{array}{c}a\\b\\c\end{array}\right]=
\left[\begin{array}{ccc}3\im\rho-\overline{\lambda}&0&0\\\xi&2\im\rho-2\overline{\lambda}&0\\0&2\xi&\im\rho-3\overline{\lambda}\end{array}\right]\left[\begin{array}{c}a\\b\\c\end{array}\right]+
\left[\begin{array}{ccc}u_1&u_2&2\im b\\u_2&u_3&\im c\\u_3&u_4&0\end{array}\right]
\left[\begin{array}{c}\eta\\\kappa\\\overline{\eta}\end{array}\right]\\
\text{ for some } u_1,u_2,u_3,u_4\in C^\infty(\mathcal{H}^2,\mathbb{C}).
\end{aligned}
\end{equation}
The identities \eqref{dII} imply that if one of $a,b,c$ vanishes identically on $\mathcal{H}^2$, they all vanish identically. Furthermore,
\begin{align}\label{detII}
\dd(ac-b^2)\equiv4(ac-b^2)(\im\rho-\overline{\lambda})\mod\{\kappa,\eta\},
\end{align}
so the determinant of the matrix \eqref{IIcoeff} is either identically zero or non-vanishing on each fiber of $\mathcal{H}^2\to M$.

\begin{definition}
The second fundamental form \eqref{II} of a Levi-nondegenerate 3-fold $M$ CR embedded in $\Q$ is given by
\begin{equation}\label{LNII}
\II=a\eta\odot\eta+2b\eta\odot\kappa+c\kappa\odot\kappa,
\end{equation}
where the coefficients are derived from \eqref{zetazero} via \eqref{IIcoeff}. The condition that $\II$ is of (sub)maximal rank on a fiber of $\mathcal{H}^2\to M$ is invariant under the action of CR symmetry group $SU(3-\delta_\epsilon,1+\delta_\epsilon)$ on $\Q$.
\end{definition}

We update the remaining MC equations over $\mathcal{H}^2$,
\begin{equation}\label{LNH2MCeq}
\begin{aligned}
\dd\kappa&=\im\eta\wedge\overline{\eta}+(\lambda+\overline{\lambda})\wedge\kappa,\\
\dd\eta&=(\lambda-\im\rho)\wedge\eta-\xi\wedge\kappa,\\
\dd\psi&=\psi\wedge(\lambda+\overline{\lambda})-\im\xi\wedge\overline{\xi}+\epsilon\im\kappa\wedge(b\overline{c}\eta-\overline{b}c\overline{\eta})-\epsilon\im|b|^2\eta\wedge\overline{\eta},\\
\dd\xi&=\psi\wedge\eta+\xi\wedge(\overline{\lambda}+\im\rho)+\epsilon\kappa\wedge(|b|^2\eta-ac\overline{\eta})-\epsilon\overline{a}b\eta\wedge\overline{\eta},\\
\dd\lambda&=\im\overline{\xi}\wedge\eta-\psi\wedge\kappa,\\
\dd\rho&=-\overline{\xi}\wedge\eta-\xi\wedge\overline{\eta}-\epsilon\im\kappa\wedge(a\overline{b}\eta+\overline{a}b\overline{\eta})+\epsilon\im|a|^2\eta\wedge\overline{\eta}.
\end{aligned}
\end{equation}
Recall from Remark \ref{MQbundlesrem} that a CR embedding $M\hookrightarrow\Q$ lifts to a mapping $\F^3\hookrightarrow\mathcal{H}$ between the Cartan bundles of $M$ and $\Q$ in a manner that pulls back $\mu$ \eqref{UMCform} to $\gamma$ \eqref{Mcarcon}. To realize the image of $\F^3$ over $M\subset\Q$, we locate $\gamma$ within $\mu|_{\mathcal{H}^2}$ using the fact that the semi-basic forms $\kappa,\eta,\overline{\eta}$ over both $\F^3$ and $\mathcal{H}^2$ encode adapted (co)framings of $M$. The equations \eqref{threeCRSE} and \eqref{LNH2MCeq} for $\dd\kappa$, $\dd\eta$ alone demonstrate 
\begin{align*}
&\lambda=-\alpha+\im\rho+a_0\kappa,
&\xi=\beta-b_0\kappa-a_0\eta,
&&\text{for some }a_0,b_0\in C^\infty(\mathcal{H}^2,\mathbb{C}),
\end{align*}
and comparing $\dd\alpha,\dd\beta$ in \eqref{threeCRSE} to $\dd\lambda,\dd\xi$ in \eqref{LNH2MCeq} validates the substitution
\begin{align*}
&\psi=-\sigma+s_0\kappa+s_1\eta+\overline{s}_1\overline{\eta},
&&\text{for some }s_0\in C^\infty(\mathcal{H}^2), s_1\in C^\infty(\mathcal{H}^2,\mathbb{C}),
\end{align*}
after which we solve for $a_0,b_0,s_0,s_1$ to bring the first five equations of \eqref{LNH2MCeq} into the form \eqref{threeCRSE}:
\begin{align*}
a_0&=-\epsilon\tfrac{\im}{4}|a|^2,
&&b_0=-\epsilon(\im\overline{a}b+\tfrac{1}{6}a\overline{u}_1),\\
s_0&=\tfrac{9}{16}|a|^4+\epsilon(\tfrac{3\im}{4}(a\overline{u}_2-\overline{a}u_2)-3|b|^2-\tfrac{1}{6}|u_1|^2)
&&s_1=-\epsilon\tfrac{1}{12}(6a\overline{b}+\im\overline{a}u_1).
\end{align*}
In conclusion, 
\begin{equation}\label{threeCRSEinH}
\begin{aligned}
\alpha&=\im\rho-\lambda-\epsilon\tfrac{\im}{4}|a|^2\kappa,\quad \quad 
\beta=\xi-\epsilon(\im\overline{a}b+\tfrac{1}{6}a\overline{u}_1)\kappa-\epsilon\tfrac{\im}{4}|a|^2\eta,\\
\sigma&=-\psi+(\tfrac{9}{16}|a|^4+\epsilon(\tfrac{3\im}{4}(a\overline{u}_2-\overline{a}u_2)-3|b|^2-\tfrac{1}{6}|u_1|^2))\kappa
-\epsilon\tfrac{1}{12}(6a\overline{b}+\im\overline{a}u_1)\eta
-\epsilon\tfrac{1}{12}(6\overline{a}b-\im a\overline{u}_1)\overline{\eta}.
\end{aligned}
\end{equation}

\vspace{\baselineskip}

Next we seek expressions for the coefficient functions $S,P$ of the curvature tensor $\dd\gamma+\gamma\wedge\gamma$ as they appear in the equations for $\dd\beta,\dd\sigma$, as well as the higher-order coefficients in \eqref{Bianchi}, \eqref{dRQUVW}, and \eqref{dR'}. To this end, we differentiate \eqref{dII} and use the identity $\dd^2=0$ to compute
\begin{equation}\label{dus}
\begin{aligned}
\dd\left[\begin{array}{c}u_1\\u_2\\u_3\\u_4\end{array}\right]&=
\left[\begin{array}{cccc}4\im\rho-\overline{\lambda}-\lambda&0&0&0\\\xi&3\im\rho-2\overline{\lambda}-\lambda&0&0\\0&2\xi&2\im\rho-3\overline{\lambda}-\lambda&0\\
0&0&3\xi&\im\rho-4\overline{\lambda}-\lambda\end{array}\right]
\left[\begin{array}{c}u_1\\u_2\\u_3\\u_4\end{array}\right]\\
&+\left[\begin{array}{cc}0&3\im a\\-a&2\im b\\-2b&\im c\\-3c&0\end{array}\right]
\left[\begin{array}{c}\psi\\\overline{\xi}\end{array}\right]
+\left(\left[\begin{array}{cccc}v_1&v_2&3\im u_2\\v_2&v_3&2\im u_3\\v_3&v_4&\im u_4\\v_4&v_5&0\end{array}\right]
-\epsilon\left[\begin{array}{cc}0&3a\\a&2b\\2b&c\\3c&0\end{array}\right]
\left[\begin{array}{ccc}0&|b|^2&\overline{a}b\\0&a\overline{b}&|a|^2\end{array}\right]\right)
\left[\begin{array}{c}\eta\\\kappa\\\overline{\eta}\end{array}\right],
\end{aligned}
\end{equation}
for some $v_1,\dots,v_5\in C^\infty(\mathcal{H}^2,\mathbb{C})$. This is sufficient to express
\begin{equation}\label{SPR}
\begin{aligned}
\left[\begin{array}{c}S\\P\\R\end{array}\right]&=
\frac{\epsilon}{6}\left[\begin{array}{ccc}
a&b&c\\
u_1&u_2&u_3\\
v_1&v_2&v_3
\end{array}\right]
\left[\begin{array}{c}
\overline{v}_1\\8\im\overline{u}_1\\-12\overline{a}
\end{array}\right]
+\frac{\epsilon}{6}\left[\begin{array}{ccc}
0&0&0\\
a& b& c\\
2u_1&2u_2&2u_3
\end{array}\right]
\left[\begin{array}{c}
-4\im\overline{v}_2\\24\overline{u}_2\\24\im\overline{b}
\end{array}\right]\\
&+\frac{\epsilon}{6}\left[\begin{array}{ccc}
0&0&0\\
0&0&0\\
a&b&c
\end{array}\right]
\left[\begin{array}{c}
-12\overline{v}_3\\-48\im\overline{u}_3\\24\overline{c}
\end{array}\right]
-\frac{1}{6}\left[\begin{array}{ccc}
0&0&0\\
0&30b&0\\
30|a|^2&54u_2&36u_1
\end{array}\right]
\left[\begin{array}{c}-\epsilon|a|^4\\\im\overline{a}^2a\\\im\overline{a}^2b
\end{array}\right]\\
&+\frac{1}{12}\left[\begin{array}{cc}
0&0\\
-20\overline{a}&0\\
72\im\overline{b}&108\im\overline{a}
\end{array}\right]
\left[\begin{array}{c}a^2\overline{u}_1\\a^2\overline{u}_2\end{array}\right]
-\left[\begin{array}{cc}
0&0\\
0&0\\
36b&5u_1
\end{array}\right]
\left[\begin{array}{l}|a|^2\overline{b}\\|a|^2\overline{u}_1\end{array}\right].
\end{aligned}
\end{equation}

\begin{rem}\label{Srem}
The prototypical example of a 3-dimensional CR manifold is a hypersurface given by a regular level set of a smooth, non-constant function $\varrho:\mathbb{C}^2\to\mathbb{R}$. The curvature coefficient $S$ for the level set depends on derivatives of $\varrho$ up to order six, which is to say that $S$ is a function of the 6-jet of $\varrho$. In the present setting, $S$ for $M\subset\Q$ is a function of the 2-jet of $\II$.
\end{rem}

The remaining functions in \eqref{Bianchi} require another derivative of $\II$, so we apply $\dd^2=0$ to \eqref{dus} to get 
\begin{align}\label{dvs}
\dd\left[\begin{array}{c}v_1\\v_2\\v_3\\v_4\\v_5\end{array}\right]&=
\left[\begin{array}{ccccc}5\im\rho-\overline{\lambda}-2\lambda&0&0&0&0\\\xi&4\im\rho-2\overline{\lambda}-2\lambda&0&0&0\\0&2\xi&3\im\rho-3\overline{\lambda}-2\lambda&0&0\\
0&0&3\xi&2\im\rho-4\overline{\lambda}-2\lambda&0\\0&0&0&4\xi&\im\rho-5\overline{\lambda}-2\lambda\end{array}\right]
\left[\begin{array}{c}v_1\\v_2\\v_3\\v_4\\v_5\end{array}\right]\\
&+\left[\begin{array}{cc}0&8\im u_1\\-2u_1&6\im u_2\\-4u_2&4\im u_3\\-6u_3&2\im u_4\\-8u_4&0\end{array}\right]
\left[\begin{array}{c}\psi\\\overline{\xi}\end{array}\right]
+\left[\begin{array}{ccc}w_1&w_2&4\im v_2\\w_2&w_3&3\im v_3\\w_3&w_4&2\im v_4\\w_4&w_5&\im v_5\\w_5&w_6&0\end{array}\right]
\left[\begin{array}{c}\eta\\\kappa\\\overline{\eta}\end{array}\right]\nonumber\\
&-\epsilon\left(\left[\begin{array}{ccc}0&0&10u_1\\0&4u_1&6u_2\\u_1&6u_2&3u_3\\3u_2&6u_3&u_4\\6u_3&4u_4&0\end{array}\right]
\left[\begin{array}{ccc}0&\overline{b}c&\overline{a}c\\0&|b|^2&\overline{a}b\\0&a\overline{b}&|a|^2\end{array}\right]
-\im\left[\begin{array}{ccc}0&0&6a\\0&3a&3b\\a&4b&c\\3b&3c&0\\6c&0&0\end{array}\right]
\left[\begin{array}{ccc}0&|c|^2&\overline{b}c\\0&b\overline{c}&|b|^2\\0&a\overline{c}&a\overline{b}\end{array}\right]\right)
\left[\begin{array}{c}\eta\\\kappa\\\overline{\eta}\end{array}\right],\nonumber
\end{align}
for some $w_1,\dots,w_6\in C^\infty(\mathcal{H}^2,\mathbb{C})$, and with this we can write
\begin{align}\label{QUVWR}
\left[\begin{array}{c}Q\\U\\V\\W\\R_1'\end{array}\right]&=
\frac{\epsilon}{6}\left[\begin{smallmatrix}
a&b&c\\
0&0&0\\
u_1&u_2&u_3\\
0&0&0\\
v_1&v_2&v_3
\end{smallmatrix}\right]
\left[\begin{array}{c}
\overline{w}_1\\10\im\overline{v}_1\\-20\overline{u}_1
\end{array}\right]
-\frac{\epsilon}{6}\left[\begin{smallmatrix}
0&0&0\\
-a&-b&-c\\
4\im a&4\im b&4\im c\\
-u_1&-u_2&-u_3\\
8\im u_1&8\im u_2&8\im u_3
\end{smallmatrix}\right]
\left[\begin{array}{c}
\overline{w}_2\\8\im\overline{v}_2\\-12\overline{u}_2\end{array}\right]
-\frac{\epsilon}{6}\left[\begin{smallmatrix}
0&0&0\\
0&0&0\\
0&0&0\\
-a& -b& -c\\
3\im a&3\im b&3\im c
\end{smallmatrix}\right]
\left[\begin{array}{c}
-4\im\overline{w}_3\\24\overline{v}_3\\24\im\overline{u}_3\end{array}\right]\\
&+\frac{\epsilon}{6}\left[\begin{smallmatrix}
0&0&0&0&0&0\\
u_2&u_3&u_4&0&0&0\\
\im u_2&\im u_3&\im u_4&0&0&0\\
v_2&v_3&v_4&u_2&u_3&u_4\\
2\im v_2&2\im v_3&2\im v_4&2\im u_2&2\im u_3&2\im u_4
\end{smallmatrix}\right]
\left[\begin{smallmatrix}
\overline{v}_1\\8\im\overline{u}_1\\-12\overline{a}\\
-4\im\overline{v}_2\\24\overline{u}_2\\24\im\overline{b}\\
\end{smallmatrix}\right]
-\frac{1}{6}\left[\begin{smallmatrix}
0&0&0&0\\
0&0&24b&0\\
0&\im45 c&\im 84 b&0\\
\tfrac{1}{2}(45\overline{a}b+5\im a\overline{u}_1)&27u_3&54u_2&4u_1\\
450\im\overline{a}b+180a\overline{u}_1&144\im u_3&288\im u_2&48\im u_1
\end{smallmatrix}\right]
\left[\begin{array}{c}-\epsilon|a|^4\\\im\overline{a}^2a\\\im\overline{a}^2b\\\im\overline{a}^2c\end{array}\right]\nonumber\\
&-\frac{1}{12}\left[\begin{smallmatrix}
0&0&0\\
0&\im \overline{a}&0\\
20\overline{u}_1&\tfrac{53}{2}\overline{a}&0\\
20\overline{u}_2&8\overline{b}&22\overline{a}\\
-180\im\overline{u}_2&-84\im\overline{b}&-156\im\overline{a}
\end{smallmatrix}\right]
\left[\begin{smallmatrix}a^2\overline{u}_1\\a^2\overline{v}_1\\a^2\overline{v}_2\end{smallmatrix}\right]
-\frac{1}{3}\left[\begin{smallmatrix}
0&0&0&0&0&0&0&0\\
0&3b&0&0&0&0&0&0\\
0&63\im b&0&0&0&0&0&0\\
-3c&9u_2&21\im b&\tfrac{5}{24}\im u_1&24\overline{a}&8\im a&0&\tfrac{13}{3}\overline{a}\\
144\im c&123\im u_2&288b&20u_1&288\im\overline{a}&204a&15\overline{u}_1&82\im\overline{a}
\end{smallmatrix}\right]
\left[\begin{smallmatrix}|a|^2\overline{b}\\|a|^2\overline{u}_1\\|a|^2\overline{u}_2\\|a|^2\overline{v}_1\\|b|^2b\\|b|^2\overline{u}_1\\|u_1|^2a\\|u_1|^2b\end{smallmatrix}\right],\nonumber
\end{align}
as well as the totally real coefficients
\begin{align}\label{RRzero}
\left[\begin{array}{c}R_0'\\R_{0}''\end{array}\right]&=
\epsilon\left[\begin{smallmatrix}-2a&-2b&-2c\\5\im a&5\im b&5\im c\end{smallmatrix}\right]
\left[\begin{smallmatrix}\overline{w}_4\\ 4\im\overline{v}_4\\ -2\overline{u}_4\end{smallmatrix}\right]
-\frac{\epsilon}{3}\left[\begin{smallmatrix}-u_1&-u_2&-u_3\\5\im u_1&5\im u_2&5\im u_3\end{smallmatrix}\right]
\left[\begin{smallmatrix}-4\im\overline{w}_3\\ 24\overline{v}_3\\ 24\im\overline{u}_3\end{smallmatrix}\right]
-\frac{\epsilon}{12}\left[\begin{smallmatrix}-2v_1&-2v_2&-2v_3\\25\im v_1&25\im v_2&25\im v_3\end{smallmatrix}\right]
\left[\begin{smallmatrix}\overline{w}_2\\ 8\im\overline{v}_2\\ -12\overline{u}_2\end{smallmatrix}\right]\\
&-\epsilon\left[\begin{smallmatrix}2u_2&2u_3&2u_4\\5\im u_2&5\im u_3&5\im u_4\end{smallmatrix}\right]
\left[\begin{smallmatrix}\overline{v}_3\\ 4\im\overline{u}_3\\ -2\overline{c}\end{smallmatrix}\right]
+\frac{\epsilon}{6}\left[\begin{smallmatrix}2v_2&2v_3&2v_4\\5\im v_2&5\im v_3&5\im v_4\end{smallmatrix}\right]
\left[\begin{smallmatrix}-4\im\overline{v}_2\\ 24\overline{u}_2\\ 24\im\overline{b}\end{smallmatrix}\right]
+\frac{\epsilon}{12}\left[\begin{smallmatrix}2w_2&2w_3&2w_4\\5\im w_2&5\im w_3&5\im w_4\end{smallmatrix}\right]
\left[\begin{smallmatrix}\overline{v}_1\\ 8\im\overline{u}_1\\ -12\overline{a}\end{smallmatrix}\right]
\nonumber\\
&+\frac{\epsilon}{6}\left(\left[\begin{smallmatrix}0&0&0\\w_1&w_2&w_3\end{smallmatrix}\right]
\left[\begin{smallmatrix}\overline{w}_1\\ 10\im\overline{v}_1\\ -20\overline{u}_1\end{smallmatrix}\right]
+|a|^4\left[\begin{smallmatrix}0\\6210|b|^2+\tfrac{2195}{2}|u_1|^2\end{smallmatrix}\right]
\right)
-\frac{|a|^2}{6}\left[\begin{smallmatrix}0\\\tfrac{1035}{2}|a|^6 +342|c|^2+1710|u_2|^2+\tfrac{217}{4}|v_1|^2\end{smallmatrix}\right]
\nonumber\\
&+\Re\left(\frac{\epsilon\overline{a}^3a}{3}\left[\begin{smallmatrix}90&90\\
\tfrac{2385}{2}\im&1650\im\end{smallmatrix}\right]
\left[\begin{smallmatrix}au_2\\bu_1\end{smallmatrix}\right]
+\frac{\overline{a}^2}{3}\left[\begin{smallmatrix}-54\im& -60\im& -6\im& -36\im& -54\im\\ 270& 390& 51& 210& 315\end{smallmatrix}\right]
\left[\begin{smallmatrix}av_3\\ bv_2\\ cv_1\\ u_1u_3\\ {u_2}^2\end{smallmatrix}\right]
+\frac{|a|^2}{3}\left[\begin{smallmatrix}-180&-30\\1350\im&375\im\end{smallmatrix}\right]
\left[\begin{smallmatrix}b\overline{u}_3\\u_1\overline{v}_2\end{smallmatrix}\right]
\right)\nonumber\\
&-\Re\left(\frac{1}{3}\left[\begin{smallmatrix}(288|b|^2+30|u_1|^2)a\overline{u}_2
+a(\overline{b}(10u_1\overline{v}_1-24\im\overline{u}_1u_2)+36b\overline{c}\overline{u}_1)\\
(720b\overline{c}+199\im\overline{u}_1v_1)\overline{a}b-\im(2880|b|^2+425|u_1|^2)a\overline{u}_2
+40(30b\overline{u}_2+u_1\overline{v}_1)\overline{a}u_1-498\im ab\overline{c}\overline{u}_1
\end{smallmatrix}\right]
\right)\nonumber\\
&-\frac{1}{6}\left[\begin{smallmatrix}0\\30|u_1|^4+1440|b|^4+832|bu_1|^2\end{smallmatrix}\right].
\nonumber
\end{align}
One last derivative of \eqref{dvs} yields
\begin{align}\label{dws}
\dd\left[\begin{array}{c}w_1\\w_2\\w_3\\w_4\\w_5\\w_6\end{array}\right]&=
\left[\begin{smallmatrix}6\im\rho-\overline{\lambda}-3\lambda&0&0&0&0&0\\
\xi&5\im\rho-2\overline{\lambda}-3\lambda&0&0&0&0\\
0&2\xi&4\im\rho-3\overline{\lambda}-3\lambda&0&0&0\\
0&0&3\xi&3\im\rho-4\overline{\lambda}-3\lambda&0&0\\
0&0&0&4\xi&2\im\rho-5\overline{\lambda}-3\lambda&0\\
0&0&0&0&5\xi&\im\rho-6\overline{\lambda}-3\lambda\end{smallmatrix}\right]
\left[\begin{array}{c}w_1\\w_2\\w_3\\w_4\\w_5\\w_6\end{array}\right]\\
&+\left[\begin{smallmatrix}0&15\im v_1\\-3v_1&12\im v_2\\-6v_2&9\im v_3\\-9v_3&6\im v_4\\-12v_4&3\im v_5\\-15v_5&0\end{smallmatrix}\right]\left[\begin{array}{c}\psi\\\overline{\xi}\end{array}\right]
+\left(\left[\begin{smallmatrix}z_1&z_2&5\im w_2\\z_2&z_3&4\im w_3\\z_3&z_4&3\im w_4\\z_4&z_5&2\im w_5\\z_5&z_6&\im w_6\\z_6&z_7&0\end{smallmatrix}\right]
-\epsilon\left[\begin{smallmatrix}10{u_1}^2\\10u_1u_2\\4u_1u_3+6{u_2}^2\\u_1u_4+9u_2u_3\\4u_2u_4+6{u_3}^2\\10u_3u_4\end{smallmatrix}\right]\left[\begin{array}{ccc}0&\overline{b}&\overline{a}\end{array}\right]\right)\left[\begin{array}{c}\eta\\\kappa\\\overline{\eta}\end{array}\right]\nonumber\\
-&\epsilon\left(
\left[\begin{smallmatrix}0&0&15v_1\\0&5v_1&10v_2\\v_1&8v_2&6v_3\\3v_2&9v_3&3v_4\\6v_3&8v_4&v_5\\10v_4&5v_5&0\end{smallmatrix}\right]\left[\begin{array}{ccc}0&\overline{b}c&\overline{a}c\\0&|b|^2&\overline{a}b\\0&a\overline{b}&|a|^2\end{array}\right]
-\im \left[\begin{smallmatrix}0&0&30u_1\\0&12u_1&18u_2\\3u_1&18u_2&9u_3\\9u_2&18u_3&3u_4\\18u_3&12u_4&0\\30u_4&0&0\end{smallmatrix}\right]\left[\begin{array}{ccc}0&|c|^2&\overline{b}c\\0&b\overline{c}&|b|^2\\0&a\overline{c}&a\overline{b}\end{array}\right]\right)
\left[\begin{array}{c}\eta\\\kappa\\\overline{\eta}\end{array}\right],\nonumber
\end{align}
for some $z_1,\dots,z_7\in C^\infty(\mathcal{H}^2,\mathbb{C})$. In the language of Remark \ref{Srem}, we have functions of the 4-jet of $\II$,
\begin{equation}\label{QUUV}
\begin{aligned}
\left[\begin{array}{c}Q'\\U_1'\\U_2'\\V'\end{array}\right]&=
\frac{\epsilon}{6}\left(\left[\begin{smallmatrix}
a&b&c\\
0&0&0\\
0&0&0\\
0&0&0\end{smallmatrix}\right]
\left[\begin{smallmatrix}
\overline{z}_1\\ 12\im\overline{w}_1\\ -30\overline{v}_1\end{smallmatrix}\right]
+
\left[\begin{smallmatrix}
0&0&0\\
a&b&c\\
0&0&0\\
u_1&u_2&u_3\end{smallmatrix}\right]
\left[\begin{smallmatrix}
\overline{z}_2\\ 10\im\overline{w}_2\\ -20\overline{v}_2\end{smallmatrix}\right]
+
\left[\begin{smallmatrix}
0&0&0\\
0&0&0\\
a&b&c\\
-4\im a&-4\im b&-4\im c\end{smallmatrix}\right]
\left[\begin{smallmatrix}
\overline{z}_3\\ 8\im\overline{w}_3\\-12\overline{v}_3\end{smallmatrix}\right]\right)\\
&+\frac{\epsilon}{6}\left(\left[\begin{smallmatrix}
0&0&0\\
u_2&u_3&u_4\\
0&0&0\\
v_2&v_3&v_4
\end{smallmatrix}\right]
\left[\begin{smallmatrix}
\overline{w}_1\\ 10\im\overline{v}_1\\ -20\overline{u}_1
\end{smallmatrix}\right]
+
\left[\begin{smallmatrix}
0&0&0\\
0&0&0\\
2u_2&2u_3&2u_4\\
-3\im u_2&-3\im u_3&-3\im u_4
\end{smallmatrix}\right]
\left[\begin{smallmatrix}
\overline{w}_2\\ 8\im\overline{v}_2\\ -12\overline{u}_2
\end{smallmatrix}\right]
+
\left[\begin{smallmatrix}
0&0&0\\
0&0&0\\
v_3&v_4&v_5\\
\im v_3&\im v_4&\im v_5
\end{smallmatrix}\right]
\left[\begin{smallmatrix}\overline{v}_1\\ 8\im\overline{u}_1\\ -12\overline{a}\end{smallmatrix}\right]\right)
\\
&+\frac{|a|^2}{36}\left(\epsilon\left[\begin{smallmatrix}
0&0&0&0&0\\
0&0&0&0&0\\
216&144&-150\im&-15&-10\\
1296\im&\tfrac{603}{2}\im&390&\im\tfrac{345}{8}&70\im
\end{smallmatrix}\right]
\left[\begin{smallmatrix}\overline{a}^2b^2\\|a|^2\overline{a}c\\|a|^2b\overline{u}_1\\|a|^2a\overline{v}_1\\a^2{\overline{u}_1}^2\end{smallmatrix}\right]
-
\left[\begin{smallmatrix}
0&0&0&0&0&0\\
30\im c&0&15b&0&\tfrac{9\im}{2}a&0\\
24u_3&-144\im c&-9\im u_2&48b&-u_1&6\im a\\
624\im u_3&-306c&\frac{243}{2}u_2&708\im b&4\im u_1&\tfrac{183}{2}a
\end{smallmatrix}\right]
\left[\begin{smallmatrix}\overline{u}_1\\\overline{u}_2\\\overline{v}_1\\\overline{v}_2\\\overline{w}_1\\\overline{w}_2\end{smallmatrix}\right]\right)\\
&-\frac{1}{36}\left(\left[\begin{smallmatrix}
0&0&0\\
0&-360c&0\\
-36\im u_4&504\im u_3&252\im u_2\\
-234u_4&-1224u_3&-612u_2
\end{smallmatrix}\right]
\left[\begin{smallmatrix}\overline{a}^2a\\\overline{a}^2b\\\overline{a}^2c\end{smallmatrix}\right]
+\left[\begin{smallmatrix}
0&0&0\\
5\im \overline{v}_1&0&0\\
4\im\overline{v}_2&21\im\overline{u}_2&6\im\overline{b}\\
136\overline{v}_2&84\overline{u}_2&24\overline{b}
\end{smallmatrix}\right]
\left[\begin{smallmatrix}a^2\overline{u}_1\\a^2\overline{v}_1\\a^2\overline{w}_1\end{smallmatrix}\right]
\right)\\
&-\frac{1}{36}\left(\left[\begin{smallmatrix}
0&0&0&0\\
20b&0&360\im\overline{a}&0\\
-8\im u_2&-32\im u_1&-288\im\overline{b}&432\im\overline{a}\\
88u_2&52u_1&288\overline{b}&-432\overline{a}
\end{smallmatrix}\right]
\left[\begin{smallmatrix}a{\overline{u}_1}^2\\b{\overline{u}_1}^2\\b^2\overline{u}_1\\b^2\overline{u}_2\end{smallmatrix}\right]
+\left[\begin{smallmatrix}
0&0&0&0\\
0&0&0&0\\
-84\overline{v}_1&288c&-5\overline{v}_1&68c\\
156\im\overline{v}_1&1008\im c&5\im\overline{v}_1&148\im c
\end{smallmatrix}\right]
\left[\begin{smallmatrix}a|b|^2\\\overline{a}|b|^2\\a|u_1|^2\\\overline{a}|u_1|^2\end{smallmatrix}\right]\right)\\
&+
\frac{1}{9}\left[\begin{smallmatrix}0\\0\\
6\im a\overline{b}c\overline{u}_1+4\im\overline{a}bu_1\overline{v}_1-72\overline{a}b\overline{u}_1u_2+18 ab\overline{u}_1\overline{u}_2\\
24a\overline{b}c\overline{u}_1-14\overline{a}bu_1\overline{v}_1-297\im\overline{a}b\overline{u}_1u_2-162\im ab\overline{u}_1\overline{u}_2
\end{smallmatrix}\right],
\end{aligned}
\end{equation}
and our collection is completed by
\begin{align}\label{WR}
\left[\begin{array}{c}W'\\R_1''\end{array}\right]&=
-\frac{\epsilon}{6}\left(4\left[\begin{smallmatrix}\im a&\im b&\im c\\3a&3b&3c\end{smallmatrix}\right]
\left[\begin{smallmatrix}\overline{z}_4\\ 6\im\overline{w}_4\\ -6\overline{v}_4\end{smallmatrix}\right]
+\left[\begin{smallmatrix}-u_1&-u_2&-u_3\\8\im u_1&8\im u_2&8\im u_3\end{smallmatrix}\right]
\left[\begin{smallmatrix}\overline{z}_3\\ 8\im\overline{w}_3\\-12\overline{v}_3\end{smallmatrix}\right]
-\left[\begin{smallmatrix}0&0&0\\v_1&v_2&v_3\end{smallmatrix}\right]
\left[\begin{smallmatrix}\overline{z}_2\\ 10\im\overline{w}_2\\ -20\overline{v}_2\end{smallmatrix}\right]\right)\\
&-\frac{\epsilon}{6}\left(4\left[\begin{smallmatrix}2\im u_2&2\im u_3&2\im u_4\\u_2&u_3&u_4\end{smallmatrix}\right]
\left[\begin{smallmatrix}\overline{w}_3\\ 6\im\overline{v}_3\\ -6\overline{u}_3\end{smallmatrix}\right]
+2\left[\begin{smallmatrix}-v_2&-v_3&-v_4\\3\im v_2&3\im v_3&3\im v_4\end{smallmatrix}\right]
\left[\begin{smallmatrix}\overline{w}_2\\ 8\im\overline{v}_2\\ -12\overline{u}_2\end{smallmatrix}\right]
-\left[\begin{smallmatrix}0&0&0\\w_2&w_3&w_4\end{smallmatrix}\right]
\left[\begin{smallmatrix}\overline{w}_1\\ 10\im\overline{v}_1\\ -20\overline{u}_1\end{smallmatrix}\right]\right)\nonumber\\
&+\frac{\epsilon}{6}\left(4\left[\begin{smallmatrix}-\im v_3&-\im v_4&-\im v_5\\2v_3&2v_4&2v_5\end{smallmatrix}\right]
\left[\begin{smallmatrix}\overline{v}_2\\ 6\im\overline{u}_2\\ -6\overline{b}\end{smallmatrix}\right]
+\left[\begin{smallmatrix}w_3&w_4&w_5\\2\im w_3&2\im w_4&2\im w_5\end{smallmatrix}\right]
\left[\begin{smallmatrix}\overline{v}_1\\ 8\im\overline{u}_1\\ -12\overline{a}\end{smallmatrix}\right]
+\overline{a}^3a\left[\begin{smallmatrix}\tfrac{99}{2}u_3 &252u_2 &4u_1\\414\im u_3& 1404\im u_2& 228\im u_1\end{smallmatrix}\right]
\left[\begin{smallmatrix}a\\b\\c\end{smallmatrix}\right]\right)\nonumber\\
&+\frac{\epsilon}{6} \left(a^3\left[\begin{smallmatrix}16\im\overline{b}\overline{u}_1& 54\im\overline{a}\overline{u}_1&16\im\overline{a}\overline{b}& \tfrac{33}{2}\im\overline{a}^2\\108\overline{b}\overline{u}_1 &\tfrac{819}{2}\overline{a}\overline{u}_1&\tfrac{261}{2}\overline{a}\overline{b}&\tfrac{399}{2}\overline{a}^2\end{smallmatrix}\right]
\left[\begin{smallmatrix}\overline{u}_1\\ \overline{u}_2\\ \overline{v}_1\\ \overline{v}_2\end{smallmatrix}\right]
+|a|^4\left[\begin{smallmatrix}15\im c_1&-57\im u_2&-81b&-\tfrac{125}{24}u_1\\360c& \tfrac{903}{2}u_2& 918\im b &\ 10\im u_1\end{smallmatrix}\right]
\left[\begin{smallmatrix}\overline{b}\\ \overline{u}_1\\ \overline{u}_2\\ \overline{v}_1\end{smallmatrix}\right]\right)\nonumber\\
&+\frac{\epsilon}{6}\left[\begin{smallmatrix}-(152a\overline{u}_1+324\im\overline{a}b)&-(15a\overline{u}_1+\tfrac{224}{3}\im\overline{a}b)& 72
\\1008\overline{a}b+546\im a\overline{u}_1&431\overline{a}b+\tfrac{45}{2}\im a\overline{u}_1& 384\im\end{smallmatrix}\right]
\left[\begin{smallmatrix}|ab|^2\\ |au_1|^2\\ \overline{a}^3b^2u_1\end{smallmatrix}\right]
+\frac{\overline{a}^2}{6}\left[\begin{smallmatrix}-24\im&-90\im& -10\im& -4\im& -126\im\\138& 360& 120& 48& 432\end{smallmatrix}\right]
\left[\begin{smallmatrix}av_4\\ bv_3\\ cv_2\\ u_1u_4\\ u_2u_3\end{smallmatrix}\right]\nonumber\\
&+\frac{a^2}{6}\left[\begin{smallmatrix}-12& -12& 0& -12& -36& -4\\81\im& 66\im& 6\im& 96\im& 168\im& 42\im\end{smallmatrix}\right]
\left[\begin{smallmatrix}\overline{a}\overline{w}_3\\ \overline{b}\overline{w}_2\\ \overline{c}\overline{w}_1\\ \overline{u}_1\overline{v}_3\\ \overline{u}_2\overline{v}_2\\ \overline{u}_3\overline{v}_1\end{smallmatrix}\right]
-\frac{|a|^2}{6}\left[\begin{smallmatrix}36\im&24& \tfrac{3}{2}\im& 24& 54\im& -24&0& \tfrac{2}{3}\im& 18\\
558&198\im& 42& 108\im& 252& 252\im&\tfrac{\im}{4}& 37& 231\im\end{smallmatrix}\right]
\left[\begin{smallmatrix}b\overline{v}_3\\ c\overline{u}_3\\ u_1\overline{w}_2\\ u_2\overline{v}_2\\ u_3\overline{u}_2\\ u_4\overline{b}\\ v_1\overline{w}_1\\ v_2\overline{v}_1\\ v_3\overline{u}_1\end{smallmatrix}\right]\nonumber\\
&-\frac{1}{6}\left(|b|^2\left[\begin{smallmatrix}96\im& -144& -144\im& -6& -40\im& 288\\
528& 432\im& 288& 8\im&\ 480& 1296\im\end{smallmatrix}\right]
+|u_1|^2\left[\begin{smallmatrix}\tfrac{14}{3}\im& -\tfrac{4}{3}&0&-\tfrac{5}{6}&-\tfrac{20}{3}\im& 22\\64& 84\im& 60& 0& 30& 164\im\end{smallmatrix}\right]\right)
\left[\begin{smallmatrix}a\overline{v}_2\\ b\overline{u}_2\\ c\overline{b}\\ u_1\overline{v}_1\\ u_2\overline{u}_1\\ u_3\overline{a}\end{smallmatrix}\right]\nonumber\\
&-\frac{1}{6}\left(a\overline{u}_1\left[\begin{smallmatrix}-\tfrac{135}{4}&24&0&-\tfrac{1}{6}&8\im&-\tfrac{4}{3}\im\\60\im&48\im&90\im&\tfrac{\im}{2}&384&26\end{smallmatrix}\right]
+b\overline{a}\left[\begin{smallmatrix}-\frac{405}{4}\im & 120\im& 36\im & -\im& 72&\tfrac{80}{3}\\\tfrac{675}{2}&0&288&12&504\im &260\im\end{smallmatrix}\right]\right)
\left[\begin{smallmatrix}|a|^6\\|c|^2\\|u_2|^2\\|v_1|^2\\b\overline{u}_3\\v_2\overline{u}_1\end{smallmatrix}\right]\nonumber\\
&-\frac{1}{6}\left(\overline{a}u_1\left[\begin{smallmatrix}4& 16\im& -\tfrac{1}{6}& -5\im&0\\
108\im& 48& 0& 0&0\end{smallmatrix}\right]
+\overline{b}a\left[\begin{smallmatrix}0&-72& \im& 4& 12\im\\0&216\im& 10& -12\im& 276\end{smallmatrix}\right]\right)
\left[\begin{smallmatrix}b\overline{v}_2\\ c\overline{u}_2\\ u_1\overline{w}_1\\ u_2\overline{v}_1\\ u_3\overline{u}_1\end{smallmatrix}\right]\nonumber\\
&-\frac{1}{6}\left[\begin{smallmatrix}\im(108b\overline{u}_2+\tfrac{25}{6}u_1\overline{v}_1)a\overline{u}_2+24(\overline{b}c+2u_2\overline{u}_1)\overline{a}u_2-(10\im a\overline{v}_1-56b\overline{u}_1)b\overline{c}\\
(60\overline{c}\overline{v}_1+504{\overline{u}_2}^2)ab
+(648\im\overline{b}c+40u_1\overline{v}_1+246\im\overline{u}_1u_2)\overline{a}u_2
+10bv_1{\overline{u}_1}^2+72\im b^2\overline{c}\overline{u}_1+30\im\overline{a}cv_1\overline{u}_1+40au_1\overline{u}_2\overline{v}_1\end{smallmatrix}\right].\nonumber
\end{align}

\vspace{\baselineskip}

\subsubsection{Levi-Flat 3-folds} \label{LFIIsec}
Our adaptation \eqref{LNfirstadapt} together with \eqref{dv0_pullback} implies that 
\begin{align}\label{zetaiseta}
\zeta-\eta=0  \text{ on }\mathcal{H}^2.
\end{align}
Using Cartan's Lemma with \eqref{UMCeq} while recalling that $\rho$ is $\mathbb{R}$-valued,
\begin{align}\label{LFIIcoeff}
&\left[\begin{array}{c}\phi_2\\2\im\rho\end{array}\right]=
-\left[\begin{array}{c}\xi\\\phi_1-\overline{\phi}_1+\lambda-\overline{\lambda}\end{array}\right]+
\left[\begin{array}{cc}a&\im b\\\im b&0\end{array}\right]\left[\begin{array}{c}\kappa\\\eta\end{array}\right]
&\text{ for some } a\in C^\infty(\mathcal{H}^2,\mathbb{C})\quad b\in C^\infty(\mathcal{H}^2).
\end{align}
To see how these functions vary on $\mathcal{H}^2$, we invoke the MC equations \eqref{UMCeq} again to differentiate
\begin{align}\label{LFdII}
\dd\left[\begin{array}{c}b\\a\end{array}\right]=-
\left[\begin{array}{cc}\phi_1+\overline{\phi}_1+\lambda+\overline{\lambda}&0\\-2\im\xi&\tfrac{1}{2}(\phi_1+\overline{\phi}_1+\lambda+5\overline{\lambda})\end{array}\right]\left[\begin{array}{c}b\\a\end{array}\right]+
\left[\begin{array}{cc}u_0&0\\u_1&\im u_0+b^2\end{array}\right]
\left[\begin{array}{c}\kappa\\\eta\end{array}\right],\\
\text{for some } u_0\in C^\infty(\mathcal{H}^2),\quad u_1\in C^\infty(\mathcal{H}^2,\mathbb{C}),\nonumber
\end{align}
and we find: if $a$ vanishes identically on each fiber $\mathcal{H}^2\to M$, then so does $b$; $b$ is either zero or non-vanishing on each fiber.

\begin{definition}
The second fundamental form \eqref{II} of a Levi-flat 3-fold $M$ CR embedded in $\Q=SU(2,2)/\mathcal{P}$ is given by
\begin{equation}\label{LFII}
\II=a\kappa\odot\kappa+2\im b\kappa\odot\eta,
\end{equation}
where the coefficients are derived from \eqref{zetaiseta} via \eqref{LFIIcoeff}. The condition that $\II$ is of (sub)maximal rank on a fiber of $\mathcal{H}^2\to M$ is invariant under the action of CR symmetry group $SU(2,2)$ on $\Q$.
\end{definition}

We update the remaining MC equations on $\mathcal{H}^2$,
\begin{equation}\label{LFHtwoMC}
\begin{aligned}
\dd\kappa&=(\lambda+\overline{\lambda})\wedge\kappa,\\
\dd\eta&=-\xi\wedge\kappa+\tfrac{1}{2}(\phi_1+\overline{\phi}_1+3\lambda-\overline{\lambda})\wedge\eta-\tfrac{\im}{2}b\kappa\wedge\eta,\\
\dd\phi_1&=-\phi_1\wedge\overline{\phi}_1+\im\xi\wedge\overline{\eta}+\im\overline{\xi}\wedge\eta-\im b\phi_1\wedge\kappa-\im a\kappa\wedge\overline{\eta}+b\eta\wedge\overline{\eta},\\
\dd\xi&=\tfrac{1}{2}(\phi_1+\overline{\phi}_1+\lambda-3\overline{\lambda})\wedge\xi+\psi\wedge\eta+(\tfrac{\im}{2}b\xi-a\overline{\phi}_1)\wedge\kappa-\im b\overline{\phi}_1\wedge\eta,\\
\dd\lambda&=-\psi\wedge\kappa+\im\overline{a}\kappa\wedge\eta-b\eta\wedge\overline{\eta},\\
\dd\psi&=-(\lambda+\overline{\lambda})\wedge\psi+\im(a\overline{\xi}-\overline{a}\xi)\wedge\kappa-b(\xi\wedge\overline{\eta}+\overline{\xi}\wedge\eta)+b\kappa\wedge(\overline{a}\eta+a\overline{\eta})+\im b^2\eta\wedge\overline{\eta},
\end{aligned}
\end{equation}
and gather more differential identities by applying $\dd^2=0$ to the exterior derivative of \eqref{LFdII}, 
\begin{align}\label{LFdus}
\dd\left[\begin{array}{c}u_0\\u_1\end{array}\right]&=
-\left[\begin{array}{cc}\phi_1+\overline{\phi}_1+2(\lambda+\overline{\lambda})&0\\-3\im\xi&\tfrac{1}{2}(\phi_1+\overline{\phi}_1+3\lambda+7\overline{\lambda})\end{array}\right]
\left[\begin{array}{c}u_0\\u_1\end{array}\right]
+
\left[\begin{array}{cc}-2b&0\\
-3a&2b^2\end{array}\right]
\left[\begin{array}{l}\psi\\\xi\end{array}\right]
\\
&+
\im b\left[\begin{array}{cc}- b& b\\
-\tfrac{1}{2}a&\tfrac{5}{2}a\end{array}\right]
\left[\begin{array}{l}\phi_1\\\overline{\phi}_1\end{array}\right]
+
\im\left[\begin{array}{cc}-2\overline{a}b&2 ab\\
-|a|^2&3 a^2\end{array}\right]
\left[\begin{array}{c}\eta\\\overline{\eta}\end{array}\right]
+
\left[\begin{array}{cc}v_0&0\\
v_1&\im v_0+\tfrac{5}{2}bu_0-\tfrac{\im}{2}b^3\end{array}\right]
\left[\begin{array}{c}\kappa\\\eta\end{array}\right],
\nonumber
\end{align} 
for some $v_0\in C^\infty(\mathcal{H}^2)$, $v_1\in C^\infty(\mathcal{H}^2,\mathbb{C})$. Differentiate \eqref{LFdus} to get
\begin{align}\label{LFdvs}
\dd\left[\begin{array}{c}v_0\\v_1\end{array}\right]&=
-\left[\begin{array}{cc}\phi_1+\overline{\phi}_1+3(\lambda+\overline{\lambda})&0\\-4\im\xi&\tfrac{1}{2}(\phi_1+\overline{\phi}_1+5\lambda+9\overline{\lambda})\end{array}\right]
\left[\begin{array}{c}v_0\\v_1\end{array}\right]
\nonumber\\
&+
\left[\begin{array}{ccc}-6u_0&-4\im\overline{a}b&4\im ab\\
-8u_1&-\tfrac{3\im}{2}b^3-4\im|a|^2+8bu_0&6\im a^2\end{array}\right]
\left[\begin{smallmatrix}\psi\\\xi\\\overline{\xi}\end{smallmatrix}\right]
\nonumber\\
&+
 \left[\begin{array}{cc}b(b^2-3\im u_0)& b(b^2+3\im u_0)\\
\tfrac{1}{2}a(b^2-\im u_0)-\im bu_1&\tfrac{1}{2}a(9b^2+11\im u_0)+3\im bu_1\end{array}\right]
\left[\begin{array}{l}\phi_1\\\overline{\phi}_1\end{array}\right]
\\
&-
\left[\begin{array}{cc}4\overline{a}b^2+5\im\overline{a}u_0+2\im b\overline{u}_1&4ab^2-5\im au_0-2\im bu_1\\
6|a|^2b+3\im\overline{a}u_1+\im a\overline{u}_1&5a^2b-10\im au_1\end{array}\right]
\left[\begin{array}{c}\eta\\\overline{\eta}\end{array}\right]
\nonumber\\
&+
\left[\begin{array}{cc}w_0&0\\
w_1&\im w_0+\tfrac{5}{2}{u_0}^2+3bv_0-\tfrac{11\im}{4}b^2u_0-\tfrac{1}{4}b^4\end{array}\right]
\left[\begin{array}{c}\kappa\\\eta\end{array}\right];
\nonumber
\end{align} 
$w_0\in C^\infty(\mathcal{H}^2)$, $w_1\in C^\infty(\mathcal{H}^2,\mathbb{C})$, and a final derivative provides
\begin{align}\label{LFdws}
\dd\left[\begin{array}{c}w_0\\w_1\end{array}\right]&=
-\left[\begin{array}{cc}\phi_1+\overline{\phi}_1+4(\lambda+\overline{\lambda})&0\\-5\im\xi&\tfrac{1}{2}(\phi_1+\overline{\phi}_1+7\lambda+11\overline{\lambda})\end{array}\right]
\left[\begin{array}{c}w_0\\w_1\end{array}\right]
\nonumber\\
&-
\left[\begin{smallmatrix}12v_0&6\overline{a}b^2+15\im\overline{a}u_0+6\im b\overline{u}_1&6ab^2-15\im au_0-6\im bu_1\\
15v_1&b(b^3+8|a|^2-13v_0)+5\im(3\overline{a}u_1+a\overline{u}_1+\tfrac{9}{4}b^2u_0)-\tfrac{21}{2}{u_0}^2&8a^2b-30\im au_1\end{smallmatrix}\right]
\left[\begin{smallmatrix}\psi\\\xi\\\overline{\xi}\end{smallmatrix}\right]
\nonumber\\
&+
 \left[\begin{array}{cc}\im b(4|a|^2+b^3-4v_0)& -\im b(4|a|^2+b^3-4v_0)\\
\im a(6|a|^2+\tfrac{1}{2}(b^3-v_0))&
-\im a(4|a|^2+6b^3-\tfrac{19}{2}v_0)\end{array}\right]
\left[\begin{array}{l}\phi_1\\\overline{\phi}_1\end{array}\right]
\nonumber\\
&+
 \left[\begin{array}{cc}3u_0(2b^2-\im u_0)& 3u_0(2b^2+\im u_0)\\
\tfrac{3}{2}(b(au_0+bu_1-\im v_1)-\im u_0u_1)&
\tfrac{1}{2}(b(45au_0+15bu_1+7\im v_1)+17\im u_0u_1)\end{array}\right]
\left[\begin{array}{l}\phi_1\\\overline{\phi}_1\end{array}\right]
\\
&-
\left[\begin{smallmatrix}3\overline{a}b(\tfrac{13}{2}u_0-\im b^2)+\overline{u}_1(5b^2+7\im u_0)+9\im\overline{a}v_0+2\im b\overline{v}_1&3ab(\tfrac{13}{2}u_0+\im b^2)+u_1(5b^2-7\im u_0)-9\im av_0-2\im bv_1\\
\tfrac{1}{2}(23|a|^2u_0+b(37\overline{a}u_1+13a\overline{u}_1))+\im(6\overline{a}v_1+a\overline{v}_1+4|u_1|^2)&
3ab(8u_1+\im ab)+\tfrac{11}{2}a^2u_0-15\im av_1-10\im{u_1}^2\end{smallmatrix}\right]
\left[\begin{array}{c}\eta\\\overline{\eta}\end{array}\right]
\nonumber\\
&+
\left[\begin{array}{cc}z_0&0\\
z_1&\im z_0-\tfrac{27\im}{4}b{u_0}^2+8u_0v_0+\tfrac{7}{2}bw_0+\tfrac{\im}{2}b^2(15|a|^2-\tfrac{17}{2}v_0)-\tfrac{19}{8}b^3u_0+\tfrac{\im}{8}b^5\end{array}\right]
\left[\begin{array}{c}\kappa\\\eta\end{array}\right],
\nonumber
\end{align} 
for $z_0\in C^\infty(\mathcal{H}^2)$, $z_1\in C^\infty(\mathcal{H}^2,\mathbb{C})$.

\vspace{\baselineskip}

\section{Classification: Levi-Nondegenerate Case}\label{LNclasssec}

We have seen that if any coefficient of the second fundamental form \eqref{LNII} vanishes identically on $\mathcal{H}^2$, then $\II=0$ everywhere over $M$. The leading coefficient $a\in C^\infty(\mathcal{H}^2,\mathbb{C})$ is particularly indicative, since it is either identically zero or non-vanishing on each fiber of $\mathcal{H}^2\to M$, so the process of adapting frames ``branches" based on the (non)degeneracy of $a$.

If $a=0\Rightarrow\II=0$ everywhere, \eqref{LNH2MCeq} subsides to the Maurer-Cartan equations of $U(2,1)$, which translate via \eqref{threeCRSEinH} to the structure equations \eqref{threeCRSE} of the flat CR 3-sphere over a quotient of $\mathcal{H}^2$ by the central action of $U(1)$ mentioned in Remark \ref{rrecycle}. In the terminology of Definition \ref{equivembed}, $M$ is equivariantly embedded in $\Q$ as an orbit of $U(2,1)\subset SU_\star$, and in this case no further reduction of $\mathcal{H}^2$ is admissible.

Otherwise, we make the generic assumption that $a\neq0$ and arrive at another branching point based on whether $M$ is flat. Locally, all flat 3-folds are equivalent to the 3-sphere, hence the question embeddability is trivial. Rather, we focus on equivariant embeddings by normalizing $\II$ and tracking which algebras emerge as the symmetries of constant-coefficient structure equations. The results are summarized in Theorem \ref{flathomogthm}. For non-flat 3-folds, we reduce $\mathcal{H}^2\to M$ to the Cartan bundle $\mathcal{F}^3\to M$ and carry out the same normalization procedure as in \S\ref{CR3sec}. The resulting structure equations classify embedded, curved 3-folds as detailed in Theorem \ref{curvedtheorem}. Equivariant embeddings are gathered into Theorem \ref{curvedhomogthm}.

\vspace{\baselineskip}

\subsection{Flat 3-folds}\label{LNflatsec}
If $M$ is flat, the coefficients $S,P$ of the curvature tensor appearing in \eqref{threeCRSE} are zero, as are the coefficients \eqref{Bianchi} of the covariant derivative of curvature, the second covariant derivative \eqref{dRQUVW}, and those we've named \eqref{dR'} within the third covariant derivative. Since $a\neq0$, we can use the expressions \eqref{SPR}, \eqref{QUVWR}, \eqref{RRzero}, \eqref{QUUV}, and \eqref{WR} of these coefficients to solve for the jet coordinates of $\II$. When more than one variable is eligible for solution -- e.g., we can solve $S=\tfrac{\epsilon}{6}(a\overline{v}_1+8\im b\overline{u}_1-12c\overline{a})=0$ for either of $c,\overline{v}_1$ -- we default to the higher-order jet:
\begin{align}\label{flatSP}
&\text{solve }S=0\text{ for }\overline{v}_1,
&\text{solve }P=0\text{ for }\overline{v}_2,
&&\text{etc.},
\end{align}
unless the higher-order variable is obtained from another equation. It is also understood that we solve the complex-conjugated equation for the conjugate coordinate. Onward,
\begin{align}\label{solveflat}
&\begin{array}{l}
\text{solve }R=0\text{ for }\Re v_3, \\ 
\text{solve }Q=0\text{ for }\overline{w}_1, \\ 
\text{solve }U=0\text{ for }\overline{w}_2, \\ 
\text{solve }V=0\text{ for }u_4, \\ 
\text{solve }W=0\text{ for }\overline{w}_3,\\
\text{solve }R_1'=0\text{ for }v_4,\\
\text{solve }R_1''=0\text{ for }w_5,
\end{array}
\begin{array}{l}
\text{solve }R_0'=0\text{ and }
R_0''=0\text{ for }\overline{w}_4,w_4,\\
\text{solve }Q'=0\text{ for }\overline{z}_1,\\
\text{solve }U_1'=0\text{ for }\overline{z}_2,\\
\text{solve }U_2'=0\text{ for }v_5,\\
\text{solve }V'=0\text{ for }\overline{z}_3,\\
\text{solve }W'=0\text{ for }\overline{z}_4.
\end{array}
\end{align}

Now we resume the process of adapting frames. The identities \eqref{dII} for $\dd a$ and $\dd b$ tell us that we can restrict to frames whose second fundamental form is diagonalized with leading coefficient 1,
\begin{align*}
\mathcal{H}^3&=\{\underline{\vv}\in\mathcal{H}^2: a(\underline{\vv})=1, b(\underline{\vv})=0\},
\end{align*}
over which we have 
\begin{equation}\label{lambdaxiflat}
\begin{aligned}
\lambda&=-3\im\rho+\overline{u}_1\overline{\eta}+\overline{u}_2\kappa,\\
\xi&=-\im c\overline{\eta}-u_2\eta-u_3\kappa.
\end{aligned}
\end{equation}
With this normalization, we have constrained the coordinates $c_1$ \eqref{P0} and $l$ \eqref{Pl} in the fibers of $\mathcal{H}^2\to M$. In particular, fixing $l$ implies that sections $M\to\mathcal{H}^3$ factor through a unique lift $M\to\hat{M}$. The real part of $c_0$ \eqref{P0} is similarly determined by the reducing to
\begin{align*}
\mathcal{H}^4&=\{\underline{\vv}\in\mathcal{H}^3: \Re u_2(\underline{\vv})=0\},
\end{align*}
as suggested by the equation \eqref{dus} for $\dd u_2$, which reveals that on $T\mathcal{H}^4$,
\begin{equation}\label{psiflat}
\begin{aligned}
\psi&=\left(\tfrac{5}{4}+7|c|^2+\tfrac{5\epsilon}{12}|u_1|^2-3{u_2}^2-\tfrac{\im}{6}u_2(8|u_1|^2+27\epsilon)-c{u_1}^2-\overline{c}{\overline{u}_1}^2+\tfrac{1}{2}u_1u_3+\tfrac{1}{2}\overline{u}_1\overline{u}_3\right)\kappa\\
&+\tfrac{\im}{4}(8\overline{c}\overline{u}_1-5\epsilon u_1+8\im u_1u_2-10\overline{u}_3)\eta
-\tfrac{\im}{4}(8cu_1-5\epsilon\overline{u}_1+8\im \overline{u}_1u_2-10u_3)\overline{\eta}.
\end{aligned}
\end{equation}

The rank of $\II$ will index our final instance of branching in the flat setting. If $c=0$ so that $\text{rank}(\II)=1$, 
\begin{align*}
&\dd c=0\Rightarrow u_1=u_3=0,
&\dd u_1=0\Rightarrow u_2=-\im\frac{\epsilon}{2},
&&\dd u_2=0\Rightarrow \Im v_3=0,
\end{align*}
and we are left with structure equations
\begin{align*}
\dd\kappa&=\im\eta\wedge\overline{\eta},\\
\dd\eta&=\im(\epsilon\kappa-4\rho)\wedge\eta,\\
\dd\rho&=0,
\end{align*}
which describe a central extension of the symmetry algebra of (VIII,C) when $\epsilon=-1$ or (IX,D) when $\epsilon=1$ (see the end of \S\ref{CR3sec}). 

On the other hand, if $\II$ has maximal rank 2 it must be that $c$ is nonvanishing, so the identity \eqref{dII} for $\dd c$ shows that we can reduce to those frames where $c$ takes values in $\mathbb{R}\setminus\{0\}$,
\begin{align*}
\mathcal{H}^5&=\{\underline{\vv}\in\mathcal{H}^4: \Im c(\underline{\vv})=0\},
\end{align*}
over which
\begin{align}\label{rhoflat}
\rho=\frac{1}{96c}(\epsilon(18c+{u_1}^2+{\overline{u}_1}^2)+6(u_1\overline{u}_3+\overline{u}_1u_3)-4c|u_1|^2)\kappa
+\frac{\im}{16c}(3u_1c-u_3)\eta
-\frac{\im}{16c}(3\overline{u}_1c-\overline{u}_3)\overline{\eta}.
\end{align}
Differentiating this and comparing to \eqref{UMCeq} provides
\begin{align}\label{flatcomp}
\epsilon(9cu_1-{u_1}^3)-18c^2\overline{u}_1+2c\overline{u}_1{u_1}^2+54c\overline{u}_3-6\overline{u}_3{u_1}^2=0.
\end{align}
Therefore, the most generic CR embedding in $\Q$ of a flat 3-fold is encoded in the structure equations
\begin{equation}\label{flatgenSE}
\begin{aligned}
\dd\kappa&=\im\eta\wedge\overline{\eta}-\kappa\wedge(u_1\eta+\overline{u}_1\overline{\eta}),\\
\dd\eta&=-\frac{c\overline{u}_1+\overline{u}_3}{4c}\eta\wedge\overline{\eta}-\im\frac{\epsilon(18c+{u_1}^2+{\overline{u}_1}^2)+6(u_1\overline{u}_3+\overline{u}_1u_3)-4c|u_1|^2-48\im cu_2}{24c}\kappa\wedge\eta-\im c\kappa\wedge\overline{\eta},
\end{aligned}
\end{equation}
where $c,u_1,u_3\in C^\infty(\mathcal{H}^5,\mathbb{C})$ and $-\im u_2,\Im v_3\in C^\infty(\mathcal{H}^5,\mathbb{R})$ satisfy \eqref{dII}, \eqref{dus}, and \eqref{dvs} subject to $a=1,b=0$, \eqref{flatSP}, \eqref{solveflat}, \eqref{lambdaxiflat}, \eqref{psiflat}, \eqref{rhoflat}, and \eqref{flatcomp}. The equations \eqref{flatgenSE} remain invariant under the action of the CR symmetry group $SU(3-\delta_\epsilon,1+\delta_\epsilon)$ on $\Q$, so they classify embedded flat 3-folds whose second fundamental form has rank 2. Moreover, when \eqref{flatgenSE} has constant coefficients, the equations classify equivariant embeddings as in Definition \ref{equivembed}. Taking $c,u_1,u_3$ to be constant implies
\begin{align*}
&c=\tfrac{1}{9}{u_1}^2,
&{u_1}^2={\overline{u}_1}^2,
&&u_2=-\tfrac{\im}{6}(3\epsilon+|u_1|^2),
&&u_3=\tfrac{1}{3}{u_1}^3,
&&\Im v_3=0,
\end{align*} 
so $u_1\neq0$ is either real or imaginary, and by $\dd u_2=0$ from \eqref{dus} we see
\begin{align*}
16\epsilon|u_1|^2+9=0.
\end{align*}
The latter only has solutions when $\epsilon=-1$, given by $u_1=\pm\tfrac{3}{4}$ or $u_1=\pm\tfrac{3}{4}\im$. Submitting the CR coframing $\kappa,\eta$ to the 1-adapted transformation
\begin{align*}
\frac{2}{9}\left[\begin{array}{rc}2|u_1|^2&0\\-9\im|u_1|^2&3u_1\end{array}\right]\left[\begin{array}{c}\kappa\\\eta\end{array}\right]
\end{align*}
brings \eqref{flatgenSE} into the form \eqref{FsixSE} for the homogeneous 3-fold ($\text{VI}_3$, E).

\vspace{\baselineskip}

\subsection{Curved 3-folds}\label{LNcurvedsec} 
Suppose the coefficients $S\in C^\infty(\F^3,\mathbb{C})$ of the curvature tensor and $a\in C^\infty(\mathcal{H}^2,\mathbb{C})$ of the second fundamental form of $M$ are non-vanishing and invoke the identity \eqref{dII} for $\dd a$ to reduce to Hermitian frames where $a$ is $\mathbb{R}$-valued:
\begin{align*}
\mathcal{H}^3=\{\underline{\vv}\in\mathcal{H}^2 : a(\underline{\vv})-\overline{a}(\underline{\vv})=0\}.
\end{align*}
On $T\mathcal{H}^3$ we have 
\begin{align*}
\rho=\frac{\im}{6a}(a(\lambda-\overline{\lambda})+(u_2-\overline{u}_2)\kappa+(u_1+2\im\overline{b})\eta-(\overline{u}_1-2\im b)\overline{\eta}),
\end{align*}
and the transformation \eqref{threeCRSEinH} is the infinitesimal version of the identification $\mathcal{H}^3=\mathcal{F}^3$ as bundles of coframes over $M$. Now we pursue the same process of reduction as in \S\ref{CR3sec}. There, $\F^4$ is defined by normalizing $S$, which is achieved in the present setting by solving the equation \eqref{SPR} $S=1$ for $\overline{v}_1$, thereby exhausting one complex degree of freedom in the fibers of $\mathcal{H}^3\supset\mathcal{H}^4=\F^4$ over $M$. Accordingly, $\lambda$ is no longer an independent 1-form, but is determined by \eqref{alphareduction} via \eqref{threeCRSEinH} and \eqref{SPR}, \eqref{QUVWR}. Next we constrain $\overline{v}_2$ in \eqref{SPR} by $P=0$ so that $\mathcal{H}^5\subset\mathcal{H}^4$ coincides with $\F^5$ and $\xi$ satisfies \eqref{betareduction} by virtue of \eqref{threeCRSEinH}, \eqref{SPR}, and \eqref{QUVWR}. Finally, the condition $\Re U=0$ fixes $\Re w_2$ on $\F^6=\mathcal{H}^6\subset\mathcal{H}^5$, where $\psi$ is subject to \eqref{sigmareduction} with coefficients \eqref{SPR}, \eqref{QUVWR}, and \eqref{QUUV}. Our construction proves the following

\begin{thm}\label{curvedtheorem}
Let $M$ be a 3-dimensional, Levi-nondegenerate CR manifold whose curvature tensor is non-vanishing. $M$ is CR embeddable in the 5-dimensional real hyperquadric $\Q=SU(3-\delta_\epsilon,1+\delta_\epsilon)/\mathcal{P}$ if and only if $M$ admits a 1-adapted CR coframing $\kappa,\eta\in\Omega^1(M,\mathbb{C})$ and functions 
\begin{align}\label{thmfns}
a,b,c,u_1,\dots,u_4,v_1,\dots,v_5,w_1,\dots,w_5,z_1,\dots,z_4\in C^\infty(M,\mathbb{C})
\end{align}
satisfying the structure equations \eqref{FsixSE} for \eqref{ABCreduction} given by \eqref{SPR}, \eqref{QUVWR}, where $S=1,P=\Re U=0$, and the differential identities \eqref{dII}, \eqref{dus}, \eqref{dvs}, and \eqref{dws} hold for \eqref{threeCRSEinH} determined by \eqref{alphareduction}, \eqref{betareduction}, and \eqref{sigmareduction} with \eqref{QUUV}. For such $M$, the structure equations \eqref{FsixSE} remain invariant under the action of $SU(3-\delta_\epsilon,1+\delta_\epsilon)$ on $\Q$.
\end{thm}

Let us implement Theorem \ref{curvedtheorem} to treat the examples of non-flat homogeneous 3-folds. In addition to $v_1,v_2$ decided by $S=1$ and $P=0$, respectively, we find most of the functions \eqref{thmfns} by following the prescription \eqref{solveflat}, except that the coefficients of the covariant derivatives of curvature take the values listed in Remark \ref{flatvsnon} instead of zero. A general embedding now depends on the existence of functions $a,b,c,u_1,u_2,u_3,\Im v_3$ which satisfy \eqref{dII} and the identities $\dd^2=0$, and these conditions become algebraic for equivariant embeddings. Namely, if $\II$ is constant, 
\begin{align*}
c&=\tfrac{2}{a}b^2-3\im Ab-Ca,
&&u_1=2\im\overline{b}+2\overline{A}a,\\
u_2&=\tfrac{2\im}{a}|b|^2+\overline{A}b+\im(|A|^2-\tfrac{B}{4})a-\epsilon\tfrac{\im}{4}a^3,
&&u_3=\tfrac{4\im}{a^2}\overline{b}b^2+\tfrac{6A}{a}|b|^2-\tfrac{\im}{2}(Bb+4C\overline{b}-4|A|^2b)-\epsilon\tfrac{\im}{2}ba^2,
\end{align*}
which in turn implies 
\begin{align}\label{Ab}
&0=A\overline{b}+\overline{A}b,
&B=\tfrac{4}{3a^2}|b|^2+\tfrac{2\im}{a}\overline{A}b+\tfrac{8}{3}|A|^2+\epsilon a^2,
&&C=-\tfrac{10}{9}\overline{A}^2+\tfrac{38\im}{9a}\overline{Ab}+\tfrac{2}{3a^2}(2\overline{b}^2+\epsilon),
\end{align}
along with $\Im v_3=0$. Polynomial relationships between the remaining quantities $a,b,A$ 
are clarified by addressing separately the possible values of $\overline{A}=\pm A$.

First consider $A=0$, yielding 
\begin{align*}
&\overline{b}=\frac{2b}{a^4},
&b(3a^8+6\epsilon b^2)=0.
\end{align*}
If $b=0$, $a\in\mathbb{R}\setminus\{0\}$ is free and the resulting structure equations describe the homogeneous models (VIII, K) when $\epsilon=-1$ or (IX, L) when $\epsilon=1$. Note that $\II$ is diagonalized with determinant $-\epsilon\tfrac{2}{3}$. The alternative ($b\neq0$) is only possible when $\epsilon=-1$; namely, $a=\pm\sqrt[4]{2}$ and $b=\pm\sqrt{2}$, which is (IX, L) for $B=\tfrac{\sqrt{2}}{3}$, and in this case $\text{rank}\II=1$. 

For $A\neq0$, \eqref{Ab} shows $\overline{b}=\tfrac{\overline{A}}{A}b$ while $a,b,A$ are governed by three polynomials of degree six or seven,
\begin{align*}
&D_{\overline{\eta}}\lrcorner a^2\dd c,
&D_\kappa\lrcorner a^2A\dd\overline{b},
&&D_\kappa\lrcorner a^2A\dd u_1,
\end{align*}
where the hook $\lrcorner$ indicates contraction with the vector field dual to the subscripted 1-form. Solutions exist only when $\epsilon=-1$, and up to signs they are given by
\begin{align*}
&A=\im\frac{4\sqrt[4]{10}}{\sqrt{5}},
&a=\frac{\sqrt[4]{10}}{\sqrt{5}},
&&b=-\sqrt{10}.
\end{align*}
The homogeneous model is therefore (VI$_{t}$, E) as in \eqref{sixEgeneral} with $\iota=\im$, $t=\tfrac{4\sqrt[4]{10}}{\sqrt{5}}$, and $m=\tfrac{t}{2}$ so that $S=1$.

\vspace{\baselineskip}

\section{Classification: Levi-Flat Case}\label{LFclasssec}

All Levi-flat 3-folds $M$ are locally CR equivalent, so embeddability $M\subset\Q$ is a question of the signature of $\Q$'s Levi form. In this respect we recall that our discussion here applies only to the real hyperquadric whose CR symmetry group is $SU(2,2)$ ($\epsilon=-1$ in the notation of \S\ref{movingframesec}), hence the embeddings of interest will be equivariant for some action of this Lie group in the sense of Definition \ref{equivembed}. Our list of homogeneous models in \S\ref{CR3sec} omitted the Bianchi algebras that serve as infinitesimal CR symmetries of Levi-flat $M$; let us record here two models with symmetry of Bianchi type V. For an appropriate choice of bases, the extension of $\mathbb{R}^2$ by $\left[\begin{smallmatrix}\sqrt[3]{3}&0\\0&-\sqrt[3]{3}\end{smallmatrix}\right]$ has structure equations
\begin{equation}\label{Vbzero}
\begin{aligned}
\dd\kappa&=\sqrt[3]{3}(\eta+\overline{\eta})\wedge\kappa,\\
\dd\eta&=\sqrt[3]{3}\eta\wedge\overline{\eta}+\tfrac{\im}{15\sqrt[3]{3}}(4\eta-\overline{\eta})\wedge\kappa,
\end{aligned}
\end{equation}
and the extension by $\left[\begin{smallmatrix}1&0\\0&2\end{smallmatrix}\right]$ is described by
\begin{equation}\label{Vbnonzero}
\begin{aligned}
\dd\kappa&=\tfrac{\im}{2}(\eta-\overline{\eta})\wedge\kappa,\\
\dd\eta&=\im\eta\wedge\overline{\eta}.
\end{aligned}
\end{equation}

As in the Levi-nondegenerate case, the process of reducing $\mathcal{H}^2\to M$ to these and other homogeneous models branches based on the rank of the second fundamental form $\II$. This section constitutes the proof of Theorem \ref{LFhomogthm}.

\vspace{\baselineskip}

\subsection{$\II=0$: Maximal Symmetry}\label{LFmaxsec} 

Let $\underline{\vv}=(\vv_0,\vv_1,\vv_2,\vv_3)$ be a Hermitian frame \eqref{hvivj} whose vectors we rearrange,
\begin{align}\label{LFUtrans}
&\left[\begin{array}{cccc}\n_1&\n_2&\n_3&\n_4\end{array}\right]=
\left[\begin{array}{cccc}\vv_0&\vv_1&\vv_2&\vv_3\end{array}\right]U,
&U=\frac{1}{\sqrt{2}}\left[\begin{array}{rrrr}0&\sqrt{2}&0&0\\1&0&0&-1\\1&0&0&1\\0&0&\sqrt{2}&0\end{array}\right].
\end{align}
The transformation $U$ is ``unitary" in the sense that $U^{-1}=\overline{U}^t$, and even though $\overline{U}^t\h U\neq\h$, we refer to the symmetry groups of both forms as $SU(2,2)$. Let $\mathcal{R}_\mathbb{C}\subset SL_2\mathbb{C}$ be the 10-dimensional parabolic subgroup that stabilizes the partial flag 
\begin{align}\label{oscflag}
\langle \n_1\rangle_\mathbb{C}\subset\langle \n_1,\n_2,\n_3\rangle_\mathbb{C}\subset\mathbb{C}^{4},
\end{align}
and name $\mathcal{R}=\mathcal{R}_{\mathbb{C}}\cap SU(2,2)$ with Lie algebra $\mathfrak{r}\subset\mathfrak{su}(2,2)$. 

The new basis \eqref{LFUtrans} consists entirely of $\h$-null vectors, and for $\underline{\vv}\in\mathcal{H}^2$ they are adapted to $M$ -- cf. \eqref{LFfirstadapt} -- in that $\n_1,\n_2,\n_3$ span $T_{\vv_0}\hat{M}$ with $\n_1$ descending to the CR bundle of $M$. If the complex curve tangent to the CR bundle of $M$ is a complex line in $\Q$ (see \cite[Example 1.5]{BryanthololorentzCR}), then the osculating flag \eqref{oscflag} is constant along $M$ and $\hat{M}$ itself is contained in the fixed subspace $\langle \n_1,\n_2,\n_3\rangle_\mathbb{C}$. Hence, one expects the extrinsic CR symmetries of such $M$ given by the action of $SU(2,2)$ on $\Q$ to lie in $\mathcal{R}$. In general, $\mu|_{\mathcal{H}^2}$ is \eqref{UMCform} subject to \eqref{trace-free}, \eqref{zetaiseta}, and \eqref{LFIIcoeff}, so transforming according to \eqref{LFUtrans} we get
\begin{align}\label{LFMCtrans}
U^{-1}\mu|_{\mathcal{H}^2}U=
\left[\begin{array}{cccc}
-\tfrac{1}{2}(\phi_1+\overline{\phi}_1+\lambda-\overline{\lambda})&\sqrt{2}\eta&\sqrt{2}\xi-\tfrac{1}{\sqrt{2}}(a\kappa+\im b\eta)&\phi_1-\overline{\phi}_1-\tfrac{\im}{2}b\kappa\\
-\tfrac{1}{\sqrt{2}}(b\overline{\eta}+\im\overline{a}\kappa)&\lambda&\psi&\sqrt{2}\im\overline{\xi}-\tfrac{1}{\sqrt{2}}(b\overline{\eta}+\im\overline{a}\kappa)\\
0&\kappa&-\overline{\lambda}&-\im\sqrt{2}\overline{\eta}\\
-\tfrac{\im}{2}b\kappa&0&-\tfrac{1}{\sqrt{2}}(a\kappa+\im b\eta)&\tfrac{1}{2}(\phi_1+\overline{\phi}_1-\lambda+\overline{\lambda})
\end{array}\right].
\end{align}
In particular, if $a=b=0$, \eqref{LFMCtrans} takes values $\mathfrak{r}$, and \eqref{LFHtwoMC} are the Maurer-Cartan equations of $\mathcal{R}$. Thus, if $\II=0$, $\mathcal{H}^2$ is locally $\mathcal{R}\subset SU(2,2)$.

\vspace{\baselineskip}

\subsection{Rank$(\II)=1$} The second fundamental form \eqref{LFII} will have rank one only if $b=0$ in \eqref{LFIIcoeff}, which in turn requires $u_0=0$ \eqref{LFdII}, $v_0=0$ \eqref{LFdus}, $w_0=0$ \eqref{LFdvs}, and $z_0=0$ \eqref{LFdws}, but it must be that $a$ is nonvanishing. The identity \eqref{LFdII} for $\dd a$ therefore implies that we can reduce to
\begin{align*}
\mathcal{H}^3=\{\underline{\vv}\in\mathcal{H}^2: a(\underline{\vv})=1\}, 
\end{align*}
and on $\mathcal{H}^3$ we have
\begin{align*}
\lambda=-\tfrac{1}{6}(\phi_1+\overline{\phi}_1)-\tfrac{1}{12}(u_1-5\overline{u}_1)\kappa.
\end{align*}
Looking to $\dd u_1$\eqref{LFdus} and $\dd v_1$ \eqref{LFdvs}, we see the opportunity to reduce further,
\begin{align*}
\mathcal{H}^4=\{\underline{\vv}\in\mathcal{H}^3: \Re u_1(\underline{\vv})=0=v_1(\underline{\vv})\}. 
\end{align*}
After renaming
\begin{align*}
u=\Im u_1\in C^\infty(\mathcal{H}^4),
\end{align*}
we can write the newly imposed constraints
\begin{align*}
\psi&=\tfrac{1}{3}(u^2\kappa-2\im\eta+2\im\overline{\eta}),\\
\xi&=\tfrac{1}{30}(40u^3+6\im w_1-9\im\overline{w}_1)\kappa+\tfrac{\im}{15}(11u\eta+u\overline{\eta}),
\end{align*}
as well as the updated identity
\begin{align*}
\dd w_1=(6\im+w_1)\phi_1+(w_1-4\im)\overline{\phi}_1+(z_1-\tfrac{80}{3}u^4-8\im uw_1+3\im u\overline{w}_1)\kappa+\tfrac{\im}{3}u^2(16\eta+38\overline{\eta}).
\end{align*}
The latter suggests a final reduction 
\begin{align*}
\mathcal{H}^5=\{\underline{\vv}\in\mathcal{H}^4: w_1(\underline{\vv})=0\}, 
\end{align*}
over which
\begin{align*}
\phi_1&=-\tfrac{\im}{30}(80u^4-9z_1+6\overline{z}_1)\kappa-\tfrac{1}{15}u^2(62\eta+73\overline{\eta}),\\
\dd u&=-\tfrac{1}{3}uz\kappa+(1-3u^3)\eta+(1-3u^3)\overline{\eta};
&&z=\Im z_1.
\end{align*}
What's left of \eqref{LFHtwoMC} reads
\begin{equation}\label{LFrankone}
\begin{aligned}
\dd\kappa&=3u^2(\eta+\overline{\eta})\wedge\kappa,\\
\dd\eta&=3u^2\eta\wedge\overline{\eta}
-\tfrac{1}{15}(4\im u+5z)\kappa\wedge\eta+\tfrac{\im}{15}u\kappa\wedge\overline{\eta},
\end{aligned}
\end{equation}
for $u,z\in C^\infty(\mathcal{H}^5)$. The equations \eqref{LFrankone}, along with differential identities necessitated by $\dd^2=0$, classify embedded Levi-flat 3-folds $M\subset\Q$ with rank$(\II)=1$ up to the action of $SU(2,2)$ on $\Q$, including the unique equivariant embedding when $\dd u=0\Rightarrow u=3^{-\tfrac{1}{3}}$, $z=0$, which is exactly the model \eqref{Vbzero}.

\vspace{\baselineskip}

\subsection{Rank$(\II)=2$} The second fundamental form \eqref{LFII} has full rank if and only if $b$ is nonvanishing; in this case the identities \eqref{LFdII}, \eqref{LFdus} for $\dd b, \dd a, \dd u_0$ show that we can reduce to
\begin{align*}
\mathcal{H}^3=\left\{\underline{\vv}\in\mathcal{H}^2: 
\begin{array}{c}
b(\underline{\vv})=\pm1\\
a(\underline{\vv})=0\\
u_0(\underline{\vv})=0
\end{array}
\right\}, 
\end{align*}
over which we have relations
\begin{align*}
\phi_1&=-\frac{1}{2}(\lambda+\overline{\lambda})+\frac{\im}{b}\psi-\frac{v_0}{2b^2}\kappa,\\
\xi&=\frac{\im}{2}b\eta+\frac{\im u_1}{2b}\kappa,
\end{align*}
as well as the identity
\begin{align*}
\dd u_1=-u_1(\lambda+3\overline{\lambda})+(v_1+\im bu_1)\kappa+\im(v_0+\tfrac{1}{2}b^3)\eta.
\end{align*}
Thus we are obliged to consider branching based on the possible values of $u_1$. 

If $u_1=0$ identically, then  $v_0=-\tfrac{1}{2}b^3$, $v_1=0$, and the remaining coefficients in \eqref{LFdvs}, \eqref{LFdws} are zero, leaving 
\begin{equation}\label{LFranktwoMC}
\begin{aligned}
\dd\kappa&=(\lambda+\overline{\lambda})\wedge\kappa,\\
\dd\eta&=(\lambda-\overline{\lambda})\wedge\eta,\\
\dd\lambda&=-\psi\wedge\kappa-b\eta\wedge\overline{\eta},\\
\dd\psi&=\psi\wedge(\lambda+\overline{\lambda}).
\end{aligned}
\end{equation} 
Note that \eqref{LFranktwoMC} are the Maurer-Cartan equations of the MC forms
\begin{align*}
&\left[\begin{array}{cc}\tfrac{1}{2}(\lambda+\overline{\lambda})&\psi\\\kappa&-\tfrac{1}{2}(\lambda+\overline{\lambda})\end{array}\right]\in\Omega^1(\mathcal{H}^3,\mathfrak{sl}_2\mathbb{R}),
&\left[\begin{array}{cc}\tfrac{1}{2}(\lambda-\overline{\lambda})&-\sqrt{b}\overline{\eta}\\\sqrt{b}\eta&-\tfrac{1}{2}(\lambda-\overline{\lambda})\end{array}\right]\in\Omega^1(\mathcal{H}^3,\mathfrak{su}(p,q)),
\end{align*}
where $(p,q)=(2,0)$ if $b=1$ and $(p,q)=(1,1)$ if $b=-1$. 

Otherwise, in any neighborhood where $u_1\neq0$ we can normalize it to define
\begin{align*}
\mathcal{H}^4=\left\{\underline{\vv}\in\mathcal{H}^3: 
\begin{array}{c}
u_1(\underline{\vv})=1\\
\Re v_1(\underline{\vv})=0
\end{array}
\right\}, 
\end{align*}
so that on $\mathcal{H}^4$,
\begin{align*}
\lambda&=-\tfrac{\im}{2}(b+v)\kappa-\tfrac{\im}{16}(b^3+2v_0)\eta-\tfrac{3\im}{16}(b^3+2v_0)\overline{\eta} &(v=\Im v_1),\\
\psi&=\tfrac{1}{16}(vb^3 + b^4+ 2bv_0 + 2v_0v + 2\im w_0 )\eta
+\tfrac{1}{16}(vb^3 + b^4+ 2bv_0 + 2v_0v - 2\im w_0 )\overline{\eta}\\
&+\frac{1}{16b}(4v^2b + 4vb^2 + 3b^3 + 2b\overline{w}_1 + 2bw_1 - 16v_0)\kappa,
\end{align*}
and the structure equations are 
\begin{equation}\label{LFr2u1}
\begin{aligned}
\dd\kappa&=\tfrac{\im}{8}(b^3+2v_0)(\eta-\overline{\eta})\wedge\kappa,\\
\dd\eta&=\tfrac{\im}{4}(b^3+2v_0)\eta\wedge\overline{\eta}-\im(b+v)\kappa\wedge\eta,
\end{aligned}
\end{equation}
with identities
\begin{equation}\label{LFr2u1id}
\begin{aligned}
\dd v_0&=w_0\kappa-\tfrac{\im}{8}(b^6+4b^3v_0+4{v_0}^2+16b)\eta+\tfrac{\im}{8}(b^6+4b^3v_0+4{v_0}^2+16b)\overline{\eta},\\
\dd v&=\tfrac{\im}{2}(\overline{w}_1-w_1)\kappa-\tfrac{\im}{8}(3vb^3+3b^4+6v_0v+4\im w_0+6bv_0)\eta
+\tfrac{\im}{8}(3vb^3+3b^4+6v_0v-4\im w_0+6bv_0)\overline{\eta}.
\end{aligned}
\end{equation}

We conclude that embedded Levi-flat 3-folds $M\subset\Q$ with rank$(\II)=2$ are classified up to the action of $SU(2,2)$ on $\Q$ by the Maurer-Cartan equations \eqref{LFranktwoMC} of $SL_2\mathbb{R}\times SU(p,q)$ -- $(p,q)=(2,0)$ or $(1,1)$ -- when $u_1=0$, or \eqref{LFr2u1} with \eqref{LFr2u1id} when $u_1\neq0$. The latter evince a homogeneous action in $SU(2,2)$ when
\begin{align*}
&w_0=0,
&v=1,
&&v_0=-\tfrac{3}{2}\text{ or }\tfrac{5}{2},
&&b=-1,
\end{align*}
in which case they coincide (up to sign) with \eqref{Vbnonzero}.

\vspace{2\baselineskip}

\bibliographystyle{amsalpha}

\bibliography{References}

\newcommand{\etalchar}[1]{$^{#1}$}
\providecommand{\bysame}{\leavevmode\hbox to3em{\hrulefill}\thinspace}
\providecommand{\MR}{\relax\ifhmode\unskip\space\fi MR }
\providecommand{\MRhref}[2]{%
  \href{http://www.ams.org/mathscinet-getitem?mr=#1}{#2}
}
\providecommand{\href}[2]{#2}
\begin{thebibliography}{SKM{\etalchar{+}}03}

\bibitem[BGG03]{BGG}
Robert Bryant, Phillip Griffiths, and Daniel Grossman, \emph{Exterior
  differential systems and {E}uler-{L}agrange partial differential equations},
  Chicago Lectures in Mathematics, University of Chicago Press, Chicago, IL,
  2003. \MR{1985469}

\bibitem[Bry82]{BryanthololorentzCR}
Robert~L. Bryant, \emph{Holomorphic curves in {L}orentzian {CR}-manifolds},
  Transactions of the American Mathematical Society \textbf{272} (1982), no.~1,
  203--221.

\bibitem[{Bry}04]{BryantCR3}
R.~L. {Bryant}, \emph{{Real hypersurfaces in unimodular complex surfaces}},
  ArXiv Mathematics e-prints (2004).

\bibitem[Car32]{CartanCR}
\'Elie Cartan, \emph{Sur la g\'eom\'etrie pseudo-conforme des hypersurfaces de
  l'espace de deux variables complexes {II}}, Ann. Scuola Norm. Sup. Pisa Cl.
  Sci. (2) \textbf{1} (1932), no.~4, 333--354. \MR{1556687}

\bibitem[CG19]{CGbook}
S.~N. {Curry} and A.~R. {Gover}, \emph{{CR} embedded submanifolds of {CR}
  manifolds}, Memoirs of the American Mathematical Society, vol. 1241,
  Providence, RI: American Mathematical Society, 2019.

\bibitem[CM74]{ChernMoser}
S.~S. Chern and J.~K. Moser, \emph{Real hypersurfaces in complex manifolds},
  Acta Math. \textbf{133} (1974), 219--271. \MR{0425155}

\bibitem[{\v{C}}S09]{CapSlovak}
Andreas {\v{C}}ap and Jan Slov\'ak, \emph{Parabolic geometries. {I}},
  Mathematical Surveys and Monographs, vol. 154, American Mathematical Society,
  Providence, RI, 2009, Background and general theory. \MR{2532439}

\bibitem[Gar89]{Gardner}
Robert~B. Gardner, \emph{The method of equivalence and its applications},
  CBMS-NSF Regional Conference Series in Applied Mathematics, vol.~58, Society
  for Industrial and Applied Mathematics (SIAM), Philadelphia, PA, 1989.
  \MR{1062197}

\bibitem[Gri74]{GriffithsMF}
P.~Griffiths, \emph{On {C}artan's method of {L}ie groups and moving frames as
  applied to uniqueness and existence questions in differential geometry}, Duke
  Math. J. \textbf{41} (1974), 775--814. \MR{0410607}

\bibitem[GS62]{GoldbergSachs}
J.~N. Goldberg and R.~K. Sachs, \emph{A theorem on {P}etrov types}, Acta Phys.
  Polon. \textbf{22} (1962), no.~suppl., 13--23. \MR{0156679}

\bibitem[IL16]{CFB}
Thomas~A. Ivey and Joseph~M. Landsberg, \emph{Cartan for beginners}, Graduate
  Studies in Mathematics, vol. 175, American Mathematical Society, Providence,
  RI, 2016, Differential geometry via moving frames and exterior differential
  systems, Second edition [of MR2003610]. \MR{3586335}

\bibitem[Jac90]{Jacobowitz}
Howard Jacobowitz, \emph{An introduction to {CR} structures}, Mathematical
  Surveys and Monographs, vol.~32, American Mathematical Society, Providence,
  RI, 1990. \MR{1067341}

\bibitem[Ker63]{Kerr}
Roy~P. Kerr, \emph{Gravitational field of a spinning mass as an example of
  algebraically special metrics}, Phys. Rev. Lett. \textbf{11} (1963),
  237--238. \MR{0156674}

\bibitem[NT88]{NurTaf}
Pawe\l{} Nurowski and Jacek Tafel, \emph{Symmetries of {C}auchy-{R}iemann
  spaces}, Lett. Math. Phys. \textbf{15} (1988), no.~1, 31--38. \MR{929784}

\bibitem[NT02]{NurTraut}
Pawe\l{} Nurowski and Andrzej Trautman, \emph{Robinson manifolds as the
  {L}orentzian analogs of {H}ermite manifolds}, Differential Geom. Appl.
  \textbf{17} (2002), no.~2-3, 175--195, 8th International Conference on
  Differential Geometry and its Applications (Opava, 2001). \MR{1925764}

\bibitem[O'N95]{Oneill}
Barrett O'Neill, \emph{The geometry of {K}err black holes}, A K Peters, Ltd.,
  Wellesley, MA, 1995. \MR{1328643}

\bibitem[Pen67]{Pen67}
R.~Penrose, \emph{Twistor algebra}, J. Mathematical Phys. \textbf{8} (1967),
  345--366. \MR{0216828}

\bibitem[PR88]{PenBook2}
Roger Penrose and Wolfgang Rindler, \emph{Spinors and space-time. {V}ol. 2},
  second ed., Cambridge Monographs on Mathematical Physics, Cambridge
  University Press, Cambridge, 1988, Spinor and twistor methods in space-time
  geometry. \MR{944085}

\bibitem[Rob61]{robnullEM}
Ivor Robinson, \emph{Null electromagnetic fields}, Journal of Mathematical
  Physics \textbf{2} (1961), no.~3, 290--291.

\bibitem[RT62]{RTGRwaves}
Ivor Robinson and Andrzej Trautman, \emph{Some spherical gravitational waves in
  general relativity}, Proceedings of the Royal Society of London A:
  Mathematical, Physical and Engineering Sciences \textbf{265} (1962),
  no.~1323, 463--473.

\bibitem[RT83]{RTflows}
\bysame, \emph{Conformal geometry of flows in n dimensions}, Journal of
  Mathematical Physics \textbf{24} (1983), no.~6, 1425--1429.

\bibitem[Sac61]{Sachs}
R.~Sachs, \emph{Gravitational waves in general relativity. {VI}. {T}he outgoing
  radiation condition}, Proc. Roy. Soc. Ser. A \textbf{264} (1961), 309--338.
  \MR{0156678}

\bibitem[SKM{\etalchar{+}}03]{exactsolutionsbook}
Hans Stephani, Dietrich Kramer, Malcolm MacCallum, Cornelius Hoenselaers, and
  Eduard Herlt, \emph{Exact solutions of einstein's field equations}, 2 ed.,
  Cambridge Monographs on Mathematical Physics, Cambridge University Press,
  2003.

\bibitem[Taf85]{Tafel}
Jacek Tafel, \emph{On the {R}obinson theorem and shearfree geodesic null
  congruences}, Lett. Math. Phys. \textbf{10} (1985), no.~1, 33--39.
  \MR{796997}

\bibitem[Tan62]{Tan62}
Noboru Tanaka, \emph{On the pseudo-conformal geometry of hypersurfaces of the
  space of {$n$}\ complex variables}, J. Math. Soc. Japan \textbf{14} (1962),
  397--429. \MR{0145555}

\bibitem[Tan79]{Tan79}
\bysame, \emph{On the equivalence problems associated with simple graded {L}ie
  algebras}, Hokkaido Math. J. \textbf{8} (1979), no.~1, 23--84. \MR{533089}

\bibitem[Tra58]{Traut58}
A.~Trautman, \emph{Radiation and boundary conditions in the theory of
  gravitation}, Bull. Acad. Polon. Sci. S\'er. Sci. Math. Astr. Phys.
  \textbf{6} (1958), 407--412. \MR{0097266}

\bibitem[Tra17]{TrautHistory}
Andrzej Trautman, \emph{Gravitational waves}, Journal of Physics: Conference
  Series \textbf{873} (2017), no.~1, 012012.

\bibitem[WW90]{WardWells}
R.~S. Ward and Raymond~O. Wells, Jr., \emph{Twistor geometry and field theory},
  Cambridge Monographs on Mathematical Physics, Cambridge University Press,
  Cambridge, 1990. \MR{1054377}

\end{thebibliography}

\end{document}